\documentclass[twoside,english,3p]{elsarticle}
\usepackage[T1]{fontenc}
\usepackage{geometry}
\usepackage{amssymb}
\usepackage{amsmath}

\usepackage{amsthm}
\usepackage{stmaryrd}
\usepackage{graphicx}
\usepackage{booktabs}
\usepackage{multirow}
\usepackage{makecell}
\usepackage{xurl}
\usepackage{esint}
\usepackage{algorithm}
\usepackage{algpseudocode}
\usepackage{rotating}
\usepackage{adjustbox}
\usepackage{nomencl}
\usepackage{multicol}
\makenomenclature
\makeatletter
\theoremstyle{plain}

\theoremstyle{boldremark} 

\ifx\proof\undefined

\providecommand{\proofname}{Proof}
\fi

\renewenvironment{proof}[1][\proofname]{%
	\par\pushQED{\qed}\normalfont%
	\topsep6\p@\@plus6\p@\relax
	\trivlist\item[\hskip\labelsep\bfseries#1\@addpunct{.}]%
	\ignorespaces
}{%
	\popQED\endtrivlist\@endpefalse
}

\journal{Elsevier}

\usepackage{hyperref}
\hypersetup{colorlinks = true, allcolors = blue}

\usepackage[nameinlink]{cleveref}

\Crefname{figure}{Fig.}{Figs.}
\Crefformat{equation}{Eq.~#2(#1)#3}
\Crefformat{section}{Section~#2#1#3}
\AtBeginDocument{%
	\let\citet\cite
}

\usepackage[labelfont=bf]{caption}
\@ifundefined{showcaptionsetup}{}{%
	\PassOptionsToPackage{caption=false}{subfig}}
\usepackage{subfig}
\makeatother

\usepackage{babel}
\providecommand{\remarkname}{Remark}
\providecommand{\theoremname}{Theorem}

\begin{document}
	
	\begin{frontmatter}{}
		
		\title{Pretrain Finite Element Method: A Pretraining and Warm-start Framework for PDEs via Physics-Informed Neural Operators}

\author[rvt,rvt3]{Yizheng Wang}

\ead{wang-yz19@tsinghua.org.cn}

\author[rvt2]{Zhongkai Hao}

\author[rvt6]{Mohammad Sadegh Eshaghi}

\author[rvt3]{Cosmin Anitescu}

\author[rvt6]{Xiaoying Zhuang}

\author[rvt3]{Timon Rabczuk}

\author[rvt]{Yinghua Liu\corref{cor1}}

\ead{yhliu@mail.tsinghua.edu.cn}
\cortext[cor1]{Corresponding author}
\address[rvt]{Department of Engineering Mechanics, Tsinghua University, Beijing 100084, China}

\address[rvt2]{Department of Electronic Engineering, Tsinghua University, Beijing 100084, China}

\address[rvt3]{Institute of Structural Mechanics, Bauhaus-Universit\"{a}t Weimar, Marienstr. 15, D-99423 Weimar, Germany}

\address[rvt6]{ Institute of Photonics, Department of Mathematics and Physics, Leibniz University Hannover, Germany}

\begin{abstract}

We propose a Pretrained Finite Element Method (PFEM), a physics-driven framework that bridges the efficiency of neural operator learning with the accuracy and robustness of classical finite element methods (FEM). PFEM consists of a physics-informed pretraining stage and an optional warm-start stage.
In the pretraining stage, a neural operator based on the Transolver architecture is trained solely from governing partial differential equations, without relying on labeled solution data. The model operates directly on unstructured point clouds, jointly encoding geometric information, material properties, and boundary conditions, and produces physically consistent initial solutions with extremely high computational efficiency. PDE constraints are enforced through explicit finite element–based differentiation, avoiding the overhead associated with automatic differentiation.
In the warm-start stage, the pretrained prediction is used as an initial guess for conventional FEM solvers, preserving their accuracy, convergence guarantees, and extrapolation capability while substantially reducing the number of iterations required to reach a prescribed tolerance.
PFEM is validated on a broad range of benchmark problems, including linear elasticity and nonlinear hyperelasticity with complex geometries, heterogeneous materials, and arbitrary boundary conditions. Numerical results demonstrate strong generalization in the pretraining stage with relative errors on the order of 1\%, and speedups of up to one order of magnitude in the warm-start stage compared to FEM with zero initial guesses.

\end{abstract}

\printnomenclature

		
\end{frontmatter}{}

\section{Introduction}
Numerous physical phenomena are modeled by partial differential equations (PDEs) \citet{loss_is_minimum_potential_energy}, and solving PDEs is central to understanding and predicting physical behavior. Consequently, PDE solvers constitute one of the most fundamental and core tasks in computational mathematics and physics \citet{PINN_review}. Traditional numerical methods often suffer from an intrinsic trade-off between accuracy and efficiency: achieving higher accuracy typically requires substantially increased computational cost \citet{kahana2023geometry}. Moreover, once the material distribution, geometric configuration, or boundary conditions change, conventional PDE solvers must be re-run from scratch \citet{wang2021learning}. Among classical numerical approaches, the finite element method (FEM) is widely regarded for its high accuracy and numerical robustness. However, for complex problems, repeatedly recomputing high-fidelity solutions can be extremely expensive \citet{belytschko2013nonlinear}.

Recently, AI for PDEs, one of the key directions in AI for Science, has emerged as a class of methodologies that leverage deep learning to solve PDEs. Broadly speaking, AI for PDEs can be categorized into three major paradigms: physics-informed neural networks (PINNs) \citet{PINN_original_paper,loss_is_minimum_potential_energy}, operator learning \citet{DeepOnet}, and physics-informed neural operators (PINO) \citet{li2024physics}. In this work, we briefly review these three approaches.

The first paradigm, PINNs, employs neural networks as universal function approximators to directly solve PDEs. However, PINNs are typically restricted to solving a single PDE instance. The network must be retrained whenever the initial or boundary conditions change \citet{yang2023context}. Although certain transfer learning strategies have been explored for PINNs \citet{wang2025transfer}, they remain challenging when the target solution undergoes significant variation. In summary, when faced with even moderately different problems, PINNs encounter limitations similar to those of traditional PDE solvers, namely the need for retraining or re-solving \citet{goswami2022physics}. While PINNs are theoretically appealing, we believe that the full potential of neural networks can only be realized by addressing families of PDEs rather than individual problem instances \citet{yang2023context}.

This motivation leads to the second paradigm: operator learning. Representative methods include DeepONet \citet{DeepOnet} and the Fourier Neural Operator (FNO) \citet{li2020fourier}. Early operator learning methods were purely data-driven and are therefore well suited for large-scale problems with abundant data. Unlike PINNs, which target a specific PDE instance, operator learning aims to learn mappings between function spaces. As a result, once trained, neural operators can rapidly predict solutions under varying boundary conditions, geometries, and material properties \citet{kovachki2023neural}. However, the primary drawback of operator learning lies in its heavy reliance on large datasets, which are often generated using high-fidelity numerical solvers or experiments. This requirement poses challenges for problems where data are scarce or expensive to obtain.

Several studies have focused on enhancing the geometric generalization capability of neural operators. For example, He et al.\ employed PointNet as a generalized feature extractor and used it as the branch network input of DeepONet \citet{he2024geom}. Zeng et al.\ integrated PointNet into the kernel integral of neural operators, following a design philosophy similar to FNO \citet{zeng2025point}. These approaches can be viewed as extensions of DeepONet and FNO. In recent years, transformer-based operator learning frameworks have demonstrated remarkable performance and often outperform DeepONet and FNO, particularly in terms of geometric generalization \citet{bi2023accurate}. Representative examples include GNOT proposed by Hao et al.\ \citet{hao2023gnot} and Transolver introduced by Wu et al.\ \citet{wu2024transolver}. While GNOT directly applies transformers to input points, it suffers from prohibitively high computational complexity. Transolver significantly alleviates this issue by introducing physics attention, thereby greatly reducing computational cost and exhibiting strong potential for practical applications \cite{nabian2025automotive}.

Nevertheless, the aforementioned neural operator methods are typically trained in a fully data-driven manner. To mitigate data dependence, the third paradigm, the physics-informed neural operator (PINO), was introduced \citet{li2024physics,eshaghi2025variational,wang2021learning}. PINO integrates governing physical equations into operator learning, substantially reducing data requirements and, in some cases, enabling training without any labeled data. Given the strong potential of transformer-based operators, particularly the efficiency of Transolver, we focus on a Transolver-based PINO framework. Notably, Transolver requires only point-wise inputs rather than structured grids, making it naturally compatible with FEM and well suited for predicting nodal displacements. To the best of the authors' knowledge, this work represents the first Transolver-based PINO formulation.

Despite the advantages of PINO in terms of rapid prediction and generalization, its accuracy still lags behind that of traditional PDE solvers, especially the FEM, which remains the industrial standard \citet{bathe2006finite,belytschko2013nonlinear,finite_element_book}. Therefore, integrating PINO with FEM to combine the efficiency of neural operators and the accuracy of classical solvers is of critical importance. Several recent studies have explored this direction. Zhang et al.\ proposed HINT \citet{zhang2024blending}, which blends DeepONet with iterative solvers by exploiting the low-frequency approximation capability of neural networks and the high-frequency refinement of classical iterations. However, HINT faces challenges when applied to nonlinear systems. Lee et al.\ \citet{lee2025neural} extended this idea to nonlinear problems and improved Newton iteration efficiency, but still relied on large datasets for operator training. Wang et al.\ \citet{wang2024homogenius} used operator predictions as initial guesses to accelerate numerical solvers for heterogeneous solid mechanics problems, but their approach was data-driven and limited to linear settings. Harandi et al.\ \citet{harandi2025spifol} applied PINO to heterogeneous solid mechanics problems using purely physics-based training, but did not investigate iterative refinement. Taghikhani et al.\ \citet{taghikhani2025neural} extended PINO to nonlinear solid mechanics and employed PINO solutions as initial guesses for iterative solvers, significantly reducing iteration counts, but relied on FNO as the backbone. Eshaghi et al.\ \citet{eshaghi2025nows} further demonstrated the effectiveness of PINO-based warm starts but within the FNO framework.

These studies collectively indicate that using PINO solutions as initial guesses can substantially accelerate traditional PDE solvers, combining the efficiency of neural operators with the accuracy of classical numerical methods, while potentially eliminating the need for training data. However, to date, there remains a lack of research on transformer-based PINO frameworks equipped with a pretraining and warm-start mechanism.

To address this gap, we propose the \textbf{Pretrained Finite Element Method} (PFEM). PFEM adopts the Transolver architecture to solve PDEs by first training a physics-informed neural operator and then using its predictions as initial guesses for FEM-based iterative solvers. PFEM consists of two stages: a pretraining stage and a warm-start stage. In the pretraining stage, a Transolver-based PINO is trained solely using governing physical equations, without any labeled data. This stage supports arbitrary spatial resolutions and requires only point-cloud inputs, rather than structured grids as in FNO, greatly enhancing flexibility, particularly for complex geometries. In the warm-start stage, the pretrained neural operator predictions are used as initial guesses in classical iterative solvers, which progressively refine the solution toward the exact one. Because the initial guesses are already physics-consistent and close to the true solution, the number of required iterations is significantly reduced, leading to substantial computational savings.

We conduct extensive generalization studies of PFEM across different geometries, boundary conditions, and material distributions, covering both linear and nonlinear iterative solvers. The results demonstrate that PFEM simultaneously achieves the efficiency of neural operator learning and the high accuracy of traditional numerical methods. Moreover, as the number of training samples increases, PFEM exhibits progressively stronger generalization capability, indicating that it is a self-improving framework. Overall, PFEM constitutes a novel pretraining and warm-start framework for computational physics.

The main contributions of PFEM can be summarized as follows:
\begin{itemize}
	\item \textbf{Point-cloud input}: All point-wise information, including spatial coordinates, geometry, boundary conditions, and material properties, is encoded as prompts and processed by the Transolver architecture.
	\item \textbf{Discretization invariance}: Training and testing resolutions can differ arbitrarily, enabling cross-resolution generalization without retraining.
	\item \textbf{Pure physics-based training}: A Transolver-based PINO is trained solely using governing physical equations, without relying on labeled data.
	\item \textbf{Warm-start via classical solvers}: PFEM predictions are used as initial guesses in iterative solvers, significantly reducing iteration counts while achieving arbitrarily high accuracy.
\end{itemize}

To the best of our knowledge, PFEM is the first physics-informed neural operator that integrates a pretraining--warm-start paradigm and enables simple encoding-based generalization across geometries, boundary conditions, and material properties, built upon the Transolver framework.

\section{Prerequisite knowledge}
In this section, we first provide a brief introduction to operator learning. We then introduce the physics-informed neural operator (PINO), which incorporates governing physical equations into the training of neural operators.

\subsection{Input and output of neural operators}

The central objective of neural operators is to learn mappings between function spaces, which has broad applications in science and engineering \citep{kovachki2023neural}. For instance, the input function may represent boundary conditions of a partial differential equation (PDE), while the output corresponds to the solution field of interest. Although functions are infinite-dimensional objects in a strict mathematical sense, in practice they are approximated using discrete samples. Fortunately, neural operators typically exhibit discretization invariance, meaning that they can generalize across different spatial resolutions and discretizations.

Traditional neural operators are trained in a purely data-driven manner. This approach is conceptually straightforward. Suppose that for different geometries $\boldsymbol{G}$, material distributions $\boldsymbol{M}$, and boundary conditions $\boldsymbol{B}$, we obtain a dataset
$\{\boldsymbol{X}^{(i)}, \boldsymbol{Y}^{(i)}\}_{i=1}^{N}$,
where $\boldsymbol{X}^{(i)}$ and $\boldsymbol{Y}^{(i)}$ denote the input and output functions, respectively, and $N$ is the total number of samples. The data-driven training of a neural operator can be formulated as
\begin{equation}
	\begin{aligned}
		\mathcal{\boldsymbol{H}}_{\boldsymbol{\theta}}(\boldsymbol{x}; \boldsymbol{X})
		&= \arg\min_{\boldsymbol{\theta}} \mathcal{L}_{\mathrm{data}}, \\
		\mathcal{L}_{\mathrm{data}}
		&= \frac{1}{N} \sum_{i=1}^{N}
		\left\| \mathcal{\boldsymbol{H}}(\boldsymbol{X}^{(i)}; \boldsymbol{\theta})
		- \boldsymbol{Y}^{(i)} \right\|^{2},
	\end{aligned}
\end{equation}
where $\boldsymbol{\theta}$ denotes the trainable parameters of the neural operator
$\mathcal{\boldsymbol{H}}(\cdot; \boldsymbol{\theta})$, and $\|\cdot\|$ represents the $\ell_2$ norm.
For simplicity, the neural operator approximation
$\mathcal{\boldsymbol{H}}(\boldsymbol{X}^{(i)}; \boldsymbol{\theta})$
is denoted as
$\mathcal{\boldsymbol{H}}_{\boldsymbol{\theta}}(\boldsymbol{x}; \boldsymbol{X}^{(i)})$,
where $\boldsymbol{x}$ represents spatial coordinates.
Both $\boldsymbol{X}^{(i)}$ and $\boldsymbol{Y}^{(i)}$ are functions sampled at discrete spatial locations.
The primary advantage of purely data-driven operator learning lies in its simplicity and flexibility.
However, its major limitation is the requirement for large amounts of high-quality training data, which are often expensive to generate and may contain significant numerical or experimental noise.

\begin{figure}[t]
	\centering
	\includegraphics[scale=0.65]{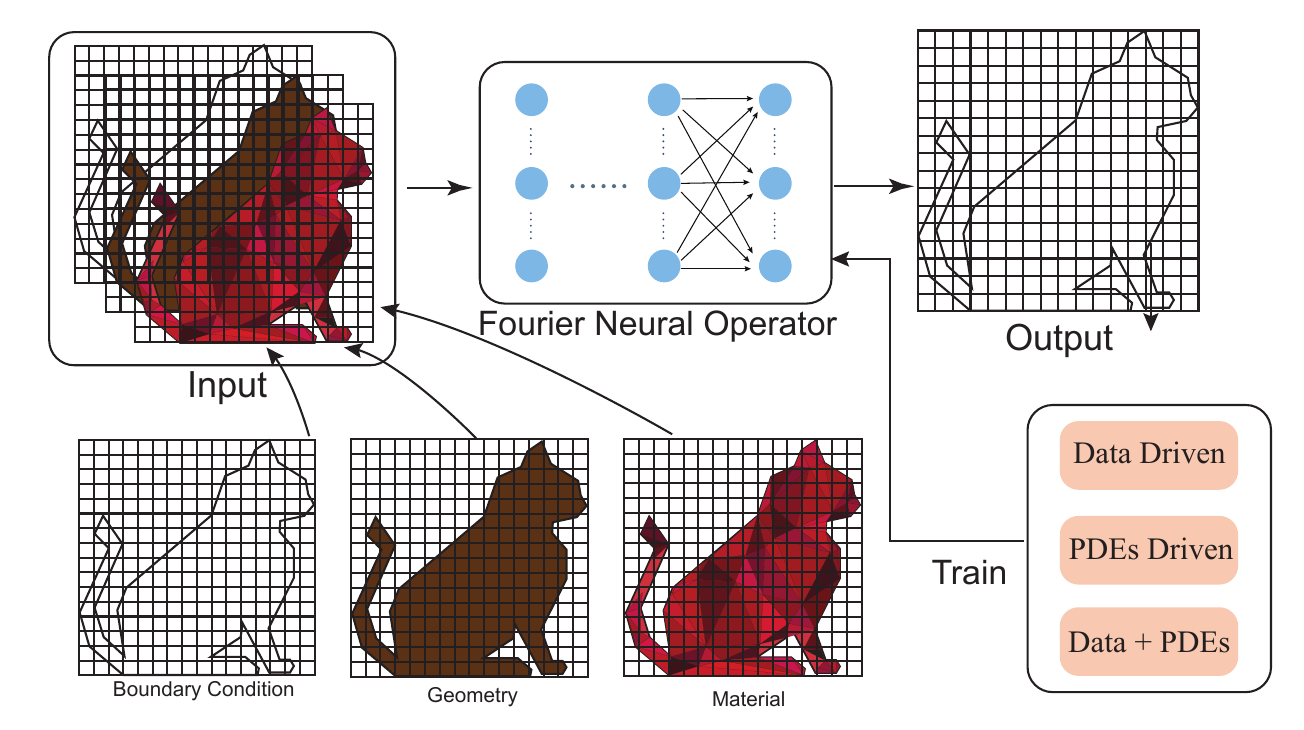}
	\caption{Illustration of the input and output representation in the Fourier Neural Operator (FNO).}
	\label{fig:Input-and-output_FNO}
\end{figure}

Among existing neural operator architectures, the Fourier Neural Operator (FNO) \citep{li2020fourier} is one of the most representative approaches.
The input and output of FNO are typically defined on structured grids and represented as high-dimensional tensors, as illustrated in \Cref{fig:Input-and-output_FNO}.
Consider a simple example in which the material distribution of a cat-shaped domain in \Cref{fig:Input-and-output_FNO} is encoded as a binary field, with a value of $1$ indicating the presence of material and $0$ otherwise.
Such a representation highlights a key limitation of FNO: its ability to accurately describe complex geometries is restricted, often leading to approximation errors.
Although several extensions, such as Geo-FNO \cite{li2022fourier}, have been proposed to improve geometric handling, generalization to highly complex geometries remains challenging.

This naturally raises the question of whether it is possible to abandon structured grid representations altogether.
Unfortunately, such a modification is difficult within the original FNO framework, as its core efficiency relies on the fast Fourier transform (FFT), which requires data defined on regular grids.
This limitation is intrinsic to the FFT-based design of FNO.
To overcome this issue, we adopt a transformer-based architecture to reformulate operator learning, which is introduced in detail in \Cref{subsec:Transolver}.
By representing all relevant information for PDE solving using point clouds, this approach removes the restriction of structured grids and enables greater flexibility in handling complex geometries.

\subsection{Physics-informed neural operator}
\label{subsec:Physics-informed-Neural-Operator}

As discussed in the previous subsection, traditional neural operators are typically trained in a fully data-driven manner, where high-fidelity solution data must be generated in advance. This requirement poses a significant challenge for many scientific and engineering applications in which data are scarce or expensive to obtain. To address this limitation, it is of great importance to train neural operators directly using governing physical equations, giving rise to the Physics-informed Neural Operator (PINO) \citep{li2024physics,wang2021learning,eshaghi2025variational}. A schematic illustration of this paradigm is shown in the ``PDEs Driven'' branch of \Cref{fig:Input-and-output_FNO}. Notably, PINO also allows for hybrid training strategies that combine data and PDE constraints, often referred to as ``Data + PDEs'' training.

Most problems in computational physics are governed by well-defined PDEs, and the core objective of computational mechanics is the numerical solution of these equations. The key idea of the ``PDEs Driven'' paradigm is therefore to embed the governing PDEs directly into the training of neural operators. In PINO, the loss function is typically defined as
\begin{equation}
	\begin{aligned}
		\mathcal{\boldsymbol{H}}_{\boldsymbol{\theta}}(\boldsymbol{x}; \boldsymbol{X})
		&= \arg\min_{\boldsymbol{\theta}} \mathcal{L}_{\mathrm{pde}}, \\
		\mathcal{L}_{\mathrm{pde}}
		&= \frac{1}{N} \sum_{i=1}^{N}
		\left\{
		\frac{1}{N_{d}} \sum_{j=1}^{N_{d}}
		\left\| \boldsymbol{P}
		\big(
		\mathcal{\boldsymbol{H}}_{\boldsymbol{\theta}}
		(\boldsymbol{x}^{(j)}; \boldsymbol{X}^{(i)})
		\big) \right\|^{2}
		+
		\frac{1}{N_{b}} \sum_{j=1}^{N_{b}}
		\left\| \boldsymbol{I}
		\big(
		\mathcal{\boldsymbol{H}}_{\boldsymbol{\theta}}
		(\boldsymbol{x}^{(j)}; \boldsymbol{X}^{(i)})
		\big) \right\|^{2}
		\right\},
	\end{aligned}
	\label{eq:PDEs_driven}
\end{equation}
where $\boldsymbol{P}$ and $\boldsymbol{I}$ denote the residual operators associated with the governing PDE in the domain and the boundary conditions, respectively. $N_d$ and $N_b$ represent the number of interior and boundary collocation points, and $N$ is the number of training samples.

For clarity, we consider a simple example of steady-state heat conduction in a two-dimensional heterogeneous medium:
\begin{equation}
	\begin{cases}
		-\nabla \cdot \big[ k(\boldsymbol{x}) \nabla T(\boldsymbol{x}) \big]
		= f(\boldsymbol{x}),
		& \boldsymbol{x} \in \Omega, \\
		T(\boldsymbol{x}) = \bar{T}(\boldsymbol{x}),
		& \boldsymbol{x} \in \Gamma^{T}, \\
		k(\boldsymbol{x}) \dfrac{\partial T(\boldsymbol{x})}{\partial \boldsymbol{n}}
		= \bar{q}(\boldsymbol{x}),
		& \boldsymbol{x} \in \Gamma^{q},
	\end{cases}
	\label{eq:poisson_equation}
\end{equation}
where $k(\boldsymbol{x})$ denotes the spatially varying thermal conductivity, $T(\boldsymbol{x})$ is the temperature field, $f(\boldsymbol{x})$ is the heat source, and $\bar{T}(\boldsymbol{x})$ and $\bar{q}(\boldsymbol{x})$ are the prescribed temperature and heat flux on the Dirichlet boundary $\Gamma^{T}$ and Neumann boundary $\Gamma^{q}$, respectively. In this setting, operator learning aims to approximate the mapping
\[
\{ k(\boldsymbol{x}), f(\boldsymbol{x}), \bar{T}(\boldsymbol{x}), \bar{q}(\boldsymbol{x}) \}
\;\xrightarrow{\;\mathcal{\boldsymbol{H}}_{\boldsymbol{\theta}}\;}\;
T(\boldsymbol{x}).
\]

The strong-form PINO loss corresponding to \Cref{eq:poisson_equation} can be written as
\begin{equation}
	\begin{aligned}
		\mathcal{L}_{\mathrm{pde\text{-}s}}
		= \frac{1}{N} \sum_{i=1}^{N}
		\Bigg\{
		& \frac{1}{N_{d}} \sum_{j=1}^{N_{d}}
		\Big\|
		\nabla \cdot \Big[
		k^{(i)}(\boldsymbol{x}^{(j)})
		\nabla
		\mathcal{\boldsymbol{H}}_{\boldsymbol{\theta}}
		(\boldsymbol{x}^{(j)};
		k^{(i)}, f^{(i)}, \bar{T}^{(i)}, \bar{q}^{(i)})
		\Big]
		+ f^{(i)}(\boldsymbol{x}^{(j)})
		\Big\|^{2}
		\\
		& + \frac{1}{N_{bT}} \sum_{j=1}^{N_{bT}}
		\Big\|
		\mathcal{\boldsymbol{H}}_{\boldsymbol{\theta}}
		(\boldsymbol{x}^{(j)};
		k^{(i)}, f^{(i)}, \bar{T}^{(i)}, \bar{q}^{(i)})
		- \bar{T}^{(i)}(\boldsymbol{x}^{(j)})
		\Big\|^{2}
		\\
		& + \frac{1}{N_{bq}} \sum_{j=1}^{N_{bq}}
		\Big\|
		k^{(i)}(\boldsymbol{x}^{(j)})
		\boldsymbol{n} \cdot
		\nabla
		\mathcal{\boldsymbol{H}}_{\boldsymbol{\theta}}
		(\boldsymbol{x}^{(j)};
		k^{(i)}, f^{(i)}, \bar{T}^{(i)}, \bar{q}^{(i)})
		- \bar{q}^{(i)}(\boldsymbol{x}^{(j)})
		\Big\|^{2}
		\Bigg\},
	\end{aligned}
	\label{eq:strong_form_PDEs}
\end{equation}
where $N_{bT}$ and $N_{bq}$ denote the number of collocation points on the Dirichlet and Neumann boundaries, respectively. The loss in \Cref{eq:strong_form_PDEs} corresponds to the strong-form formulation of PINO \citep{li2024physics}.

Alternatively, when a variational principle is available, the loss function can be formulated in an energy-based manner, leading to the variational physics-informed neural operator (VINO) \citep{eshaghi2025variational}:
\begin{equation}
	\begin{aligned}
		\mathcal{L}_{\mathrm{pde\text{-}v}}
		= \frac{1}{N} \sum_{i=1}^{N}
		\Bigg\{
		& \int_{\Omega}
		\frac{1}{2}
		k^{(i)}
		\nabla
		\mathcal{\boldsymbol{H}}_{\boldsymbol{\theta}}
		(\boldsymbol{x};
		k^{(i)}, f^{(i)}, \bar{T}^{(i)}, \bar{q}^{(i)})
		\cdot
		\nabla
		\mathcal{\boldsymbol{H}}_{\boldsymbol{\theta}}
		(\boldsymbol{x};
		k^{(i)}, f^{(i)}, \bar{T}^{(i)}, \bar{q}^{(i)})
		\, \mathrm{d}\Omega
		\\
		& - \int_{\Omega}
		f^{(i)}
		\mathcal{\boldsymbol{H}}_{\boldsymbol{\theta}}
		(\boldsymbol{x};
		k^{(i)}, f^{(i)}, \bar{T}^{(i)}, \bar{q}^{(i)})
		\, \mathrm{d}\Omega
		\\
		& - \int_{\Gamma^{q}}
		\bar{q}^{(i)}
		\mathcal{\boldsymbol{H}}_{\boldsymbol{\theta}}
		(\boldsymbol{x};
		k^{(i)}, f^{(i)}, \bar{T}^{(i)}, \bar{q}^{(i)})
		\, \mathrm{d}\Gamma
		\Bigg\}.
	\end{aligned}
	\label{eq:energy_form_PDEs}
\end{equation}

The choice between the strong-form loss $\mathcal{L}_{\mathrm{pde\text{-}s}}$ and the variational loss $\mathcal{L}_{\mathrm{pde\text{-}v}}$ depends on the nature of the governing PDE. In general, an energy formulation exists only for steady-state problems, and not all PDEs admit a variational principle; see Appendix~A of \citet{wang2025physics} for a detailed discussion. Empirically, VINO often exhibits better efficiency and accuracy than strong-form PINO, primarily because the variational formulation involves derivatives of lower order, which leads to improved numerical stability during training.

\section{Method}
\label{sec:Method}

\subsection{Transolver}
\label{subsec:Transolver}

\begin{figure}
	\begin{centering}
		\includegraphics[scale=0.45]{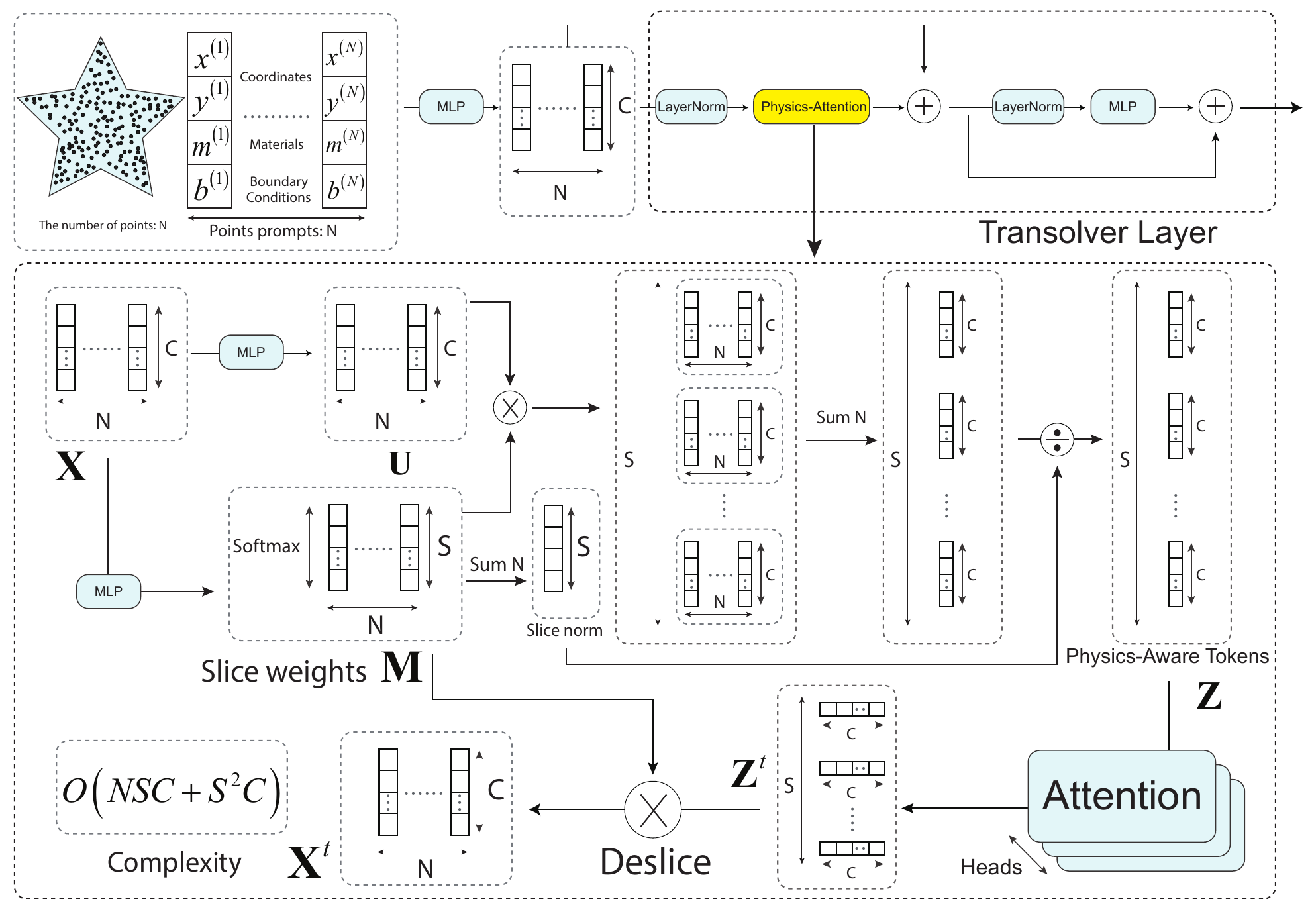}
		\par\end{centering}
	\caption{The architecture of Transolver layer \label{fig:Transolver} }
\end{figure}

The core idea of Transolver is to perform attention on a small set of physics-aware tokens with low computational complexity, rather than directly attending to all input points. The overall workflow is illustrated in \Cref{fig:Transolver}. In the original GNOT framework \citep{hao2023gnot}, self-attention is applied directly to point-wise features, resulting in a quadratic computational complexity of $\mathcal{O}(N^{2})$, where $N$ denotes the number of input points. In contrast, Transolver \citep{wu2024transolver} performs attention on physics-aware tokens, reducing the complexity to $\mathcal{O}(N + S^{2})$, where $S$ is the number of physics-aware tokens. Since typically $S \ll N$, Transolver exhibits substantial scalability and efficiency advantages for complex large-scale problems.

We now describe the Transolver architecture in detail. First, all points in the computational domain are encoded as
$\{x^{(i)}, y^{(i)}, m^{(i)}, b^{(i)}\}_{i=1}^{N}$,
where $x^{(i)}$ and $y^{(i)}$ denote the spatial coordinates, $m^{(i)}$ represents material information at the corresponding location (e.g., elastic modulus), and $b^{(i)}$ encodes boundary condition information. Notably, these points are not required to lie on a structured grid. Each point is mapped through a multilayer perceptron (MLP), followed by a LayerNorm operation to improve training stability and normalize feature statistics. This yields a point-wise feature representation
$
\boldsymbol{X} \in \mathbb{R}^{N \times C},
$
where $C$ denotes the number of feature channels.

The core component of Transolver is the physics-attention module. The feature matrix $\boldsymbol{X}$ is processed by two separate MLPs to produce
$
\boldsymbol{U} \in \mathbb{R}^{N \times C}
$
and
$
\boldsymbol{M} \in \mathbb{R}^{N \times S},
$
where $S$ denotes the number of physics-aware tokens. The slice weights $\boldsymbol{M}$ are normalized along the $S$ dimension using a softmax operation, such that they can be interpreted as attention weights. The physics-aware tokens
$
\boldsymbol{Z} \in \mathbb{R}^{S \times C}
$
are then obtained via a weighted aggregation of $\boldsymbol{U}$:
\begin{equation}
	Z_{JK}
	=
	\frac{\sum_{I=1}^{N} M_{IJ} U_{IK}}
	{\sum_{I=1}^{N} M_{IJ}} .
\end{equation}
Since $S \ll N$, the dimensionality of $\boldsymbol{Z}$ is significantly smaller than that of $\boldsymbol{X}$.

Multi-head self-attention is subsequently applied along the token dimension of $\boldsymbol{Z}$:
\begin{equation}
	\boldsymbol{Q}, \boldsymbol{K}, \boldsymbol{V}
	= \mathrm{MLP}(\boldsymbol{Z}),
	\qquad
	\boldsymbol{Z}^{t}
	=
	\mathrm{softmax}(\boldsymbol{Q}\boldsymbol{K}^{\mathsf{T}})\boldsymbol{V},
\end{equation}
yielding the transformed physics-aware tokens
$\boldsymbol{Z}^{t} \in \mathbb{R}^{S \times C}$.
The final physics-attention output
$\boldsymbol{X}^{t} \in \mathbb{R}^{N \times C}$
is obtained via a deslicing operation using the same weights $\boldsymbol{M}$:
\begin{equation}
	X^{t}_{IK}
	=
	\sum_{J=1}^{S} M_{IJ} Z^{t}_{JK}.
\end{equation}

The remaining operations follow the standard Transolver architecture shown in \Cref{fig:Transolver}. The essential contribution of the physics-attention module is that it reduces the quadratic complexity $\mathcal{O}(N^{2})$ of point-wise self-attention to $\mathcal{O}(N + S^{2})$ by operating on a compact set of physics-aware tokens. At the same time, Transolver retains strong approximation capability. As demonstrated in \citet{wu2024transolver}, Transolver exhibits excellent performance on three-dimensional problems with complex geometries, including benchmarks such as ShapeNet Car.

Multiple Transolver layers can be stacked in a manner analogous to Fourier layers in the Fourier Neural Operator (FNO), which typically improves expressiveness and generalization at the cost of increased computational overhead. The number of Transolver layers and the token size $S$ can therefore be adjusted according to the complexity of the target problem: more complex physics generally requires deeper architectures and larger $S$. From a physical perspective, the physics-aware tokens can be interpreted as a coarse representation of the underlying physical regions, replacing dense point-wise descriptions. This design is motivated by the observation that physical fields often exhibit spatial continuity and correlation, rendering point-wise information highly redundant.

\subsection{Pretrained Finite Element Method}
\label{subsec:PFEM}

\begin{figure}[t]
	\centering
	\includegraphics[scale=0.55]{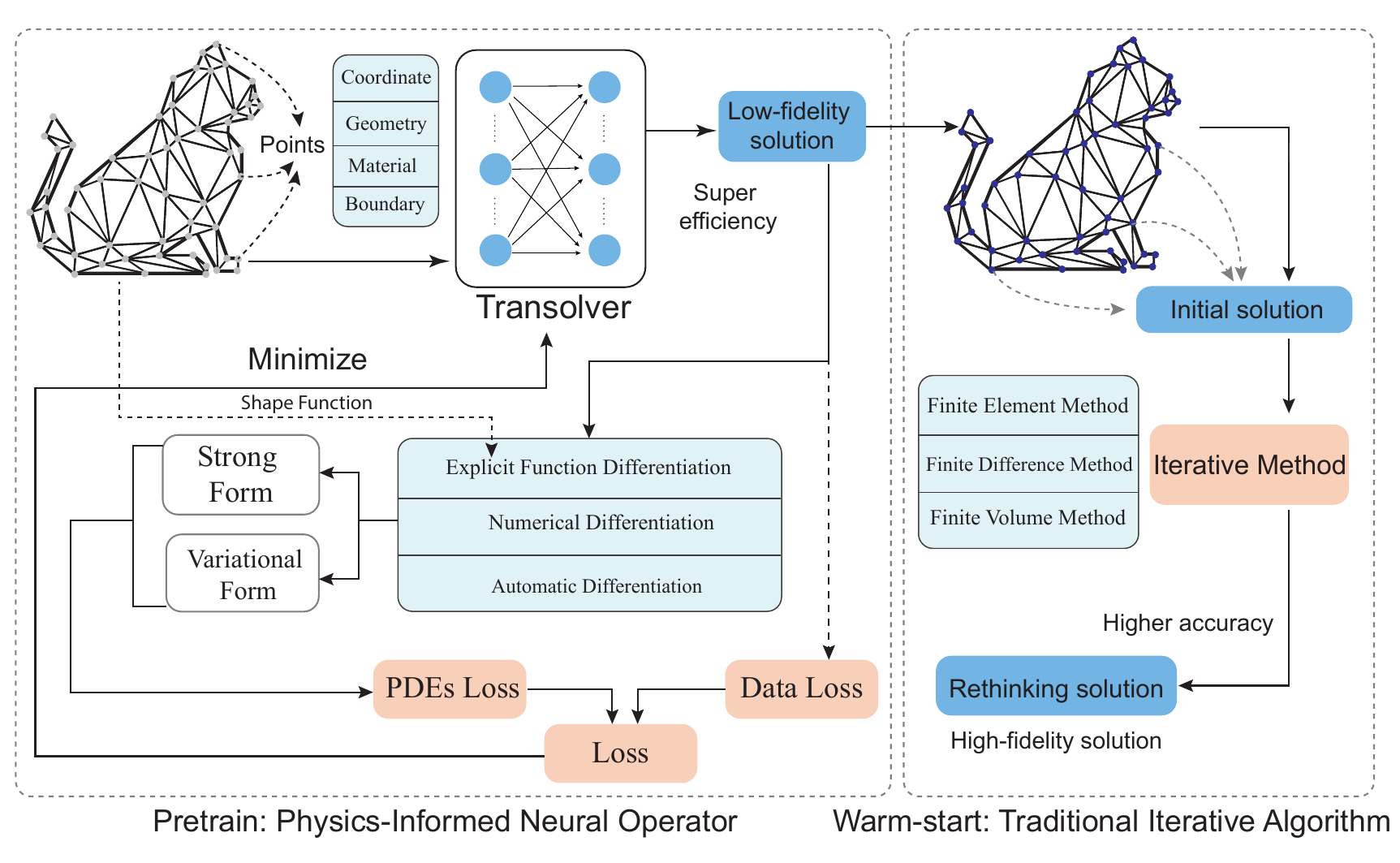}
	\caption{Schematic illustration of the Pretrained Finite Element Method (PFEM). PFEM consists of a pretraining stage and a warm-start stage. In the pretraining stage, a Transolver-based physics-informed neural operator (PINO) is trained using governing physical equations; the dashed line indicates that data-driven training is optional but not required. This stage yields a low-fidelity solution with extremely high efficiency. Spatial derivatives are constructed via explicit function differentiation, and either strong-form or variational PDE losses are employed. In the warm-start stage, classical numerical solvers (e.g., FEM, finite difference, or finite volume methods) are cast into iterative forms, and the neural operator prediction is used as the initial guess to obtain a high-fidelity solution.}
	\label{fig:PFEM}
\end{figure}

Owing to its strong approximation capability, particularly for point-cloud representations, Transolver naturally lends itself to a physics-driven pretraining strategy. Motivated by this observation, we propose the \textbf{Pretrained Finite Element Method} (PFEM), whose core idea is to train a Transolver-based neural operator purely using governing physical equations and subsequently use its predictions as initial guesses for classical finite element iterations. An overview of PFEM is illustrated in \Cref{fig:PFEM}.

PFEM consists of two stages: a pretraining stage and a warm-start stage. Notably, both stages rely exclusively on physical equations and do not require labeled solution data. In the following, we describe these two stages in detail.

\subsubsection{Pretraining}

The pretraining stage of PFEM is essentially a physics-informed neural operator (PINO) with Transolver as the backbone. While PINO can, in principle, be trained using a combination of data and physics, in this work we focus exclusively on \textbf{purely physics-driven training}.

The workflow begins by discretizing the computational domain using a mesh. All nodal information is then converted into a point-cloud representation
$
\{x^{(i)}, y^{(i)}, m^{(i)}, b^{(i)}\}_{i=1}^{N},
$
where $x^{(i)}$ and $y^{(i)}$ denote spatial coordinates, $m^{(i)}$ represents material properties at the corresponding location (e.g., elastic modulus), and $b^{(i)}$ encodes boundary condition information. Unlike conventional neural operators that require structured grids, PFEM operates directly on irregular point clouds.

Following the Transolver architecture described in \Cref{subsec:Transolver}, the point cloud is mapped to the quantities of interest at the same spatial locations,
$
\{x^{(i)}, y^{(i)}\}_{i=1}^{N}
\;\longrightarrow\;
\{ f_1^{(i)}, f_2^{(i)}, \ldots, f_n^{(i)} \}_{i=1}^{N},
$
which constitute a low-fidelity solution. To enforce physical consistency, spatial derivatives of the predicted fields are constructed directly from this low-fidelity solution via explicit function differentiation. Specifically, since the input originates from a finite element discretization, derivatives can be explicitly computed using standard finite element shape functions. Compared to automatic differentiation, explicit function differentiation significantly reduces computational cost, as discussed in detail in \ref{sec:computational-graph}.

Finally, appropriate PDE loss functions $\mathcal{L}_{\mathrm{pde}}$ are selected, typically in either strong or variational form, as introduced in \Cref{subsec:Physics-informed-Neural-Operator}. The trainable parameters of the Transolver backbone are then updated by minimizing the chosen PDEs loss (strong or variational form), completing the pretraining stage of PFEM.

\subsubsection{Warm-start}

\begin{figure}[t]
	\centering
	\includegraphics[scale=0.50]{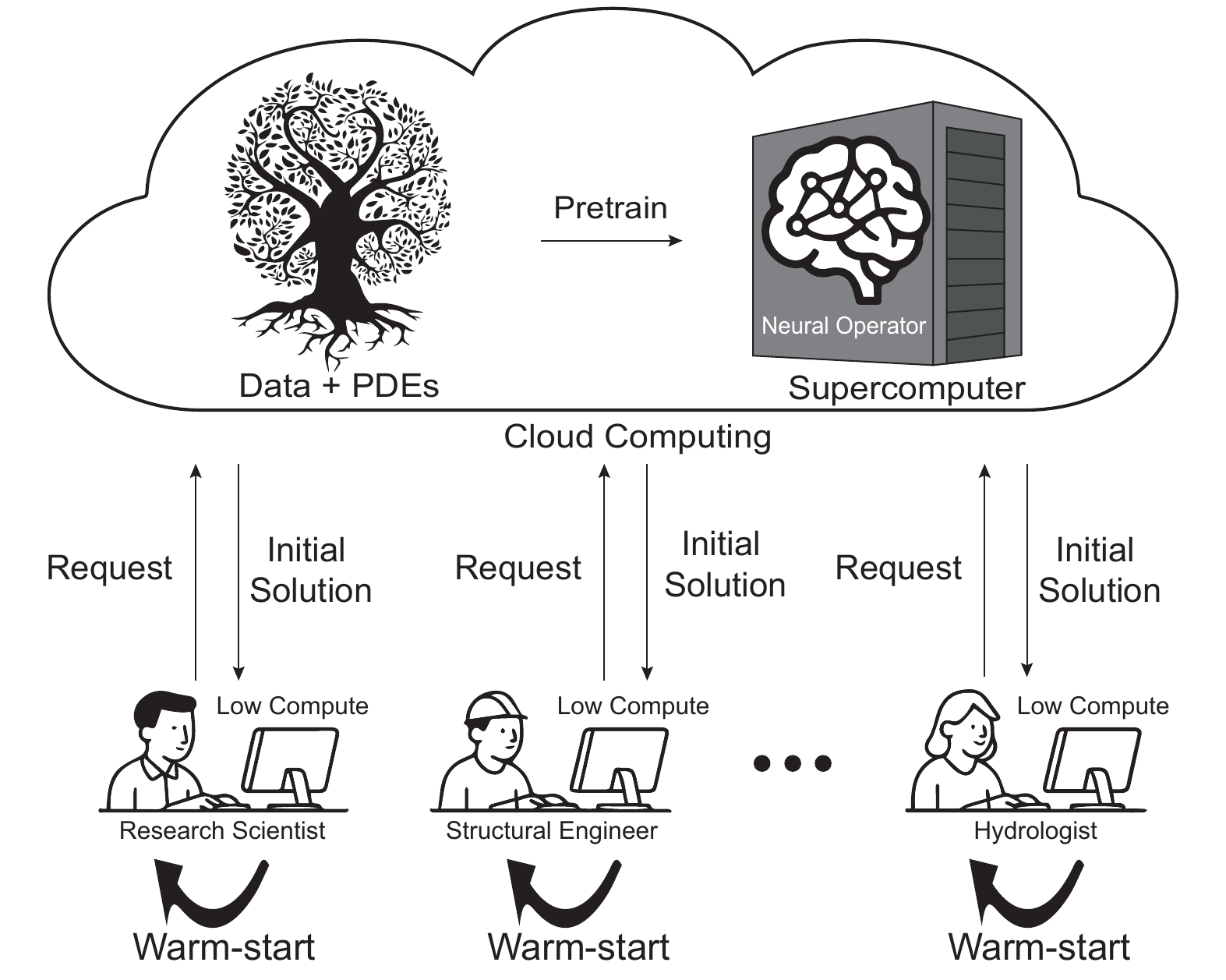}
	\caption{Schematic illustration of the warmstart phase of PFEM.}
	\label{fig:Warmstart_PFEM}
\end{figure}

The core idea of the warm-start phase of PFEM is to use the output of the pretraining stage as the initial guess (warm-start) of an iterative numerical solver and to subsequently refine this solution through classical iterations, as illustrated in \Cref{fig:Warmstart_PFEM}. In this paradigm, the pretraining stage of PFEM can be regarded as the core component deployed on a cloud-computing platform. When researchers from different application domains encounter specific PDE problems, they can submit their problem descriptions to the platform, which rapidly provides an initial solution via the pretrained model. Each user then performs local iterative refinement using their own computational resources to efficiently obtain a high-accuracy numerical solution. This ``cloud-based pretraining + local warm-start'' paradigm substantially improves overall computational efficiency.

From a technical perspective, the warm-start stage of PFEM is fundamentally based on iterative numerical algorithms. In this work, we focus on iterative solvers within the finite element framework. Both linear and nonlinear finite element iterations are reviewed in \ref{sec:Iterative-FEM}. The key idea of PFEM warm-start is to set the pretrained neural operator prediction $\boldsymbol{U}^{\mathrm{NO}}$ as the initial guess $\boldsymbol{U}^{(0)}$ of the iterative solver:
\begin{equation}
	\boldsymbol{U}^{(k+1)}
	=
	\boldsymbol{\phi}
	\big(
	\boldsymbol{U}^{(k)};
	\boldsymbol{K}^{(k)},
	\boldsymbol{F}^{(k)},
	\boldsymbol{U}^{(0)}=\boldsymbol{U}^{\mathrm{NO}}
	\big),
	\label{eq:iteraitive_linear_TINO}
\end{equation}
where $\boldsymbol{U}$ denotes the nodal values of the unknown field (e.g., displacement), $\boldsymbol{F}$ represents the nodal external force vector determined by boundary conditions, and $\boldsymbol{K}$ is the global stiffness matrix depending on geometry and material properties. The operator $\boldsymbol{\phi}(\cdot)$ denotes the specific iterative algorithm. The iteration is terminated once the convergence criterion
\begin{equation}
	\big\| \boldsymbol{r}(\boldsymbol{U}^{(k+1)}) \big\| < \mathrm{tol}
	\label{eq:tol_TINO}
\end{equation}
is satisfied, where $\boldsymbol{r}$ is the residual measure. In computational mechanics, the residual typically takes the form
$\boldsymbol{r} = \boldsymbol{f}^{\mathrm{int}}(\boldsymbol{U}) - \boldsymbol{f}^{\mathrm{ext}}(\boldsymbol{U})$,
with $\boldsymbol{f}^{\mathrm{int}}$ and $\boldsymbol{f}^{\mathrm{ext}}$ denoting the internal and external nodal force vectors, respectively. Here, $\mathrm{tol}$ denotes the prescribed convergence threshold.

Since the pretrained PFEM model can provide predictions at arbitrary spatial locations, the mesh nodes used in the warm-start stage can be directly queried to obtain nodal predictions. These predictions are then used as the initial guess for the finite element iterations. Conceptually, the pretraining and warm-start stages of PFEM are analogous to the pretraining and reinforcement learning stages of large language models \citep{guo2025deepseek}. Because the pretraining stage yields a high-quality, physics-consistent initial solution, the number of iterations required in the warm-start stage can be significantly reduced compared to traditional approaches that employ random or zero initial vectors. A reduction in iteration count directly translates into lower computational cost.

From a theoretical standpoint, the warm-start solution can achieve arbitrary accuracy, provided that the convergence threshold $\mathrm{tol}$ is chosen sufficiently small. It is worth emphasizing that the warm-start stage is optional and can be applied depending on accuracy requirements. To determine whether warm-start is necessary, we propose evaluating the initial solution provided by the pretraining stage by substituting it into \Cref{eq:tol_TINO} and examining the residual norm $\|\boldsymbol{r}\|$. If
$
\|\boldsymbol{r}\| < \mathrm{tol}_{\mathrm{fine}},
$
the warm-start stage can be skipped, where $\mathrm{tol}_{\mathrm{fine}}$ is a user-defined threshold for deciding whether further refinement is required.

If a moderate level of error is acceptable, the pretrained PFEM model can be used directly as a fast solver. When high-precision solutions are desired, the warm-start stage can be activated at the cost of additional computational effort. It is worth noting that for linear problems, a high-quality initial guess reduces the upper bound on the number of required iterations; a rigorous analysis is provided in \ref{sec:PFEM_warmstart_benefit}. For nonlinear problems, the quality of the initial guess is equally critical. Taking Newton's method as an example, if the initial solution lies within the attraction basin of the exact solution, the iteration enters a local fast-convergence regime with a convergence order approaching quadratic, thereby significantly improving overall computational efficiency.

\section{Results}
\label{sec:Result}

\section{Results}
\label{sec:Result}

\begin{table}[t]
	\centering
	\caption{Performance of PFEM on different benchmark problems. The relative error of PFEM refers to the error obtained in the pretraining stage. The reported PFEM time corresponds only to the neural operator prediction time, excluding the warm-start cost. The FEM time denotes the runtime of direct numerical solvers. The iteration numbers compare the warm-start iterations of PFEM with those of FEM using a zero initial guess. All iteration counts are evaluated under a tolerance of $\mathrm{tol}=10^{-3}$, corresponding to a relative error of $10^{-3}$. "G", "M", "B" refer to different geometries, materials, and boundary conditions}
	\label{tab:PFEM_arch_diff_problem}
	
	\begin{adjustbox}{max width=\textwidth}
		\begin{tabular}{cccccccc}
			\toprule
			Problem & Layers & Tokens & Learning rate &
			Relative error (mean/std) &
			Time: s (PFEM/FEM) &
			Iteration (Speedup: FEM/PFEM) &
			Optimizer \\
			\midrule
			Elastic plate (G)
			& 3 & 64 & 0.002 & 0.0113 / 0.0344
			& 0.0257 / 5.117
			& 6.24 = 519.88 / 83.32
			& Adam \\
			
			Elastic plate (G+M)
			& 4 & 128 & 0.002 & 0.0262 / 0.0183
			& 0.0633 / 5.826
			& 6.43 = 834.84 / 129.69
			& Adam \\
			
			Elastic plate (G+M+B)
			& 4 & 128 & 0.002 & 0.0209 / 0.0157
			& 0.0694 / 5.984
			& 9.18 = 910.08 / 99.06
			& Adam \\
			
			Hyperelastic beam (M+B)
			& 4 & 128 & 0.002 & 0.0661 / 0.0453
			& 0.0914 / 36.147
			& 3.66 = 4.90 / 1.34
			& Adam \\
			
			Hyperelastic Cook (G+M+B)
			& 4 & 128 & 0.002 & 0.0304 / 0.0460
			& 0.122 / 93.752
			& 5.41 = 8.38 / 1.55
			& Adam \\

			3D Homogenization (Geometry+material) & 4 & 64 & 0.0005 & 0.0949/0.0384 & 0.074/150.0 & 2.727=40.9/15.0 & Adam \\
			\bottomrule
		\end{tabular}
	\end{adjustbox}
\end{table}

\subsection{Linear elastic plate with arbitrary holes}

\begin{figure}[t]
	\centering
	\includegraphics[scale=0.55]{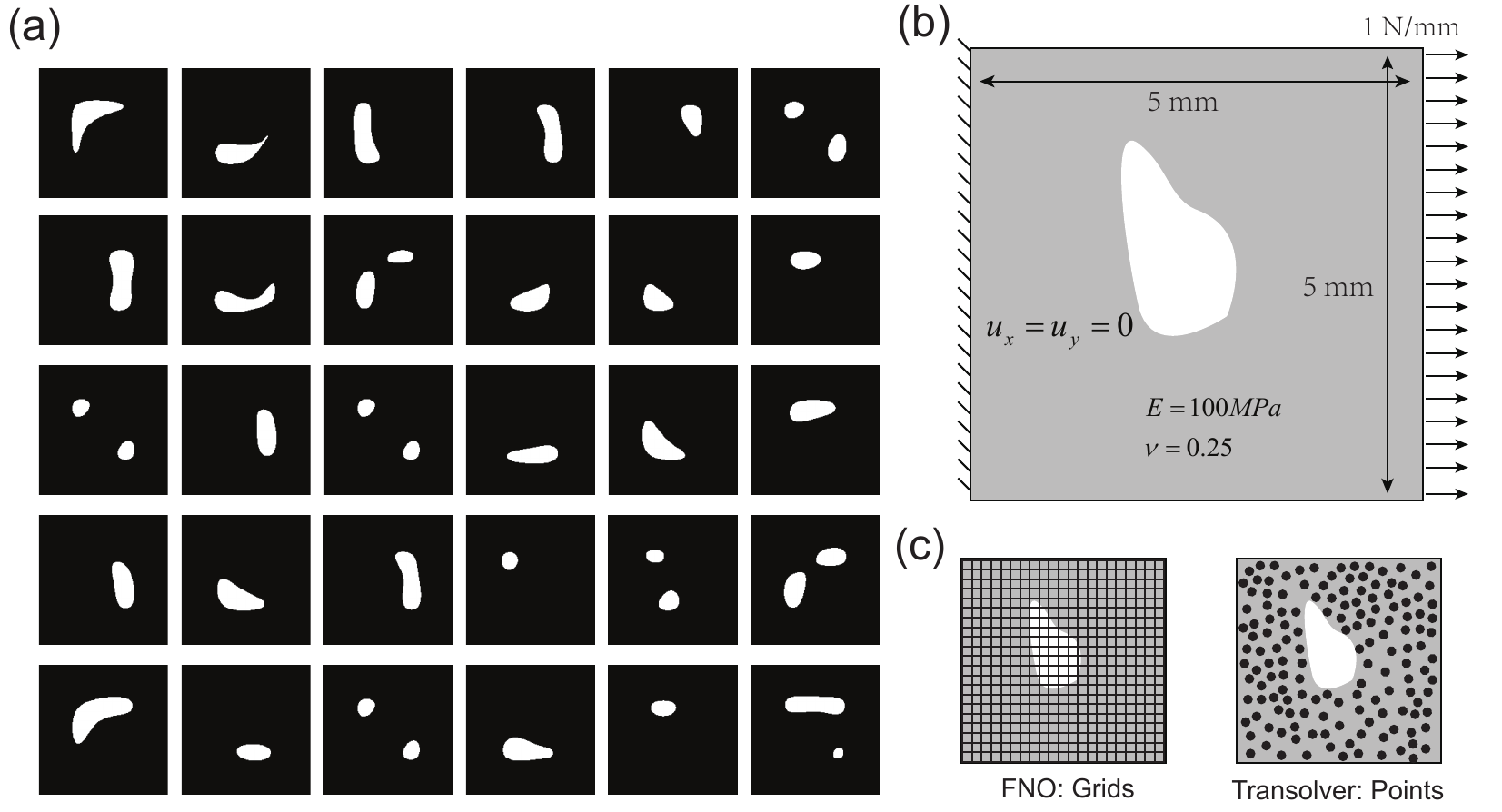}
	\caption{Elastic plate with arbitrary holes: (a) randomly generated hole geometries in the dataset; (b) boundary condition configuration; (c) comparison of input representations, where FNO requires structured grids while Transolver operates directly on point clouds.}
	\label{fig:intro_elasticity_hole}
\end{figure}

We first consider a linear elastic square plate with randomly distributed holes of arbitrary shapes. The plate domain is $\Omega = [0,5]^2$, and holes are randomly generated within the region $[1,4]^2$ using a Gaussian random field, resulting in highly irregular geometries, as illustrated in \Cref{fig:intro_elasticity_hole}a. The boundary conditions are prescribed as follows: the left boundary is fully clamped ($u_x = u_y = 0$), the right boundary is subjected to a constant traction ($t_x = 1$, $t_y = 0$), and the top and bottom boundaries are traction-free. The material parameters are set to Young's modulus $E = 100~\mathrm{MPa}$ and Poisson's ratio $\nu = 0.25$, and a plane stress assumption is adopted, as shown in \Cref{fig:intro_elasticity_hole}b.

The input of PFEM consists of the geometry, material distribution, and boundary conditions,
$\boldsymbol{G}, \boldsymbol{M}, \boldsymbol{B} = \{\boldsymbol{X}^{(i)}\}_{i=1}^{N}$,
where
$
\boldsymbol{X}^{(i)} = \{ x_j^{(i)}, y_j^{(i)}, E_j^{(i)}, \nu_j^{(i)}, t_{y,j}^{(i)} \}_{j=1}^{N^{(i)}}.
$
The output is the displacement field at the corresponding spatial locations,
$
\boldsymbol{U} = \{ u_j^{(i)}, v_j^{(i)} \}_{i=1}^{N},
$
where $u$ and $v$ denote the displacements in the $x$- and $y$-directions, respectively.

We first evaluate the geometric generalization capability of PFEM in the pretraining stage by fixing the material and boundary conditions ($E=100~\mathrm{MPa}$, $\nu=0.25$, $u_x=u_y=0$, $t_x=1$, $t_y=0$) and varying only the hole geometries. High-fidelity reference solutions are generated using FEM for all geometry configurations. The dataset consists of 1000 training samples and 300 testing samples. It is important to emphasize that no labeled solution data are used during training; only randomly generated geometries following the Gaussian random field distribution are provided.

Since PFEM is based on Transolver, it directly accepts point-cloud inputs, whereas conventional PINO methods based on FNO (e.g., VINO \citep{eshaghi2025variational}) require structured grid inputs. The input representations of PFEM and VINO are compared in \Cref{fig:intro_elasticity_hole}c. Due to the varying hole geometries, the FEM meshes differ across samples, each containing approximately $3.7 \times 10^{4}$ triangular elements. Mesh convergence studies are performed for all cases, ensuring the reliability of the FEM reference solutions.

During training, FEM nodal coordinates with different resolutions are used as PFEM inputs to predict the displacement fields in both spatial directions. The network is trained using an energy-based physics loss:
\begin{equation}
	\begin{aligned}
		\mathcal{L}
		&=
		\int_{\Omega}
		\big( \Psi - \boldsymbol{f} \cdot \boldsymbol{u} \big)\, \mathrm{d}\Omega
		-
		\int_{\Gamma^{\boldsymbol{t}}}
		\bar{\boldsymbol{t}} \cdot \boldsymbol{u}\, \mathrm{d}\Gamma, \\
		\Psi
		&=
		\frac{1}{2}\sigma_{ij}\varepsilon_{ij}, \\
		\sigma_{ij}
		&=
		\frac{E}{1+\nu}\varepsilon_{ij}
		+
		\frac{E\nu}{(1+\nu)(1-2\nu)}
		\varepsilon_{kk}\delta_{ij}, \\
		\varepsilon_{ij}
		&=
		\frac{1}{2}(u_{i,j} + u_{j,i}), \\
		\boldsymbol{u}
		&=
		\frac{\boldsymbol{x}}{L}
		\odot
		\mathcal{\boldsymbol{H}}_{\boldsymbol{\theta}}(\boldsymbol{x}; \boldsymbol{G}, \boldsymbol{\theta}),
	\end{aligned}
\end{equation}
where $\mathcal{L}$ denotes the loss function, $\Psi$ is the strain energy density, $\boldsymbol{\sigma}$ and $\boldsymbol{\varepsilon}$ are the stress and strain tensors, respectively, $\boldsymbol{f}$ is the body force, and $\bar{\boldsymbol{t}}$ is the prescribed traction on the Neumann boundary $\Gamma^{\boldsymbol{t}}$. The operator $\mathcal{\boldsymbol{H}}_{\boldsymbol{\theta}}$ denotes the Transolver-based neural operator, which takes the geometric information
$\boldsymbol{G}=\{\boldsymbol{X}^{(i)}\}_{i=1}^{N}$ with
$\boldsymbol{X}^{(i)}=\{x_j^{(i)}, y_j^{(i)}\}_{j=1}^{N^{(i)}}$
as input and outputs the displacement field at the corresponding spatial locations.

\begin{figure}[t]
	\centering
	\includegraphics[scale=0.55]{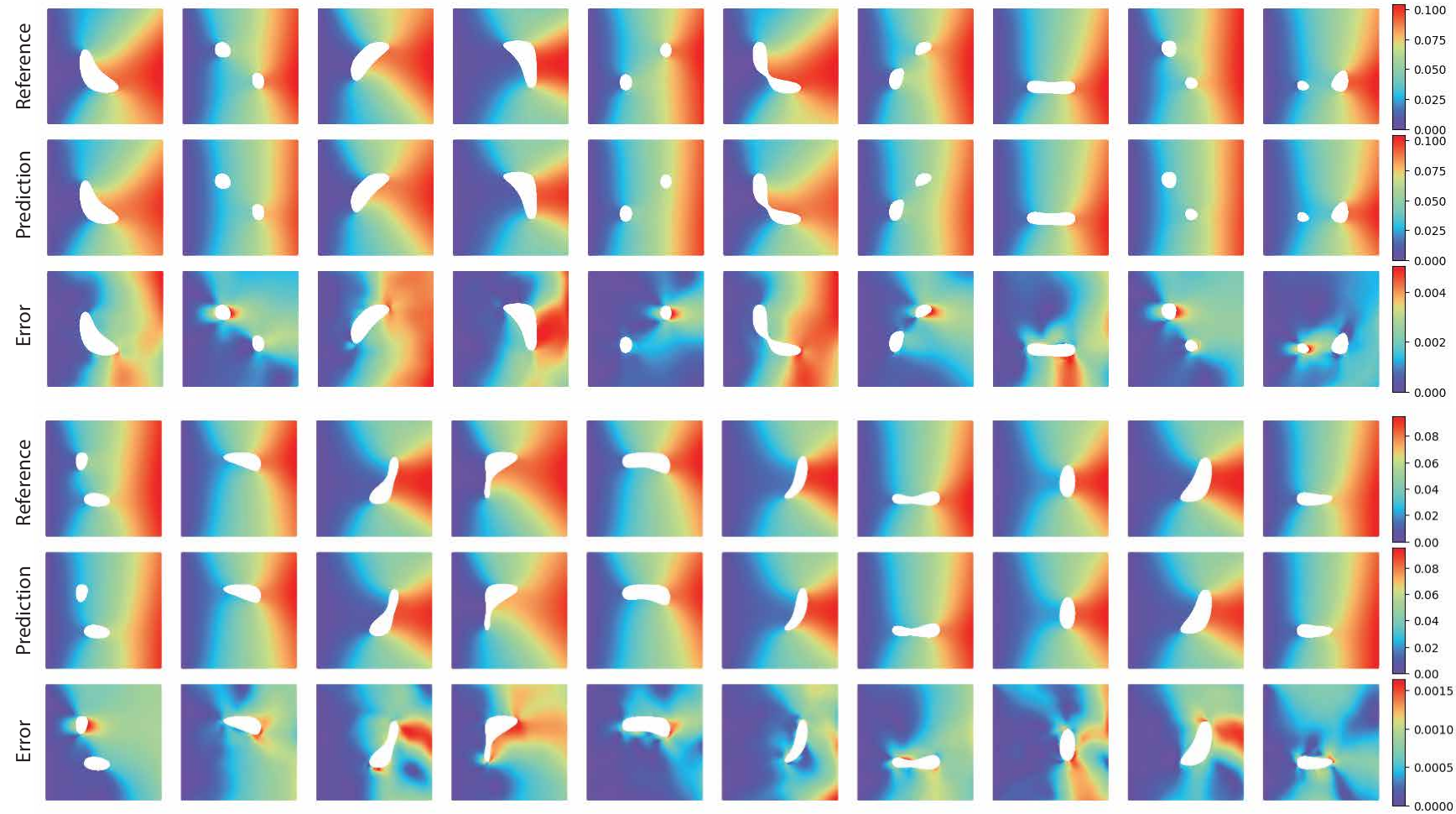}
	\caption{Absolute displacement magnitude predicted by PFEM for elastic plates with arbitrary hole geometries in the pretraining stage.}
	\label{fig:PFEM_Elasticity_hole_contourf_geo}
\end{figure}

\begin{figure}[t]
	\centering
	\includegraphics[scale=0.55]{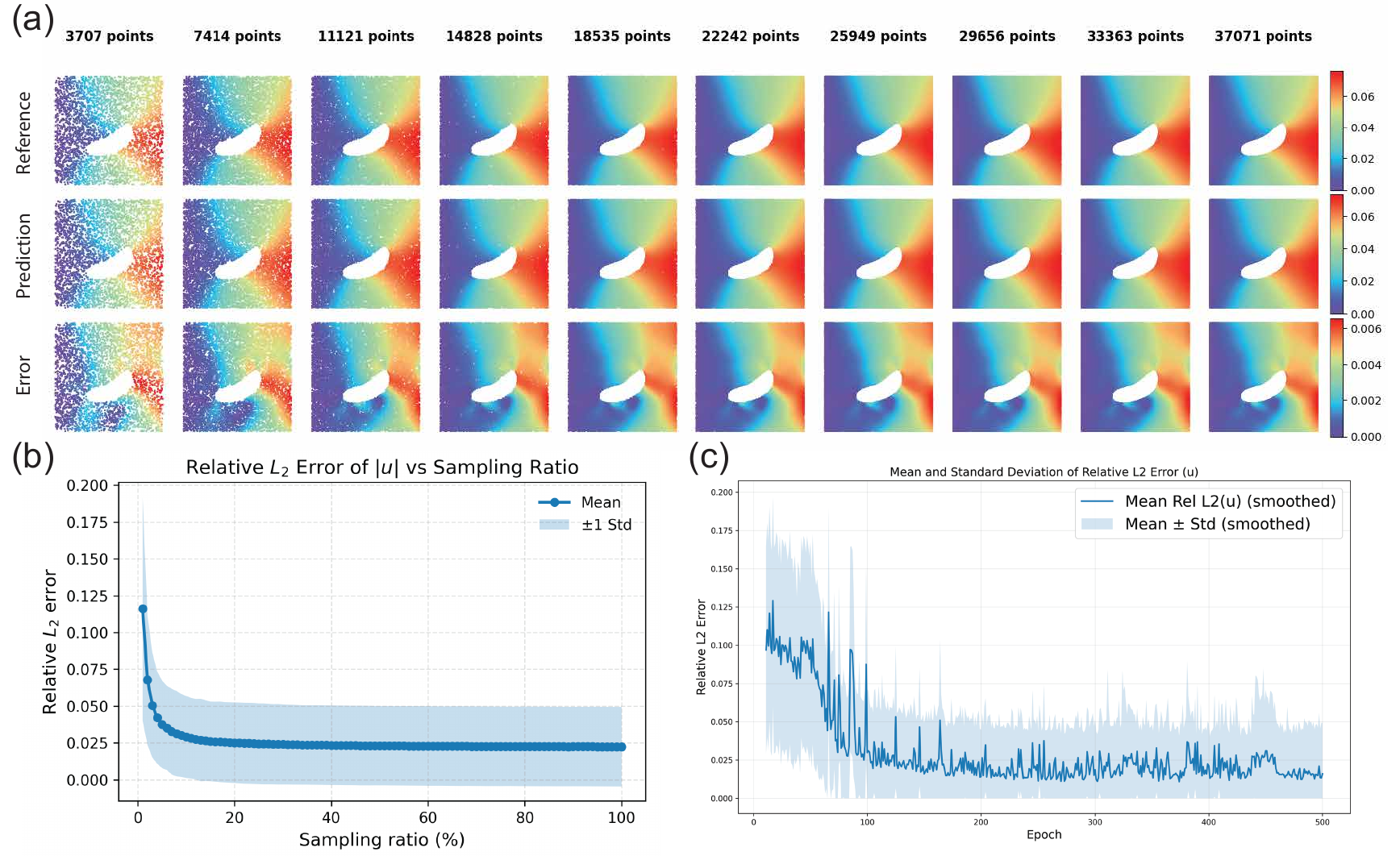}
	\caption{Performance of PFEM under different input resolutions for the elastic plate with holes: (a) predicted displacement magnitude using different proportions of test points (from $10\%$ to $100\%$ with an interval of $10\%$); (b) error distribution of $u_x$ under different point densities (from $1\%$ to $100\%$ with an interval of $1\%$); (c) convergence history of the test error in $u_x$ during training.}
	\label{fig:PFEM_Elasticity_hole_different_points}
\end{figure}

\Cref{fig:PFEM_Elasticity_hole_contourf_geo} presents the predictions obtained in the pretraining stage of PFEM for elastic plates with different hole geometries. PFEM demonstrates strong generalization capability across highly irregular geometries, with a relative error of approximately $1\%$. Quantitative statistics of the mean error and standard deviation over the test set are reported in \Cref{tab:PFEM_arch_diff_problem}. The convergence behavior of the test error during training is shown in \Cref{fig:PFEM_Elasticity_hole_different_points}c, where PFEM converges within approximately 100 epochs. The detailed model architecture and training hyperparameters are summarized in \Cref{tab:PFEM_arch_diff_problem}.

In \Cref{fig:PFEM_Elasticity_hole_different_points}a and b, we further evaluate the robustness of PFEM by testing the same geometry with different numbers of input points. The test inputs consist of point clouds with varying resolutions, while the underlying geometry remains fixed. The predicted displacement fields are nearly identical across different point densities, and the corresponding errors remain consistently low. This demonstrates the robustness of PFEM with respect to the number of input points. Notably, the point resolutions used during training and testing are completely different, indicating that PFEM possesses discretization invariance, i.e., resolution-independent generalization capability.

Since the Transolver backbone of PFEM involves several architectural hyperparameters, we further investigate the influence of the number of layers and the number of physics-aware tokens (slice size). \Cref{tab:TINO_Elasticity_layer_slice} summarizes the accuracy and efficiency of PFEM under different Transolver configurations. As the number of layers and tokens increases, both the prediction accuracy and robustness exhibit a clear improvement trend, at the expense of increased computational cost.

\begin{figure}
	\begin{centering}
		\includegraphics[scale=0.55]{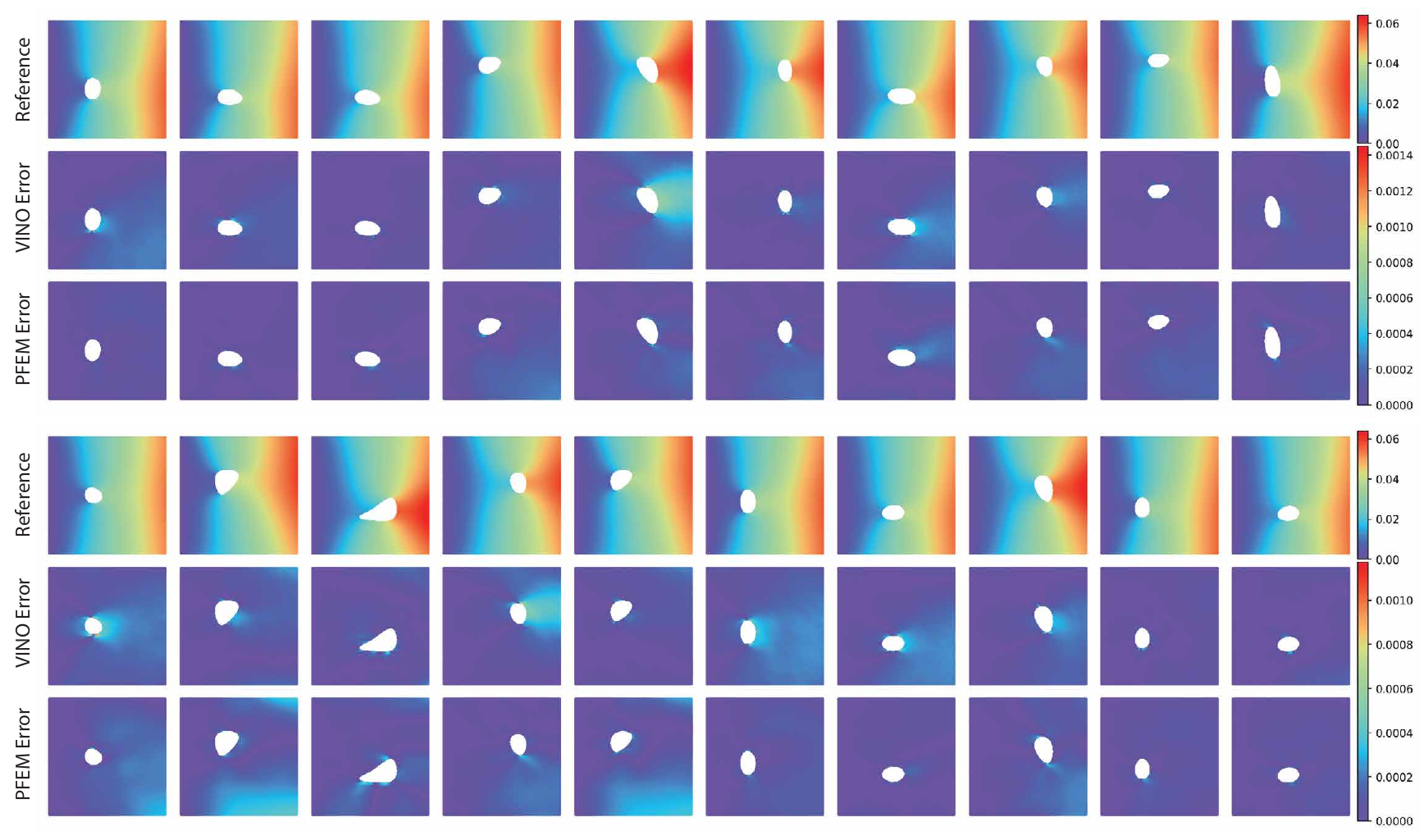}
		\par\end{centering}
	\caption{Comparison of absolute displacement field predictions between PFEM and VINO.}
	\label{fig:PFEM_VINO_Elasticity_hole_contourf_compare}
\end{figure}

For comparison, we also consider VINO, which adopts FNO as the backbone architecture, shown as \Cref{fig:PFEM_VINO_Elasticity_hole_contourf_compare}. The prediction accuracy of Transolver-based PFEM is comparable to that of FNO-based VINO; interested readers are referred to the original VINO paper for detailed results on this benchmark \citep{eshaghi2025variational}. However, while FNO requires structured grid inputs, Transolver operates directly on unstructured point clouds and achieves comparable accuracy. This highlights the superior flexibility of PFEM with respect to input representations, especially for problems involving complex geometries.

\begin{table}[t]
	\centering
	\caption{Accuracy and efficiency of PFEM with different Transolver architectures. The reported time corresponds to the training time per epoch, where each epoch consists of 1000 training samples.}
	\label{tab:TINO_Elasticity_layer_slice}
	\begin{tabular}{ccccc}
		\toprule
		Layers & Tokens & Parameters &
		Relative error $u_x$ (mean/std) & Time (s/epoch) \\
		\midrule
		\multirow{3}{*}{3}
		& 32  & 287{,}482 & 0.0131 / 0.0320 & 60.41 \\
		& 64  & 289{,}114 & 0.0113 / 0.0344 & 81.36 \\
		& 128 & 292{,}378 & 0.00999 / 0.0293 & 123.8 \\
		\multirow{3}{*}{4}
		& 32  & 371{,}874 & 0.0164 / 0.0361 & 80.71 \\
		& 64  & 374{,}050 & 0.0114 / 0.0302 & 107.1 \\
		& 128 & 378{,}402 & 0.0110 / 0.0290 & 161.1 \\
		\multirow{3}{*}{5}
		& 32  & 456{,}266 & 0.0113 / 0.0280 & 97.4 \\
		& 64  & 458{,}986 & 0.0106 / 0.0262 & 131.6 \\
		& 128 & 464{,}426 & \textbf{0.00995} / \textbf{0.0255} & 197.8 \\
		\bottomrule
	\end{tabular}
\end{table}

The above experiments demonstrate the performance of PFEM in the pretraining stage for generalization across different geometries. We next investigate the effectiveness of the warm-start stage of PFEM. All warm-start experiments are conducted on 300 test cases using the conjugate gradient (CG) method as the iterative solver. The prediction obtained from the pretraining stage of PFEM is employed as the initial guess of the iterative algorithm.

It is worth emphasizing that the warm-start stage is not mandatory. Whether warm-start is required is determined by a prescribed threshold $\mathrm{tol}_{\mathrm{fine}}$. Specifically, the prediction from the pretraining stage is substituted into the residual
$
\boldsymbol{r} = \boldsymbol{f}^{\mathrm{int}}(\boldsymbol{U}) - \boldsymbol{f}^{\mathrm{ext}}(\boldsymbol{U}).
$
If $\|\boldsymbol{r}\| < \mathrm{tol}_{\mathrm{fine}}$, the PFEM prediction is considered sufficiently accurate and no further iteration is required. Otherwise, the pretrained solution is used as the initial guess for the finite element iterative solver.

\begin{figure}[t]
	\centering
	\includegraphics[scale=0.55]{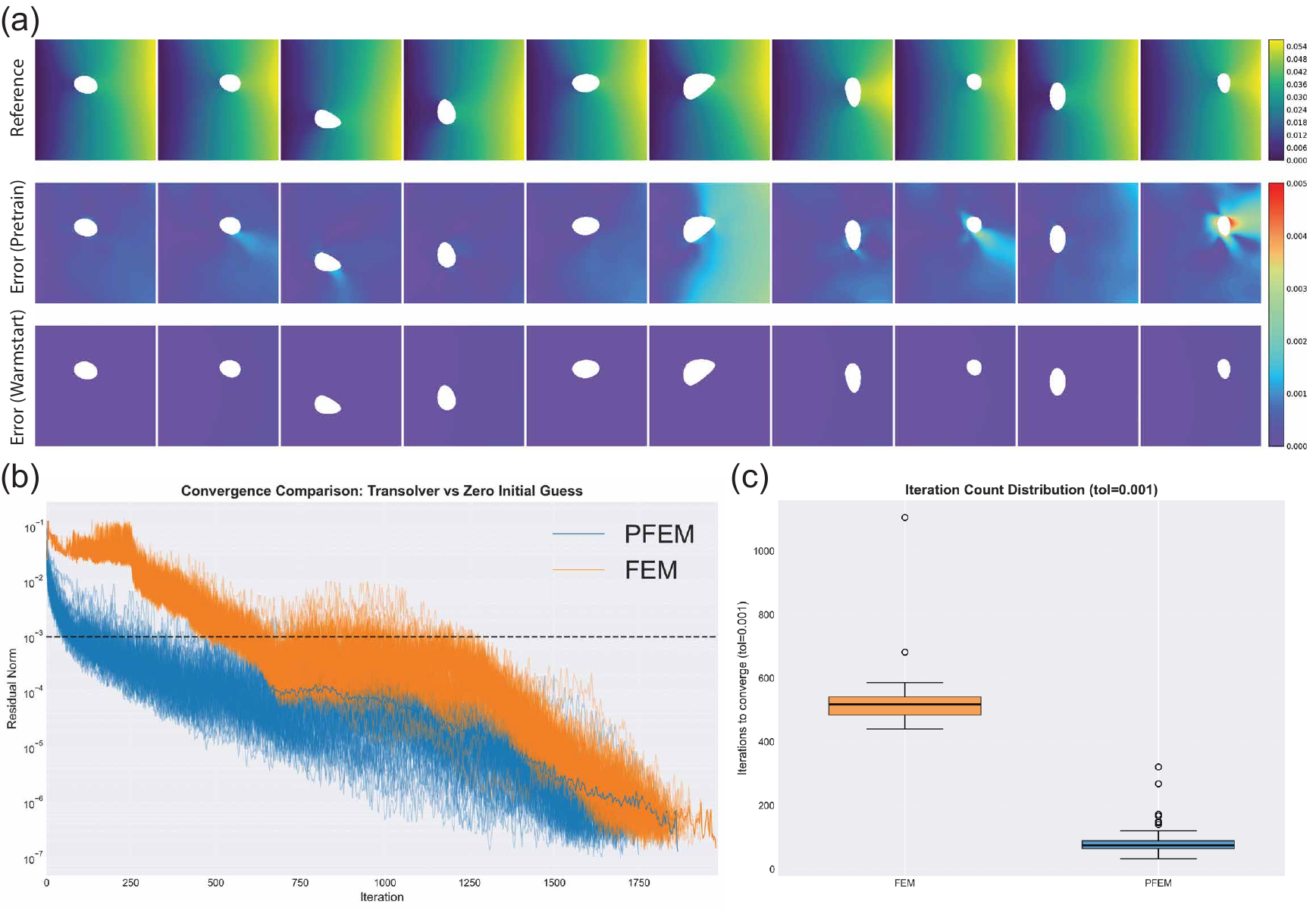}
	\caption{Warm-start results of PFEM for elastic plates with arbitrary hole geometries: (a) comparison of error contours between the pretraining and warm-start stages; (b) convergence histories of PFEM and FEM; (c) comparison of the number of iterations required by PFEM and FEM.}
	\label{fig:PFEM_Elasticity_hole_contourf_geo_warmstart}
\end{figure}

\Cref{fig:PFEM_Elasticity_hole_contourf_geo_warmstart} presents the warm-start results for elastic plates with arbitrary hole geometries. In the traditional FEM baseline, the zero vector is used as the initial guess. It can be clearly observed that PFEM significantly reduces the number of iterations required for convergence. Since the same iterative algorithm (CG) is used for both PFEM and FEM, the convergence rates are identical; the performance gain of PFEM stems solely from the improved initial guess.

For a convergence tolerance of $\mathrm{tol}=10^{-3}$, corresponding to a relative error of approximately $10^{-3}$, PFEM achieves a speedup factor of $6.24$ compared to FEM. Specifically, the average number of iterations required by PFEM is $83.32$, whereas FEM requires $519.88$ iterations. When the tolerance is  $\mathrm{tol}=10^{-6}$, the resulting relative error reaches the order of $10^{-7}$.

The above example highlights the strong approximation capability of PFEM for complex geometries. We further increase the difficulty by evaluating the generalization performance of PFEM with respect to both geometry and material variations. In this setting, Young’s modulus and Poisson’s ratio are sampled according to Gaussian random fields. To account for the increased complexity arising from simultaneous geometric and material variations, the training set is expanded to 2000 samples, and the test set contains 500 samples. All cases pass mesh convergence studies, ensuring the reliability of the FEM reference solutions. Notably, no labeled displacement fields are used during training; the training data consist solely of geometry and material distributions generated from Gaussian random fields.

\begin{figure}[t]
	\centering
	\includegraphics[scale=0.55]{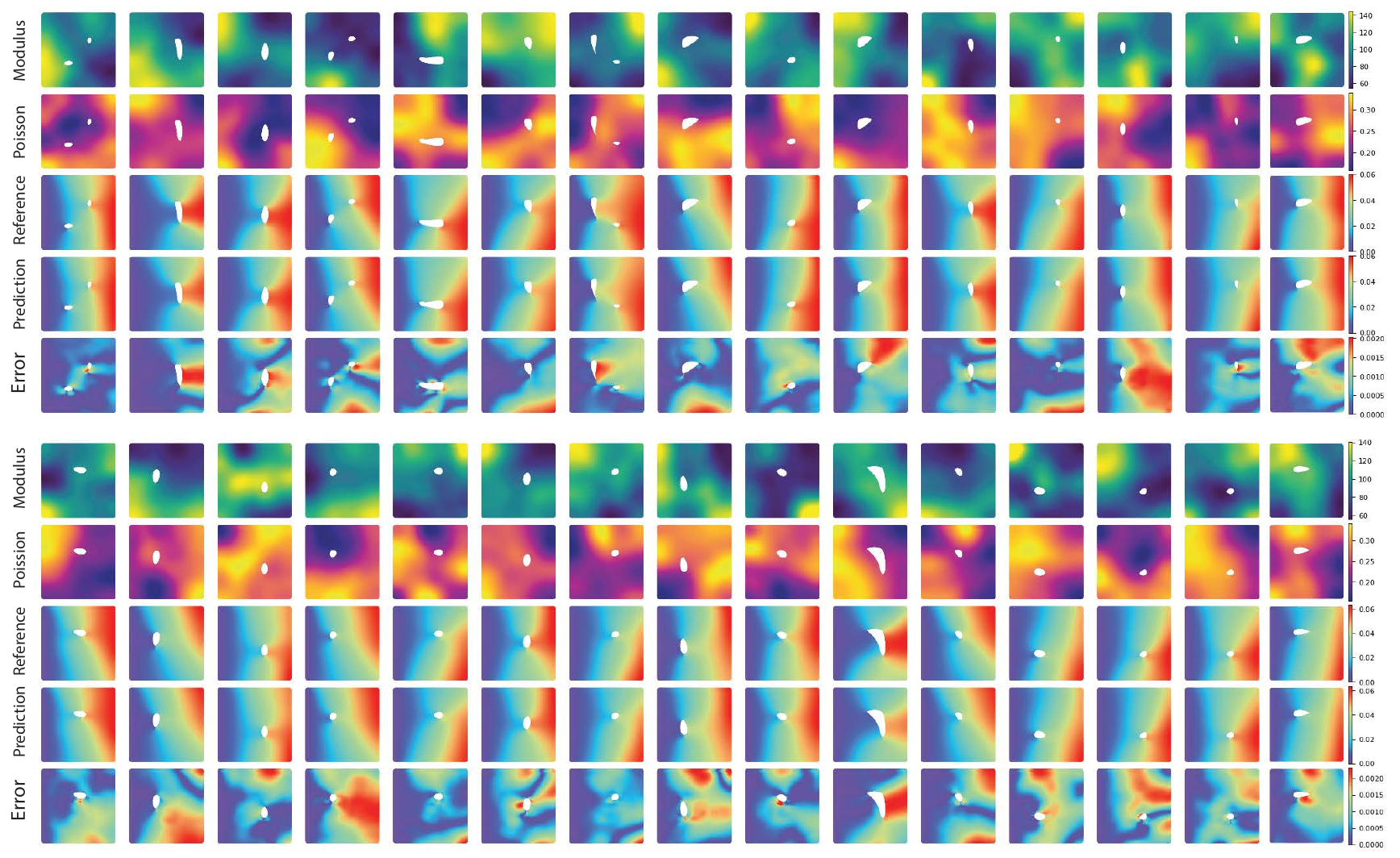}
	\caption{Predicted absolute displacement magnitude by PFEM for elastic plates with arbitrary hole geometries and heterogeneous material distributions.}
	\label{fig:PFEM_Elasticity_hole_contourf_geo_mater}
\end{figure}

The input of PFEM includes geometric and material information,
$\boldsymbol{G}, \boldsymbol{M} = \{\boldsymbol{X}^{(i)}\}_{i=1}^{N}$,
where
$
\boldsymbol{X}^{(i)} =
\{ x_j^{(i)}, y_j^{(i)}, E_j^{(i)}, \nu_j^{(i)} \}_{j=1}^{N^{(i)}}.
$
The output is the displacement field at the corresponding spatial locations,
$\boldsymbol{U} = \{ u_j^{(i)}, v_j^{(i)} \}_{i=1}^{N}$.
The pretraining results are shown in \Cref{fig:PFEM_Elasticity_hole_contourf_geo_mater}, which demonstrate that PFEM achieves a test error of approximately $1\%$ for joint geometry--material generalization; detailed quantitative results are provided in \Cref{tab:PFEM_arch_diff_problem}.

\begin{figure}[t]
	\centering
	\includegraphics[scale=0.55]{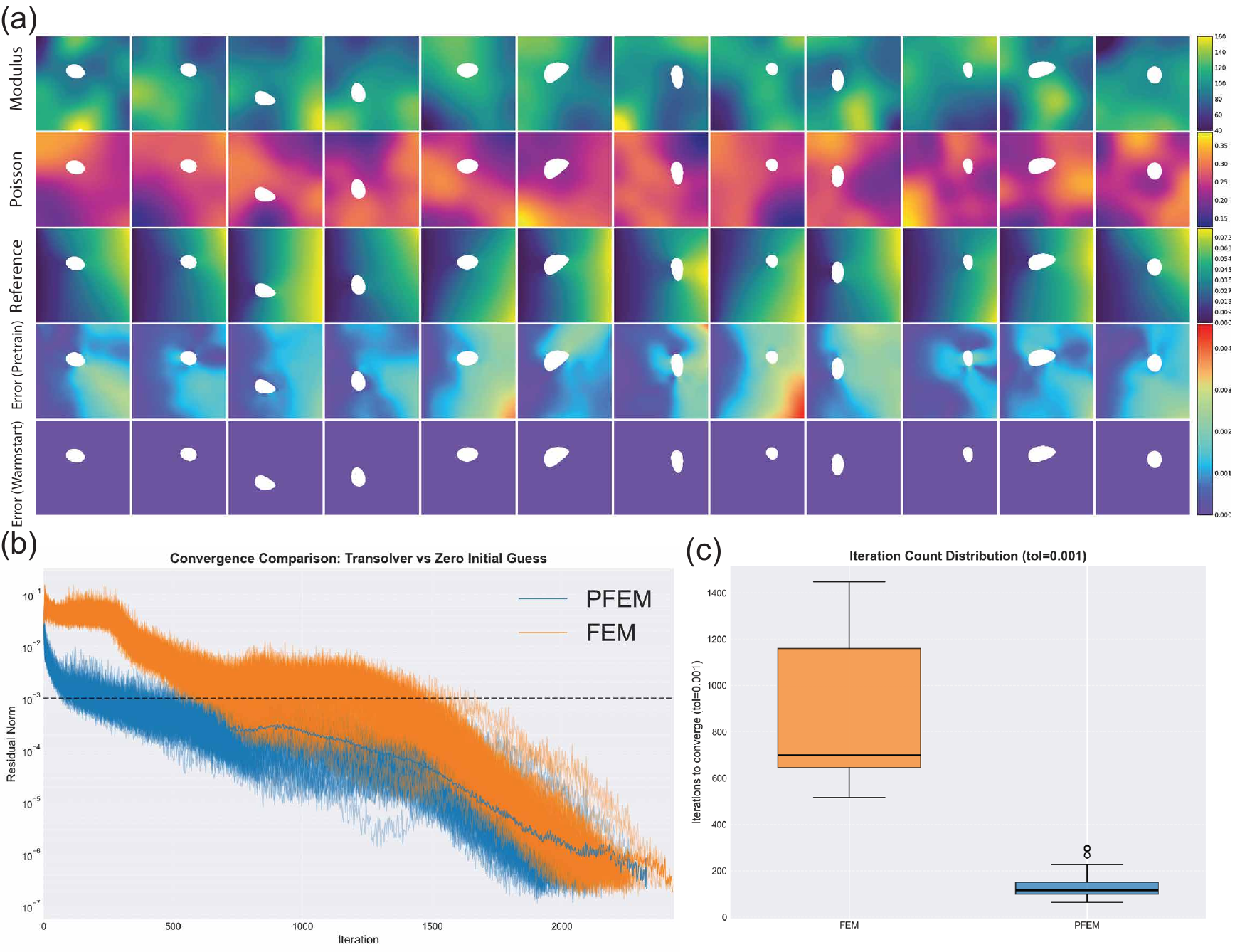}
	\caption{Warm-start results of PFEM for elastic plates with arbitrary hole geometries and heterogeneous materials: (a) comparison of error contours between the pretraining and warm-start stages; (b) convergence histories of PFEM and FEM; (c) comparison of the number of iterations required by PFEM and FEM.}
	\label{fig:PFEM_Elasticity_hole_contourf_geo_mater_warmstart}
\end{figure}

We further evaluate the warm-start performance of PFEM under combined geometric and material variations. Warm-start experiments are conducted on 500 test cases using the conjugate gradient method, with the PFEM pretraining prediction serving as the initial guess. As shown in \Cref{fig:PFEM_Elasticity_hole_contourf_geo_mater_warmstart}, PFEM again significantly reduces the number of iterations required for convergence. For a tolerance of $\mathrm{tol}=10^{-3}$, corresponding to a relative error of approximately $10^{-3}$, PFEM achieves a speedup factor of $6.43$ over the conventional FEM solver.

We further increase the level of difficulty by evaluating the generalization capability of PFEM with respect to geometry, material properties, and boundary conditions simultaneously. All three factors are sampled from Gaussian random fields. The training set consists of 1800 samples, and the test set contains 200 samples. All problems pass mesh convergence studies, ensuring the reliability of the FEM reference solutions. It is important to emphasize that the training set does not contain any labeled displacement fields; only geometry, material distributions, and boundary conditions are provided.

The input of PFEM includes geometric, material, and boundary information,
$\boldsymbol{G}, \boldsymbol{M}, \boldsymbol{B} = \{\boldsymbol{X}^{(i)}\}_{i=1}^{N}$,
where
$
\boldsymbol{X}^{(i)} =
\{ x_j^{(i)}, y_j^{(i)}, E_j^{(i)}, \nu_j^{(i)}, t_{y,j}^{(i)} \}_{j=1}^{N^{(i)}}.
$
The output is the displacement field at the corresponding spatial locations,
$\boldsymbol{U} = \{ u_j^{(i)}, v_j^{(i)} \}_{i=1}^{N}$.
The geometric and material variations follow the same setting as in the previous experiments, while the boundary condition is generalized by sampling the traction $t_y$ on the right boundary from a Gaussian random field.

\begin{figure}[t]
	\centering
	\includegraphics[scale=0.55]{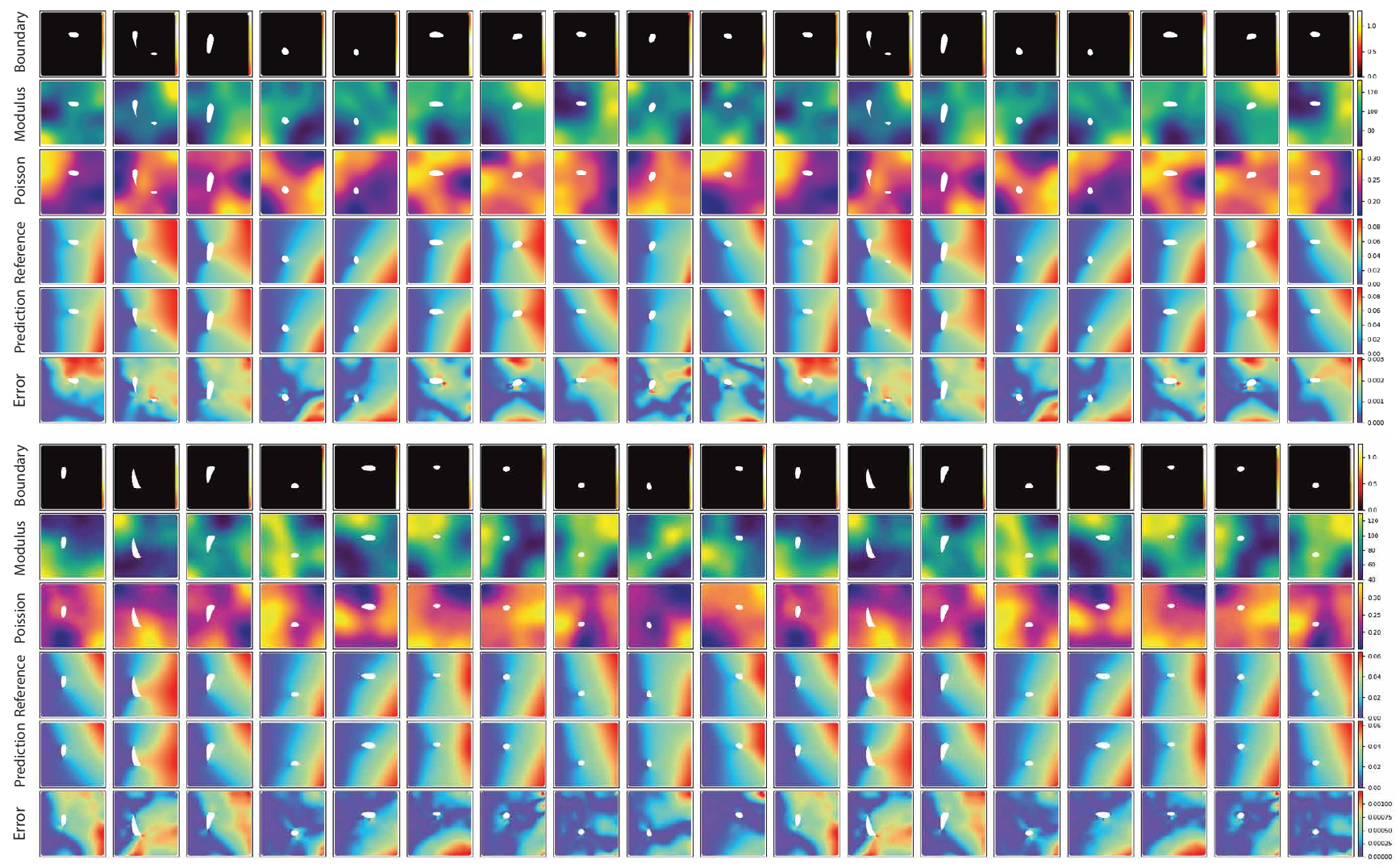}
	\caption{Predicted absolute displacement magnitude by PFEM for elastic plates with arbitrary hole geometries, heterogeneous material distributions, and varying boundary conditions.}
	\label{fig:PFEM_Elasticity_hole_contourf_geo_mater_boundary}
\end{figure}

The pretraining results are shown in \Cref{fig:PFEM_Elasticity_hole_contourf_geo_mater_boundary}. PFEM demonstrates strong predictive performance across arbitrary geometries, heterogeneous materials, and varying boundary conditions, with a relative error of approximately $1\%$. Detailed quantitative results are reported in \Cref{tab:PFEM_arch_diff_problem}. Notably, although the size of the training set is smaller than that used in the geometry--material generalization experiments, the test error does not deteriorate. This observation indicates that PFEM possesses excellent generalization capability and can maintain stable and reliable predictions even with a relatively small amount of training data.

\begin{figure}[t]
	\centering
	\includegraphics[scale=0.55]{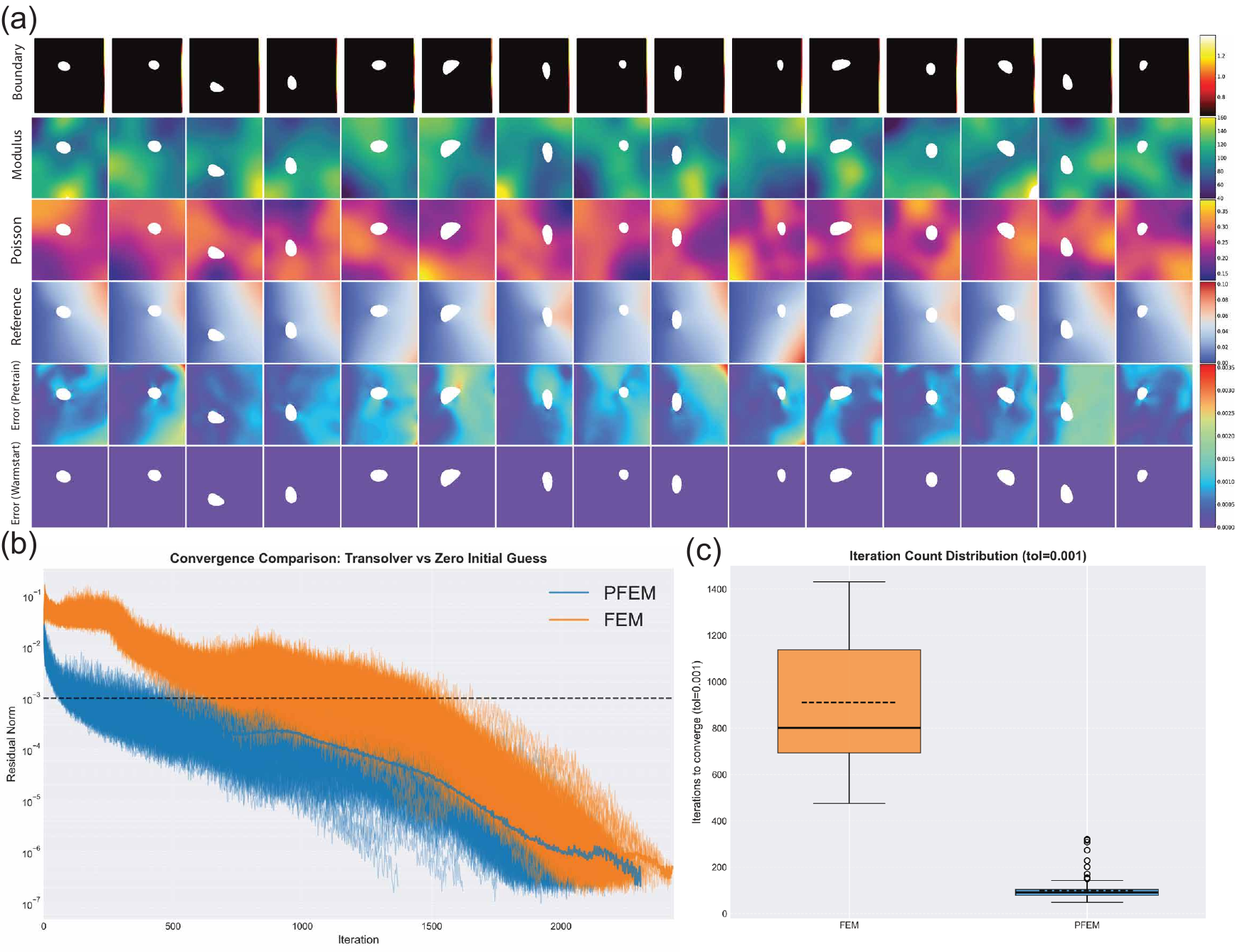}
	\caption{Warm-start results of PFEM for elastic plates with arbitrary hole geometries, heterogeneous materials, and varying boundary conditions: (a) comparison of error contours between the pretraining and warm-start stages; (b) convergence histories of PFEM and FEM; (c) comparison of the number of iterations required by PFEM and FEM.}
	\label{fig:PFEM_Elasticity_hole_contourf_geo_mater_boundary_warmstart}
\end{figure}

We next evaluate the warm-start performance of PFEM under simultaneous variations in geometry, material properties, and boundary conditions. Warm-start experiments are conducted on 200 test cases using the conjugate gradient method, with the PFEM pretraining prediction employed as the initial guess of the iterative solver. The warm-start results are presented in \Cref{fig:PFEM_Elasticity_hole_contourf_geo_mater_boundary_warmstart}. PFEM again leads to a substantial reduction in the number of iterations required for convergence.

For a convergence tolerance of $\mathrm{tol}=10^{-3}$, corresponding to a relative error of $10^{-3}$, PFEM achieves a speedup factor of $9.18$ compared to the conventional FEM solver. When the tolerance is further tightened to $\mathrm{tol}=10^{-6}$, the resulting relative error reaches the order of $10^{-7}$, demonstrating that PFEM preserves high accuracy while significantly accelerating convergence.

This subsection demonstrates that PFEM can effectively generalize across geometry, material properties, and boundary conditions for linear elasticity problems involving complex geometries, without requiring any special treatment beyond point-cloud inputs. Moreover, when combined with the warm-start stage, PFEM substantially reduces the number of iterations required by classical solvers, leading to significant computational acceleration.

\subsection{Nonlinear hyperelasticity}

\begin{figure}[t]
	\centering
	\includegraphics[scale=0.55]{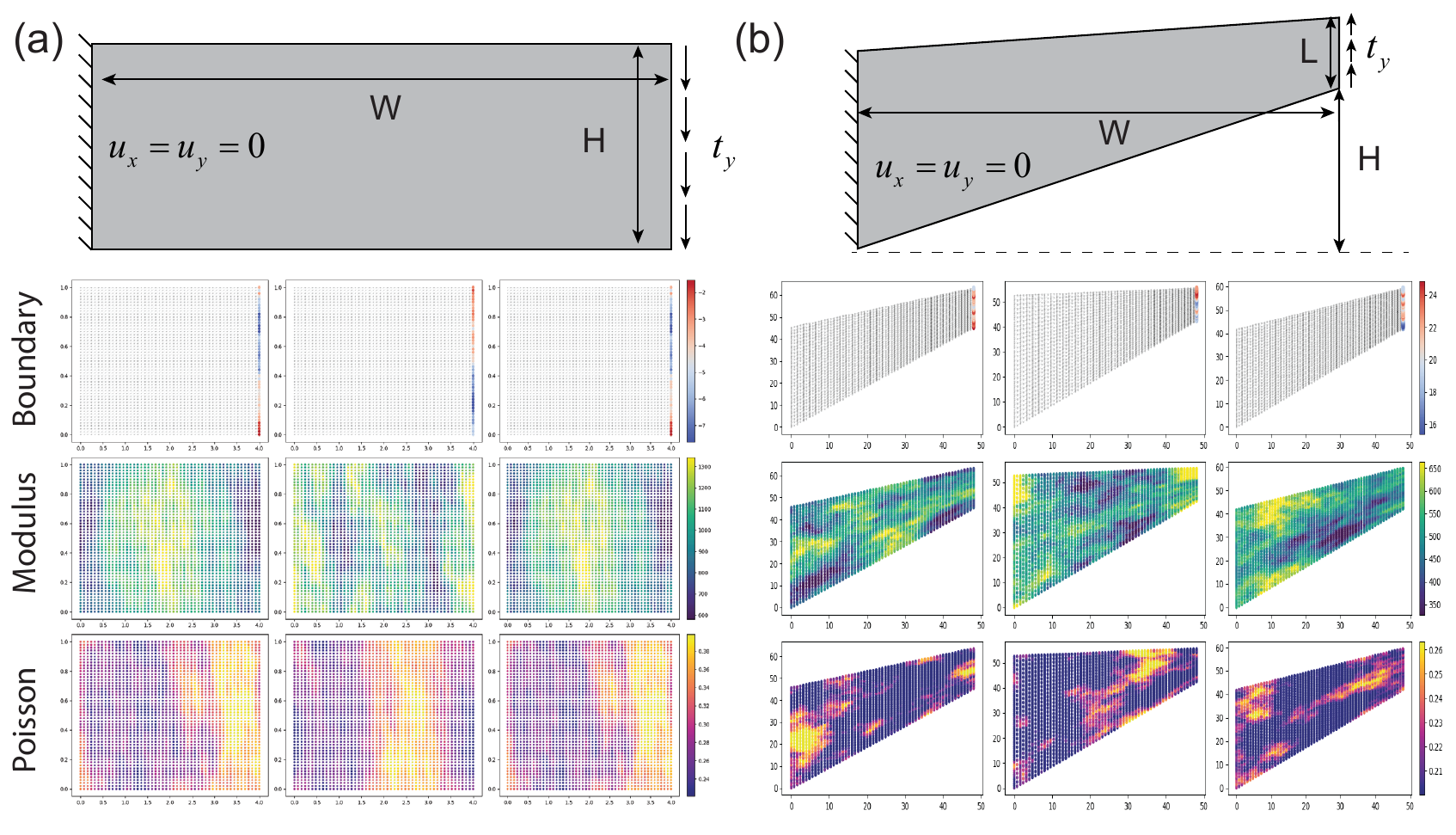}
	\caption{Nonlinear hyperelastic benchmarks: (a) hyperelastic beam with the left end fully clamped and a downward traction $t_y$ applied at the right end, with varying material properties and boundary conditions; (b) hyperelastic Cook’s membrane problem with varying geometry, material properties, and boundary conditions.}
	\label{fig:intro_hyper}
\end{figure}

The preceding experiments focus on linear elasticity problems. We next investigate the performance of PFEM on nonlinear problems. In general, owing to the strong nonlinear representation capability of neural operators, no fundamental distinction is expected when handling linear versus nonlinear partial differential equations. To validate this claim, we consider two classical nonlinear benchmarks in computational mechanics: a hyperelastic beam and Cook’s membrane problem, both formulated under plane strain conditions, as illustrated in \Cref{fig:intro_hyper}.

For hyperelasticity, the physics-informed loss function is defined as
\begin{equation}
	\begin{aligned}
		\mathcal{L}
		&=
		\int_{\Omega}
		\big(
		\Psi - \boldsymbol{f} \cdot \boldsymbol{u}
		\big)\, \mathrm{d}\Omega
		-
		\int_{\Gamma^{\boldsymbol{t}}}
		\bar{\boldsymbol{t}} \cdot \boldsymbol{u}\, \mathrm{d}\Gamma, \\
		\Psi
		&=
		\frac{1}{2}\lambda (\ln J)^2
		-
		\mu \ln J
		+
		\frac{1}{2}\mu (I_1 - 2), \\
		J
		&=
		\det(\boldsymbol{F}), \\
		I_1
		&=
		\mathrm{tr}(\boldsymbol{C}), \\
		\boldsymbol{C}
		&=
		\boldsymbol{F}^{\mathsf{T}} \boldsymbol{F}, \\
		\boldsymbol{u} & =\frac{\boldsymbol{x}}{L}*\mathcal{\boldsymbol{H}}_{\boldsymbol{\theta}}(\boldsymbol{x};\boldsymbol{G},\boldsymbol{\theta})
	\end{aligned}
\end{equation}
where $\Psi$ denotes the strain energy density of the Neo-Hookean material. The Lamé parameters $\lambda$ and $\mu$ are given by
\begin{equation}
	\begin{cases}
		\lambda = \dfrac{\nu E}{(1+\nu)(1-2\nu)}, \\[6pt]
		\mu = \dfrac{E}{2(1+\nu)},
	\end{cases}
\end{equation}
with $E$ and $\nu$ denoting Young’s modulus and Poisson’s ratio, respectively. The tensor $\boldsymbol{C}$ is the right Cauchy--Green deformation tensor, and $\boldsymbol{F}$ is the deformation gradient defined as
\begin{equation}
	\boldsymbol{F}
	=
	\dfrac{\partial \boldsymbol{x}}{\partial \boldsymbol{X}},
\end{equation}
where $\boldsymbol{x}$ and $\boldsymbol{X}$ denote the spatial and material coordinates, respectively. Their relationship is given by $\boldsymbol{x} = \boldsymbol{X} + \boldsymbol{u}$, with $\boldsymbol{u}$ being the displacement field.

The input of PFEM for the nonlinear hyperelastic problems consists of geometric information, material properties, and boundary conditions,
$
\boldsymbol{G}, \boldsymbol{M}, \boldsymbol{B}
=
\{ \boldsymbol{X}^{(i)} \}_{i=1}^{N},
$
where
$
\boldsymbol{X}^{(i)}
=
\{ x_j^{(i)}, y_j^{(i)}, E_j^{(i)}, \nu_j^{(i)}, t_{y,j}^{(i)} \}_{j=1}^{N^{(i)}}.
$
The output is the displacement field evaluated at the corresponding spatial locations,
$
\boldsymbol{U}
=
\{ \boldsymbol{u}^{(i)} \}_{i=1}^{N},
\qquad
\boldsymbol{u}^{(i)}
=
\{ u_j^{(i)}, v_j^{(i)} \}_{j=1}^{N^{(i)}}.
$

For the hyperelastic beam problem, the geometry is fixed with width $W=4$ and height $H=1$. The Young’s modulus, Poisson’s ratio, and traction boundary conditions are sampled from Gaussian random fields, as illustrated in \Cref{fig:intro_hyper}a. The geometry is kept unchanged for all samples, and a structured quadrilateral mesh is employed. In this setting, we compare the Transolver backbone adopted in PFEM with the Fourier Neural Operator (FNO) used in VINO \citep{eshaghi2025variational}, since FNO typically requires inputs defined on regular grids.

The dataset consists of 1000 samples, with 800 used for training and 200 for testing. All cases pass mesh convergence studies, ensuring the reliability of the FEM reference solutions. It is important to note that the training set does not contain labeled displacement fields; only material properties and boundary conditions are provided, and the neural operators are trained purely using governing physical equations. A uniform $50 \times 50$ mesh of linear quadrilateral elements (2500 elements in total) is used for all simulations.

The input of PFEM includes geometric information, material distributions, and boundary conditions,
$
\boldsymbol{G}, \boldsymbol{M}, \boldsymbol{B}
=
\{ \boldsymbol{X}^{(i)} \}_{i=1}^{N},
\quad
\boldsymbol{X}^{(i)}
=
\{ x_j^{(i)}, y_j^{(i)}, E_j^{(i)}, \nu_j^{(i)}, t_{y,j}^{(i)} \}_{j=1}^{N^{(i)}}.
$
The output is the displacement field at the corresponding spatial locations,
$
\boldsymbol{U}
=
\{ u_j^{(i)}, v_j^{(i)} \}_{i=1}^{N}.
$
The boundary condition is specified by applying a Gaussian-random-field traction $t_y$ on the right boundary, while all other boundaries are traction-free.

\begin{figure}[t]
	\centering
	\includegraphics[scale=0.55]{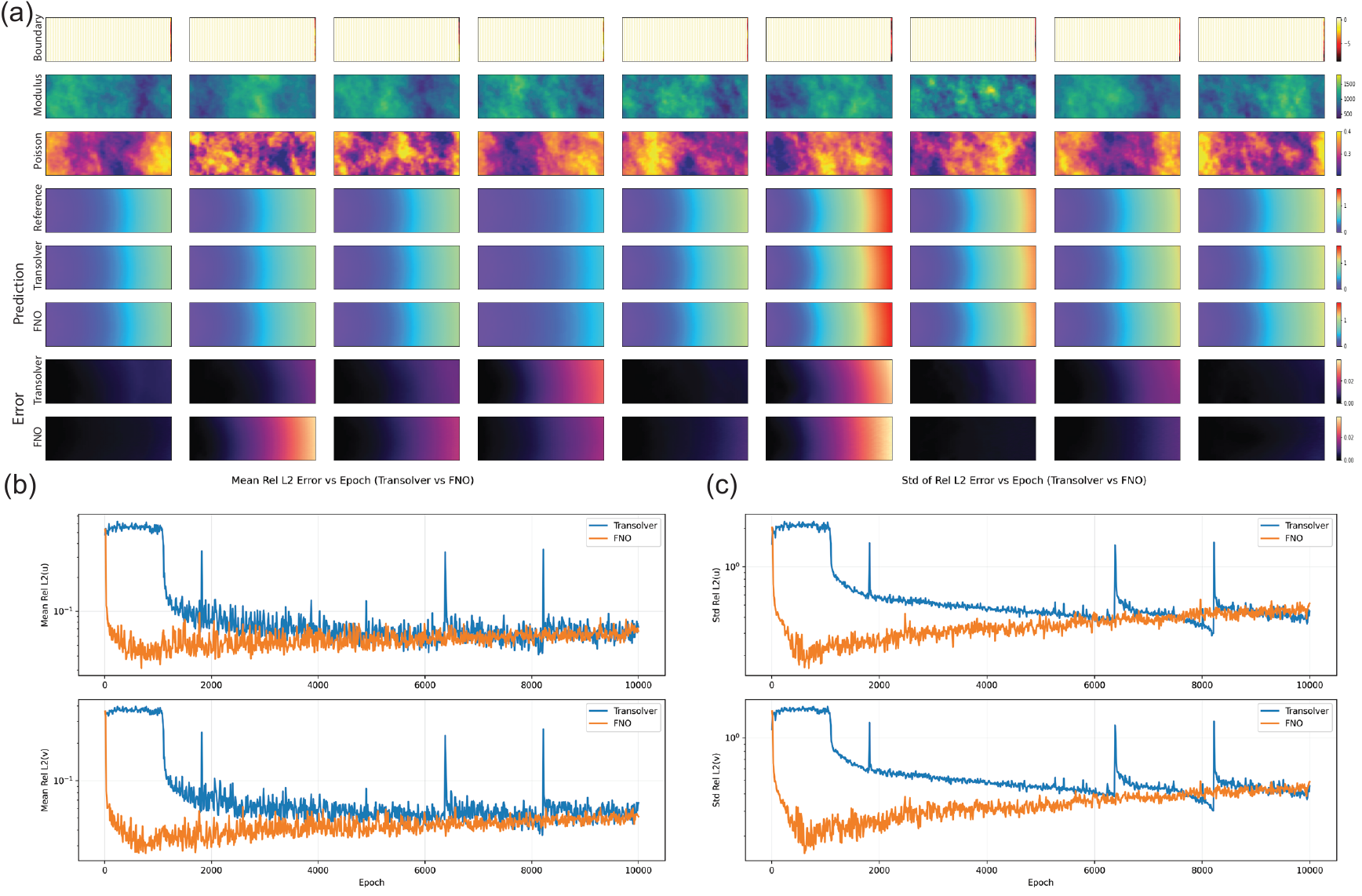}
	\caption{Comparison between Transolver (PFEM) and FNO for the hyperelastic beam problem: (a) from top to bottom: boundary condition, Young’s modulus field, Poisson’s ratio field, FEM reference solution, PFEM prediction, FNO prediction, PFEM error contour, and FNO error contour; (b) evolution of the mean relative error during training; (c) evolution of the standard deviation of the relative error during training.}
	\label{fig:FNO_transolver_compare_beam}
\end{figure}

\Cref{fig:FNO_transolver_compare_beam} presents a comparison between Transolver and FNO for the hyperelastic beam problem. Although the convergence speed of Transolver is slightly slower than that of FNO during training, both approaches achieve comparable relative errors after convergence, with Transolver occasionally exhibiting slightly better accuracy. The fluctuations observed in the error mean and standard deviation arise from the learning-rate decay strategy, where the learning rate is reduced by a factor of $0.9$ every 1000 training steps.

It is worth emphasizing that Transolver offers significantly greater flexibility in input representation and geometric handling, as it does not rely on structured grid inputs. Consequently, Transolver is particularly well suited for problems involving complex geometries, multi-scale discretizations, and unstructured meshes. This characteristic highlights its strong potential for future research on complex engineering structures and realistic geometric configurations.

We next present the warm-start results of PFEM for the nonlinear hyperelastic beam problem. We record the number of Newton iterations required by PFEM and by the conventional FEM solver over 150 test cases. In PFEM, the prediction provided by the Transolver backbone is used as the initial guess for the Newton iteration, whereas in the FEM baseline a zero initial vector is employed.

For nonlinear problems, the standard FEM procedure typically applies the external load incrementally and performs a Newton iteration at each load step; a detailed description is provided in \ref{sec:Iterative-FEM}. Since PFEM is capable of directly predicting the displacement field corresponding to the final loading state, the load step can be set to one, and the PFEM prediction can be directly used as the initial guess for the final Newton iteration. To ensure a fair comparison, both PFEM and FEM are evaluated using a single load step in the following experiments.

\begin{figure}[t]
	\centering
	\includegraphics[scale=0.55]{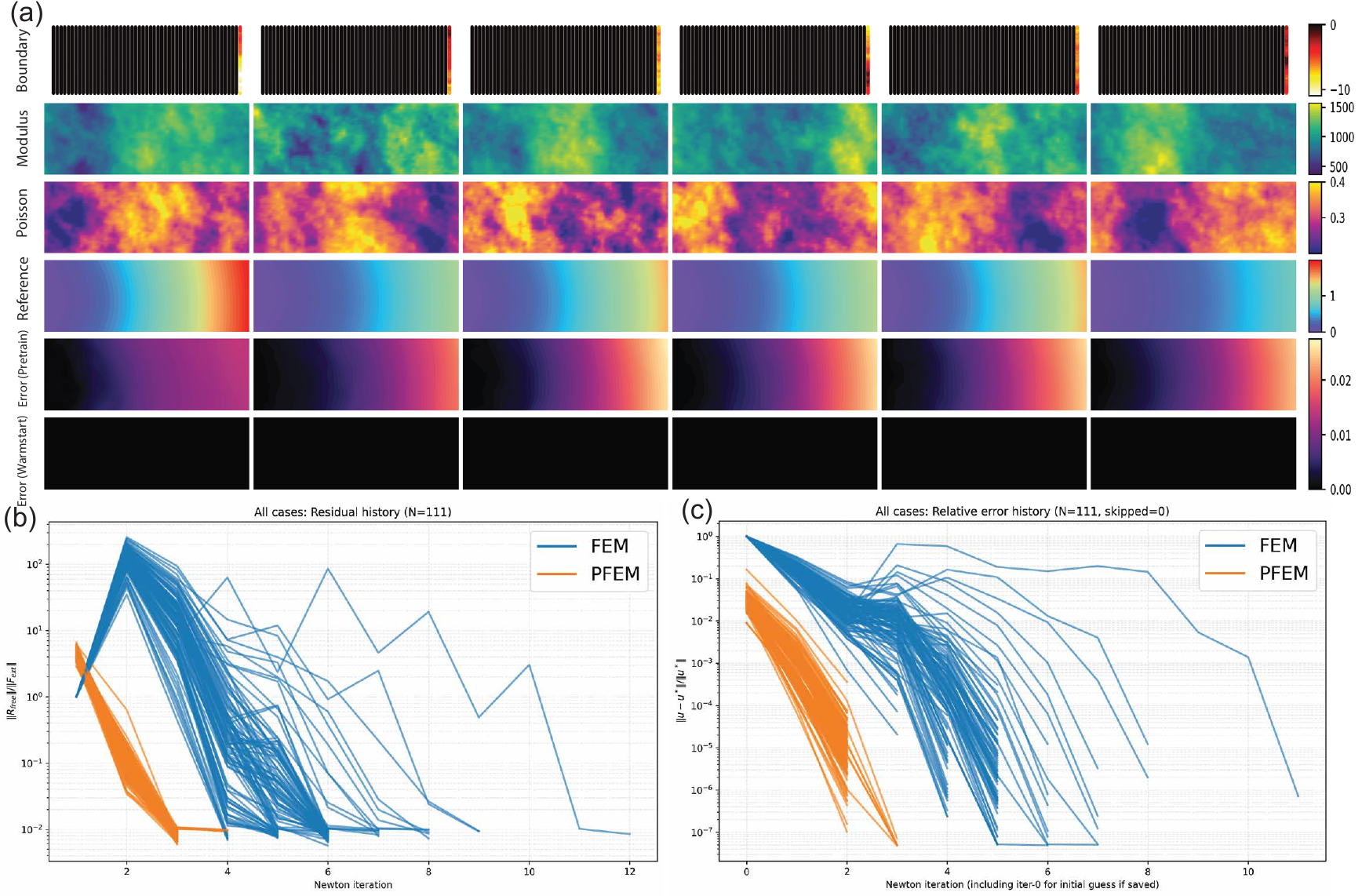}
	\caption{Warm-start results of PFEM for the hyperelastic beam problem: (a) comparison of error contours between the pretraining and warm-start stages; (b) convergence history of the Newton residual for PFEM and FEM; (c) convergence history of the relative error for PFEM and FEM.}
	\label{fig:PFEM_hyper_beam_contourf_geo_mater_boundary_warmstart}
\end{figure}

The warm-start results for the hyperelastic beam with varying material properties and boundary conditions are shown in \Cref{fig:PFEM_hyper_beam_contourf_geo_mater_boundary_warmstart}. PFEM leads to a substantial reduction in the number of Newton iterations compared to the conventional FEM solver, as also summarized in \Cref{tab:PFEM_arch_diff_problem}. The convergence history shown in \Cref{fig:PFEM_hyper_beam_contourf_geo_mater_boundary_warmstart}b is plotted in terms of the Newton residual
$
\|\boldsymbol{r}\| = \|\boldsymbol{f}^{\mathrm{int}}(\boldsymbol{U}) - \boldsymbol{f}^{\mathrm{ext}}(\boldsymbol{U})\|.
$
In \Cref{fig:PFEM_hyper_beam_contourf_geo_mater_boundary_warmstart}c, iteration index zero corresponds to the relative error of the initial guess provided by the Transolver prediction.

It is worth noting that the initial Newton residual associated with PFEM can be larger than that of FEM, while the corresponding relative error with respect to the reference solution is significantly smaller. This observation highlights an important distinction in nonlinear iterative solvers: a smaller Newton residual does not necessarily imply a smaller relative solution error. Despite this difference in the initial residual behavior, PFEM rapidly enters the local convergence regime of Newton’s method, leading to a pronounced reduction in the total number of iterations required to reach convergence.

We further investigate the potential of PFEM on the Cook’s membrane problem with irregular geometries, which serves as a classical nonlinear benchmark in computational mechanics. It is evident that handling such problems using FNO-based approaches is inconvenient due to the requirement of structured grid inputs. In contrast, PFEM naturally accommodates irregular geometries through point-cloud representations.

The dataset consists of 1000 samples, among which 800 are used for training and 200 for testing. All cases pass mesh convergence studies, ensuring the reliability of the FEM reference solutions. The input of PFEM includes geometric information, material distributions, and boundary conditions,
$
\boldsymbol{G}, \boldsymbol{M}, \boldsymbol{B}
=
\{ \boldsymbol{X}^{(i)} \}_{i=1}^{N},
\qquad
\boldsymbol{X}^{(i)}
=
\{ x_j^{(i)}, y_j^{(i)}, E_j^{(i)}, \nu_j^{(i)}, t_{y,j}^{(i)} \}_{j=1}^{N^{(i)}}.
$
The output is the displacement field evaluated at the corresponding spatial locations,
$
\boldsymbol{U}
=
\{ u_j^{(i)}, v_j^{(i)} \}_{i=1}^{N}.
$

The boundary condition is defined by prescribing a Gaussian-random-field traction $t_y$ on the right boundary, while all other boundaries are traction-free. The geometry is varied by sampling the clamped length on the left boundary as well as the length and position of the traction boundary on the right side from Gaussian random fields. These geometric variations are illustrated in \Cref{fig:intro_hyper}b. Notably, the training set does not contain any labeled displacement fields; only geometry, material properties, and boundary conditions are provided, and the neural operator is trained purely using governing physical equations.

\begin{figure}[t]
	\centering
	\includegraphics[scale=0.55]{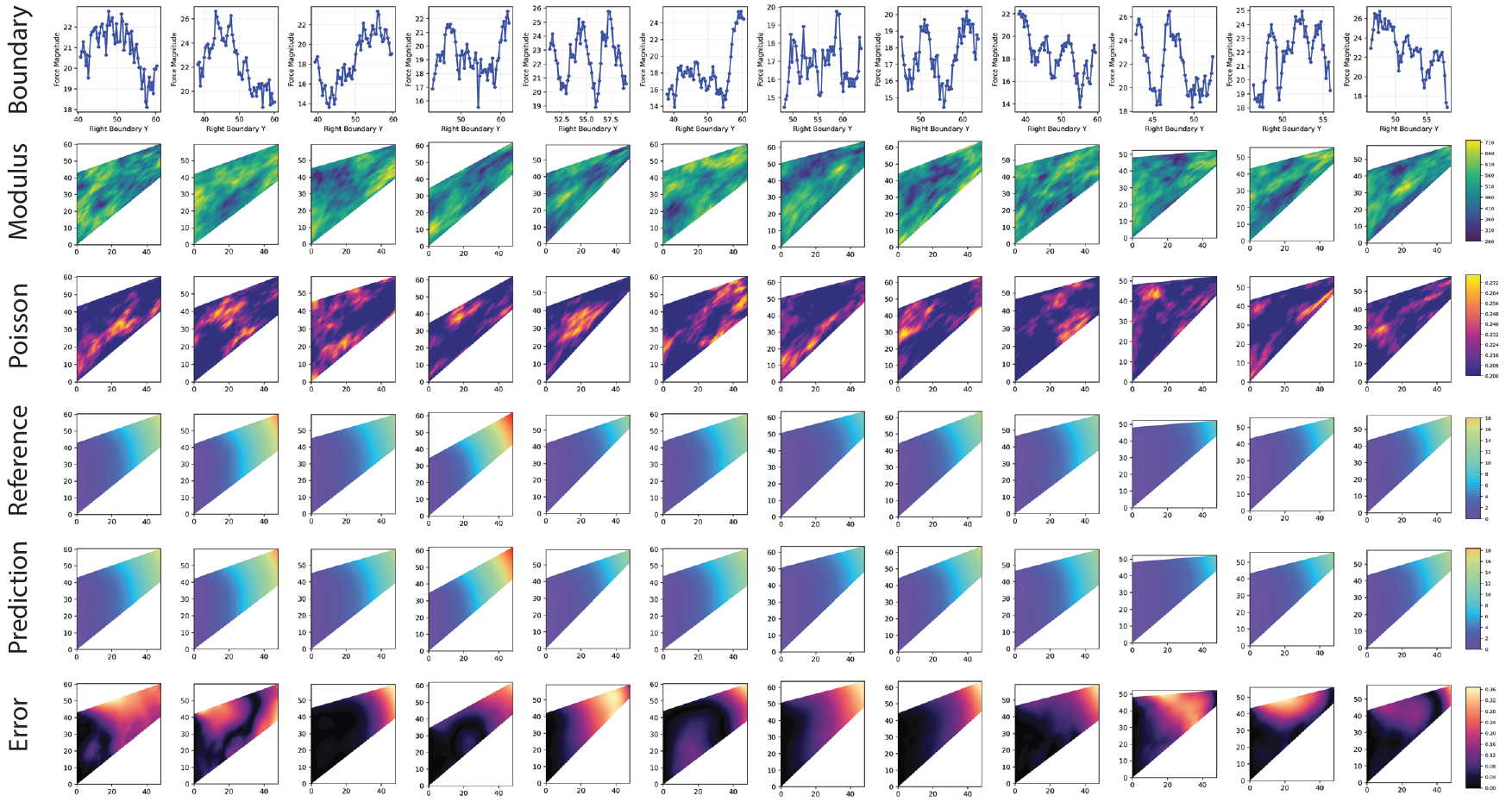}
	\caption{Predicted absolute displacement magnitude by PFEM for the Cook’s membrane problem with varying geometries, material distributions, and boundary conditions. From top to bottom: boundary condition, Young’s modulus field, Poisson’s ratio field, FEM reference solution, PFEM prediction, and PFEM error contour. In the boundary-condition plot, the horizontal axis denotes the $y$-coordinate on the right boundary of the Cook’s membrane, and the vertical axis represents the corresponding traction $t_y$.}
	\label{fig:PFEM_Hyper_cook}
\end{figure}

A uniform $30 \times 30$ mesh consisting of 900 quadratic Q8 serendipity quadrilateral elements is employed. The pretraining results of PFEM are shown in \Cref{fig:PFEM_Hyper_cook}. PFEM demonstrates strong generalization capability across irregular geometries, heterogeneous material distributions, and varying boundary conditions, while requiring only point-cloud inputs. This highlights the flexibility and robustness of PFEM for complex nonlinear problems involving realistic geometries.

We next present the warm-start results of PFEM for the Cook’s membrane problem. We record the number of Newton iterations required by PFEM and by the conventional FEM solver over 150 test cases. In PFEM, the Transolver prediction is used as the initial guess for the Newton iteration, whereas in the FEM baseline a zero initial vector is employed.

\begin{figure}[t]
	\centering
	\includegraphics[scale=0.55]{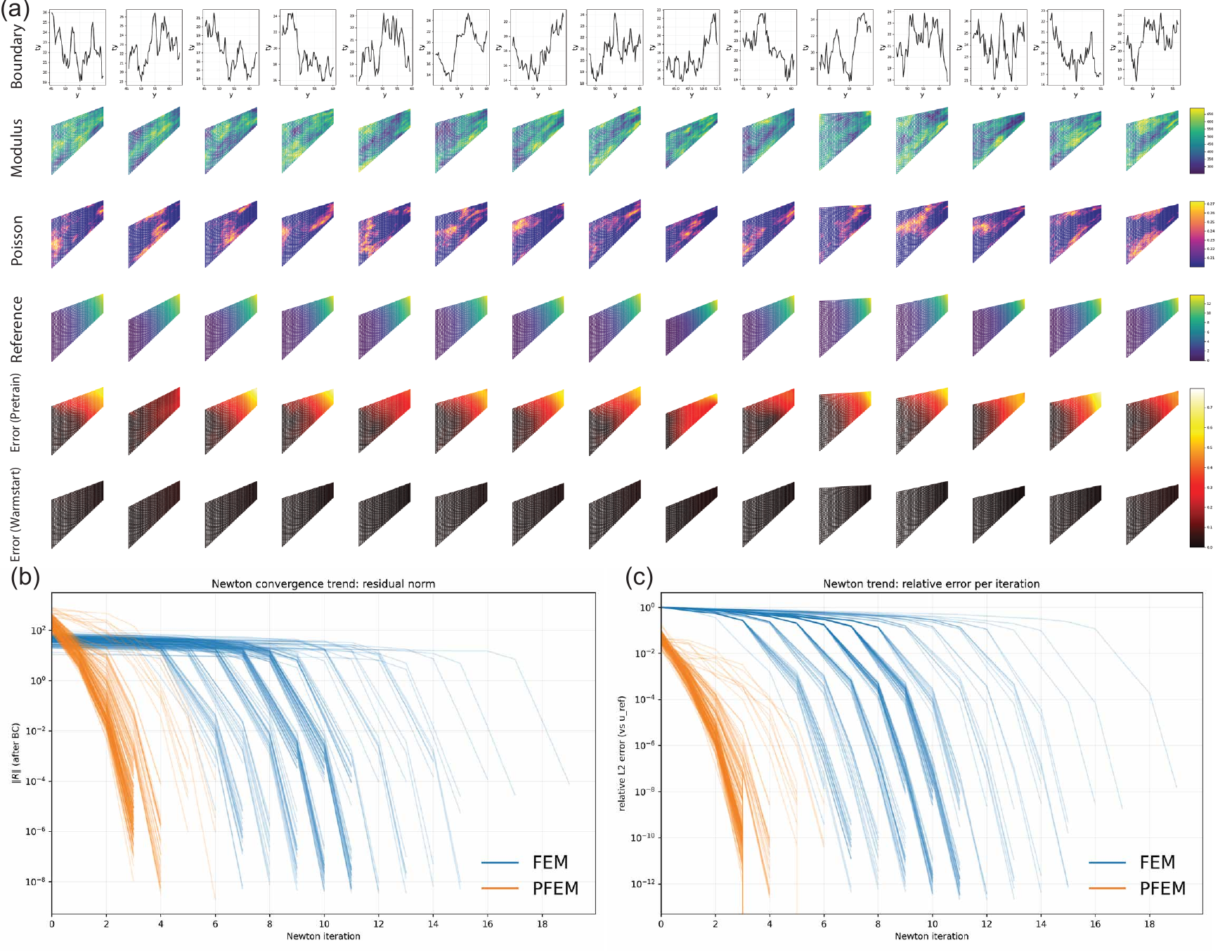}
	\caption{Warm-start results of PFEM for the Cook’s membrane problem: (a) comparison of error contours between the pretraining and warm-start stages; (b) convergence histories of the Newton residual for PFEM and FEM; (c) convergence histories of the relative error for PFEM and FEM.}
	\label{fig:PFEM_cook_contourf_geo_mater_boundary_warmstart}
\end{figure}

The warm-start results under simultaneous variations in geometry, material properties, and boundary conditions are shown in \Cref{fig:PFEM_cook_contourf_geo_mater_boundary_warmstart}. PFEM again leads to a substantial reduction in the number of Newton iterations compared to the traditional FEM solver, as also summarized in \Cref{tab:PFEM_arch_diff_problem}. These results demonstrate that PFEM effectively accelerates nonlinear solvers by providing high-quality initial guesses, even for problems involving irregular geometries and strong nonlinearities.

This subsection demonstrates that PFEM can effectively generalize across geometry, material properties, and boundary conditions for nonlinear problems involving complex and irregular geometries, without requiring any special treatment beyond point-cloud inputs. Moreover, when combined with the warm-start stage, PFEM significantly reduces the number of Newton iterations required to solve nonlinear problems, leading to substantial computational acceleration.

\subsection{Three-dimensional solid mechanics homogenization}

Homogenization plays a fundamental role in multiscale modeling in solid mechanics
\citep{hashin1983analysis,guedes1990preprocessing,hassani1998review,nguyen2025universal}, 
and is closely related to the concept of the representative volume element (RVE). 
Homogenization aims to determine the effective macroscopic material properties by solving boundary value problems, typically subject to periodic boundary conditions, defined on a RVE.
Therefore, homogenization serves as an essential tool for studying multiscale physical phenomena.

Traditional numerical homogenization methods rely heavily on finite element analysis and often require substantial computational resources, particularly for problems involving complex geometries, heterogeneous materials, and high-resolution discretizations. 
Recently, operator learning approaches have attracted increasing attention for modeling multiscale mechanical behavior in solid mechanics 
\citep{wang2026pretraining,nguyen2025universal,harandi2025spifol}. 
The pretraining stage of PFEM employs the Transolver model, which belongs to the class of operator learning. Therefore, in this subsection we investigate the performance of PFEM for three-dimensional homogenization problems in solid mechanics.

Various homogenization frameworks exist in the literature. 
Here, we adopt the classical numerical homogenization approach proposed by Andreassen et al.~\citep{andreassen2014determine}, which aims to compute the effective macroscopic elastic tensor of periodic composite materials. 
The detailed homogenization formulation is in \citep{wang2026pretraining}. 
In \citep{wang2026pretraining}, the Fourier neural operator (FNO) was trained in a data-driven manner. 
However, our manuscript trains the Transolver purely using physics-based constraints without data, following the PFEM framework proposed in this study.

The loss function is derived from the minimum potential energy principle of linear elasticity:

\begin{equation}
	\begin{aligned}
		\mathcal{L} &= \frac{1}{6}\sum_{IJ}^{6}\Pi^{(IJ)},\\
		\Pi^{(IJ)} &= \int_{\Omega}\frac{1}{2}
		E_{klpq}\varepsilon_{kl}(\boldsymbol{X}^{(v)(IJ)})
		\varepsilon_{pq}(\boldsymbol{X}^{(v)(IJ)})\, d\Omega
		-\int_{\Omega}E_{klpq}\varepsilon_{kl}^{(0)(IJ)}\, d\Omega ,
	\end{aligned}
	\label{eq:homo_loss}
\end{equation}
where $\boldsymbol{X}^{(v)}$ denotes the fluctuation displacement field, which is the target output of the neural operator. 
$\boldsymbol{\varepsilon}^{(0)(IJ)}$ represents the six prescribed macroscopic loading cases used in the homogenization procedure, corresponding to the six unit strain fields. 
Further details of the homogenization formulation are provided in \ref{sec:homo}.

\begin{figure}
	\begin{centering}
		\includegraphics[scale=0.45]{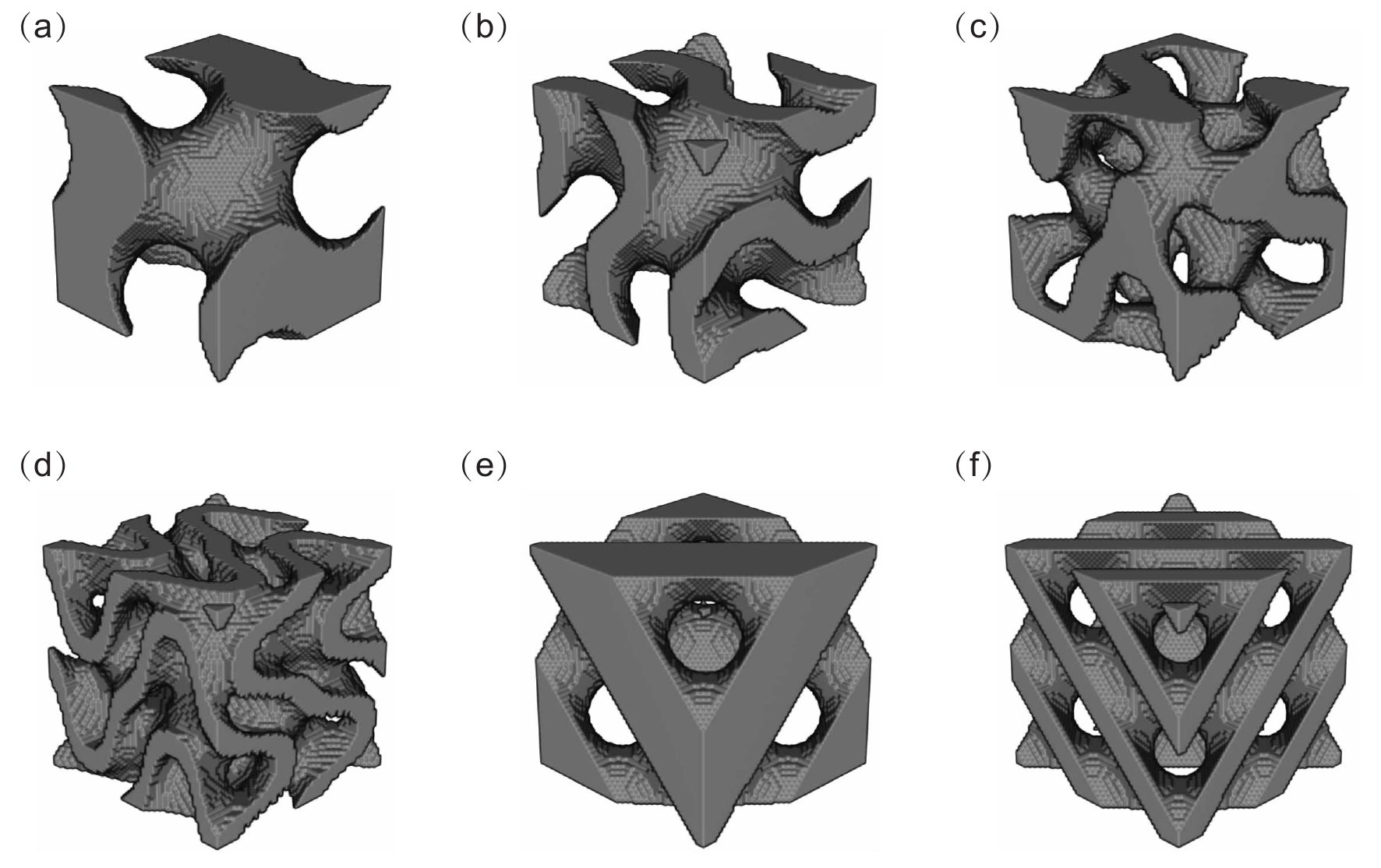}
		\par\end{centering}
	\caption{Six different TPMS structures:
		(a) Schoen Gyroid: ``solid networks''; 
		(b) Schoen Gyroid: ``sheet networks''; 
		(c) Schwarz Diamond: ``solid networks''; 
		(d) Schwarz Diamond: ``sheet networks''; 
		(e) Fischer Koch S: ``solid networks''; 
		(f) Fischer Koch S: ``sheet networks''. 
		\label{fig:6tpms}}
\end{figure}

The dataset consists of three-dimensional TPMS structures with a resolution of $64 \times 64 \times 64$. 
Six types of TPMS geometries are considered: 
Schoen Gyroid (``solid networks'' and ``sheet networks''), 
Schwarz Diamond (``solid networks'' and ``sheet networks''), 
and Fischer Koch S (``solid networks'' and ``sheet networks''), 
as illustrated in \Cref{fig:6tpms}. 
For each TPMS type, 600 structures are generated, resulting in a total of 3600 samples. 
The total dataset size is approximately 103\,GB and is publicly available at
\url{https://1drv.ms/f/c/f999b24db97d074b/UgBLB325TbKZIID5YoUEAAAAAK57y3r7MVu56Ds}.

Both the geometry and the material distributions are randomly generated. 
The geometric volume fractions range from 0.26 to 0.66, while Poisson's ratios vary between 0.1 and 0.4. 
Details of the geometric volume fraction generation are described in \ref{sec:TPMS}. 
Both parameters follow uniform distributions. 
Since the present study focuses on linear elasticity homogenization, Young's modulus is set to unity during training, which simplifies the computation of the effective elastic tensor through linear scaling. 
Among the 3600 samples, 3000 are used for training and 600 for testing. 
It is worth noting that the $64^3$ resolution grid represents the entire cubic domain, including both material and void regions. However, when training the Transolver, only the point cloud corresponding to material points is used. As a result, the number of input points varies across different samples.

\begin{figure}
	\begin{centering}
		\includegraphics[scale=0.36]{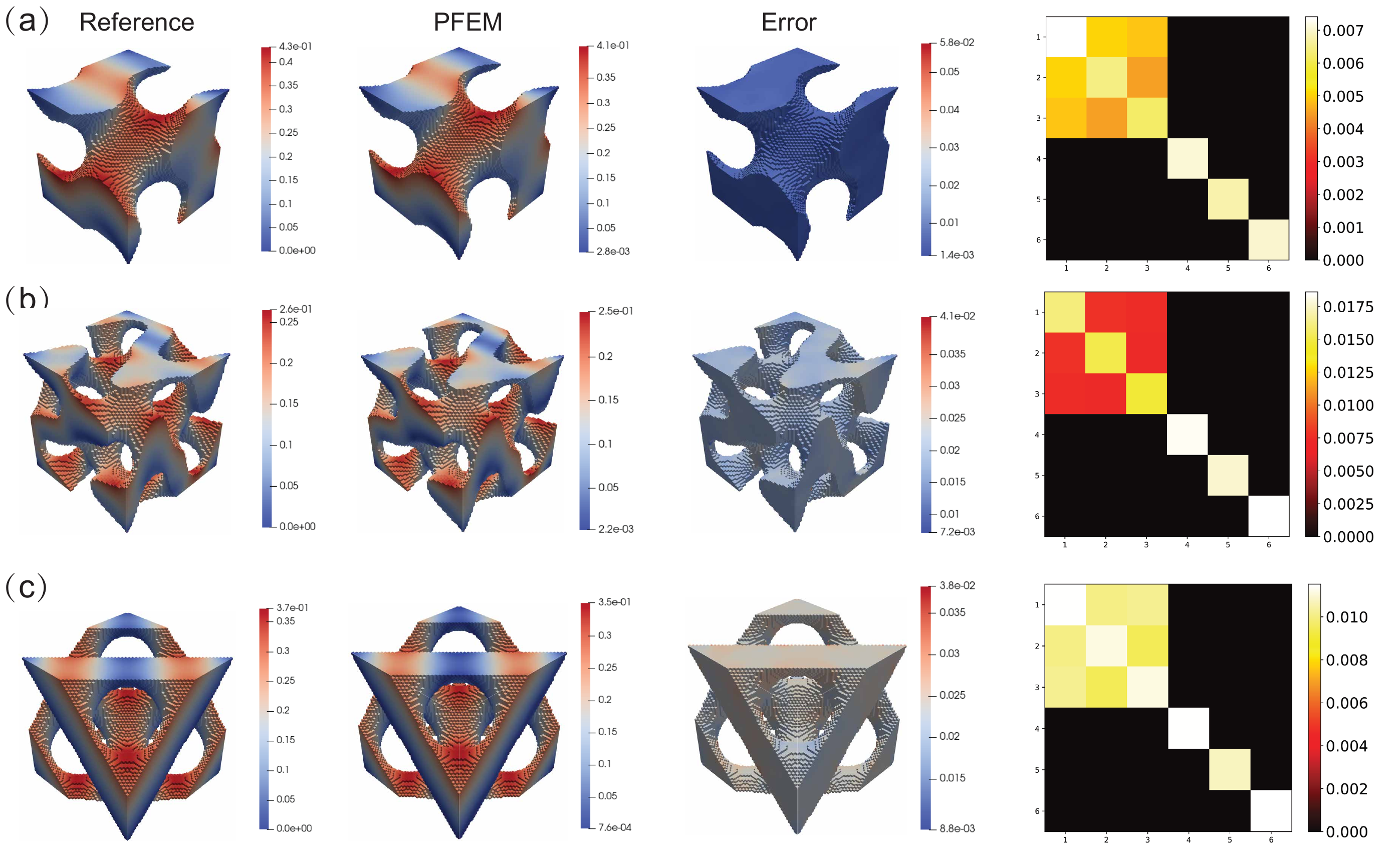}
		\par\end{centering}
	\caption{Absolute displacement predictions of PFEM in the pretraining stage for TPMS ``solid network'' structures:
		(a) Schoen Gyroid, 
		(b) Schwarz Diamond, 
		(c) Fischer Koch S. 
		The first column shows the FEM reference solution, the second column shows PFEM predictions, the third column presents the error fields, and the last column reports the relative error of the effective elastic tensor.
		\label{fig:tpms_solid}}
\end{figure}

\begin{figure}
	\begin{centering}
		\includegraphics[scale=0.36]{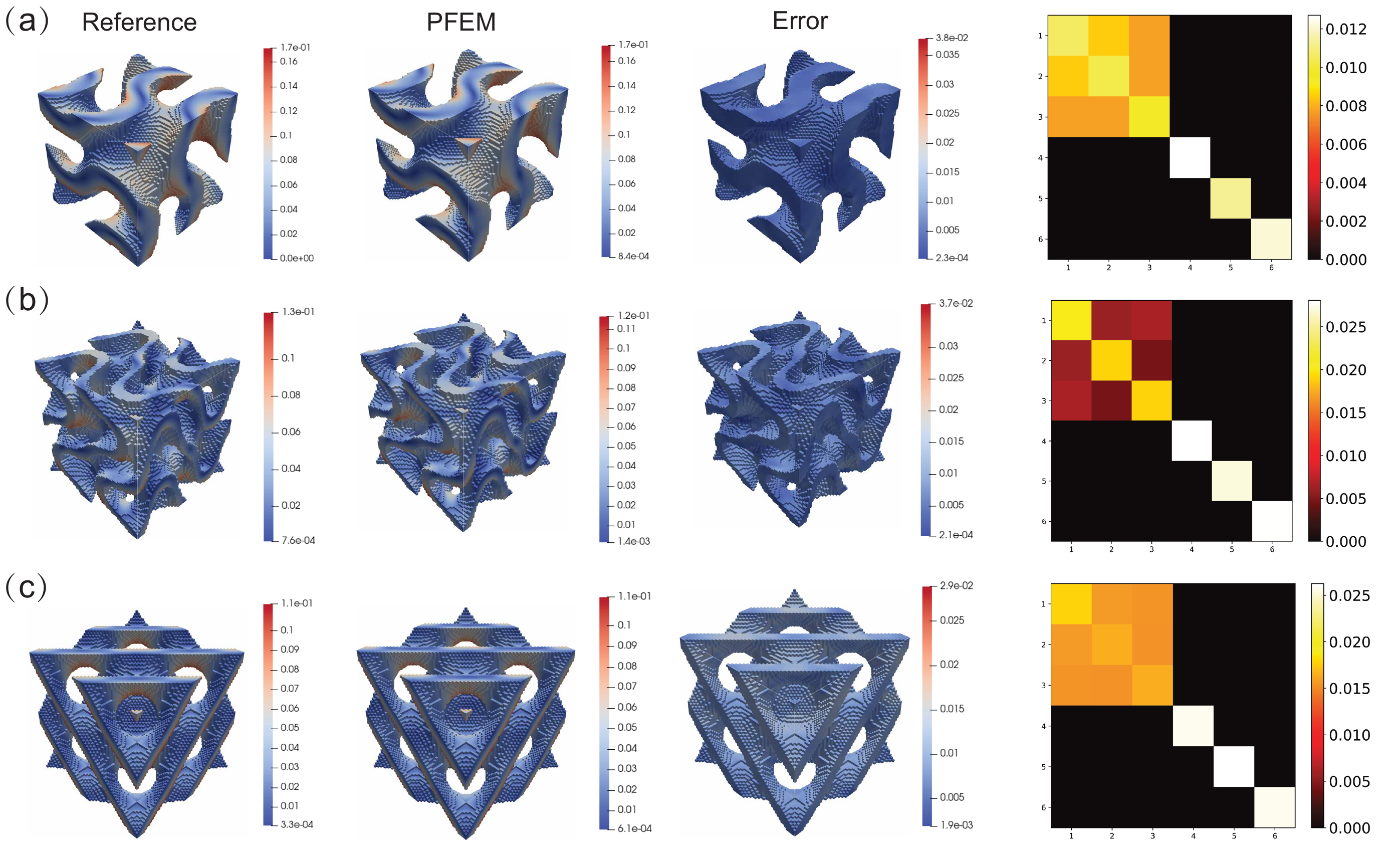}
		\par\end{centering}
	\caption{Absolute displacement predictions of PFEM in the pretraining stage for TPMS ``sheet network'' structures:
		(a) Schoen Gyroid, 
		(b) Schwarz Diamond, 
		(c) Fischer Koch S. 
		The columns correspond to FEM reference solutions, PFEM predictions, error fields, and relative errors of the effective elastic tensor, respectively.
		\label{fig:tpms_sheet}}
\end{figure}

The architectural details and training configurations of the Transolver are summarized in \Cref{tab:PFEM_arch_diff_problem}. 
\Cref{fig:tpms_solid} and \Cref{fig:tpms_sheet} present the prediction results of PFEM, where the finite element solution serves as the reference. 
The FEM simulations employ $64 \times 64 \times 64 = 262{,}144$ linear elements, and the average computational time for a single FEM sample is approximately 150\,s.

The results show that PFEM in the pretraining stage predicts displacement fields that agree closely with the FEM reference solutions, with overall relative errors below $10\%$, as summarized in \Cref{tab:PFEM_arch_diff_problem}. 
For clarity, only the absolute displacement field corresponding to the first loading case is shown. 
In practice, six loading cases are required in the homogenization procedure, as detailed in \ref{sec:homo}.

The input to the Transolver consists of four features: the three spatial coordinates and the Poisson's ratio. 
The output dimension is 18, corresponding to three displacement components under six loading cases. 
Notably, PFEM only requires the coordinates of material points in the TPMS structures as input. 
This differs from FNO-based approaches, which typically require dense 3D grid representations where void regions are filled with zeros \citep{wang2026pretraining}. 
Consequently, PFEM benefits from the point-cloud-based representation of Transolver, which provides a natural advantage for handling complex geometries.
Specifically, the input of PFEM consists of the geometry and material information $\boldsymbol{G},\boldsymbol{M}=\{\boldsymbol{X}^{(i)}\}_{i=1}^{N}$, where $\boldsymbol{X}^{(i)}=\{x_{j}^{(i)},y_{j}^{(i)},z_{j}^{(i)},\nu_{j}^{(i)}\}_{j}^{N^{(i)}}$. The output is the displacement at the corresponding coordinate points $\boldsymbol{U}=\{\boldsymbol{u}^{(i)(1)},\boldsymbol{u}^{(i)(2)},\boldsymbol{u}^{(i)(3)},\boldsymbol{u}^{(i)(4)},\boldsymbol{u}^{(i)(5)},\boldsymbol{u}^{(i)(6)}\}_{i=1}^{N}$, where $\boldsymbol{u}^{(i)(M)}=\{u_{j}^{(i)(M)},v_{j}^{(i)(M)},w_{j}^{(i)(M)}\}_{j}^{N^{(i)}}$, and $M=1,\cdots,6$ corresponds to the six loading cases in the homogenization process.

After obtaining the displacement fields, the effective elastic tensor of size $6 \times 6$ is computed, as shown in the last column of \Cref{fig:tpms_solid} and \Cref{fig:tpms_sheet}. 
PFEM in the pretraining stage demonstrates excellent predictive capability for the effective elastic tensor, with maximum relative errors below $3\%$. 
Meanwhile, PFEM achieves significantly higher computational efficiency than conventional FEM methods, as summarized in \Cref{tab:PFEM_arch_diff_problem}.

\begin{figure}
	\begin{centering}
		\includegraphics[scale=0.55]{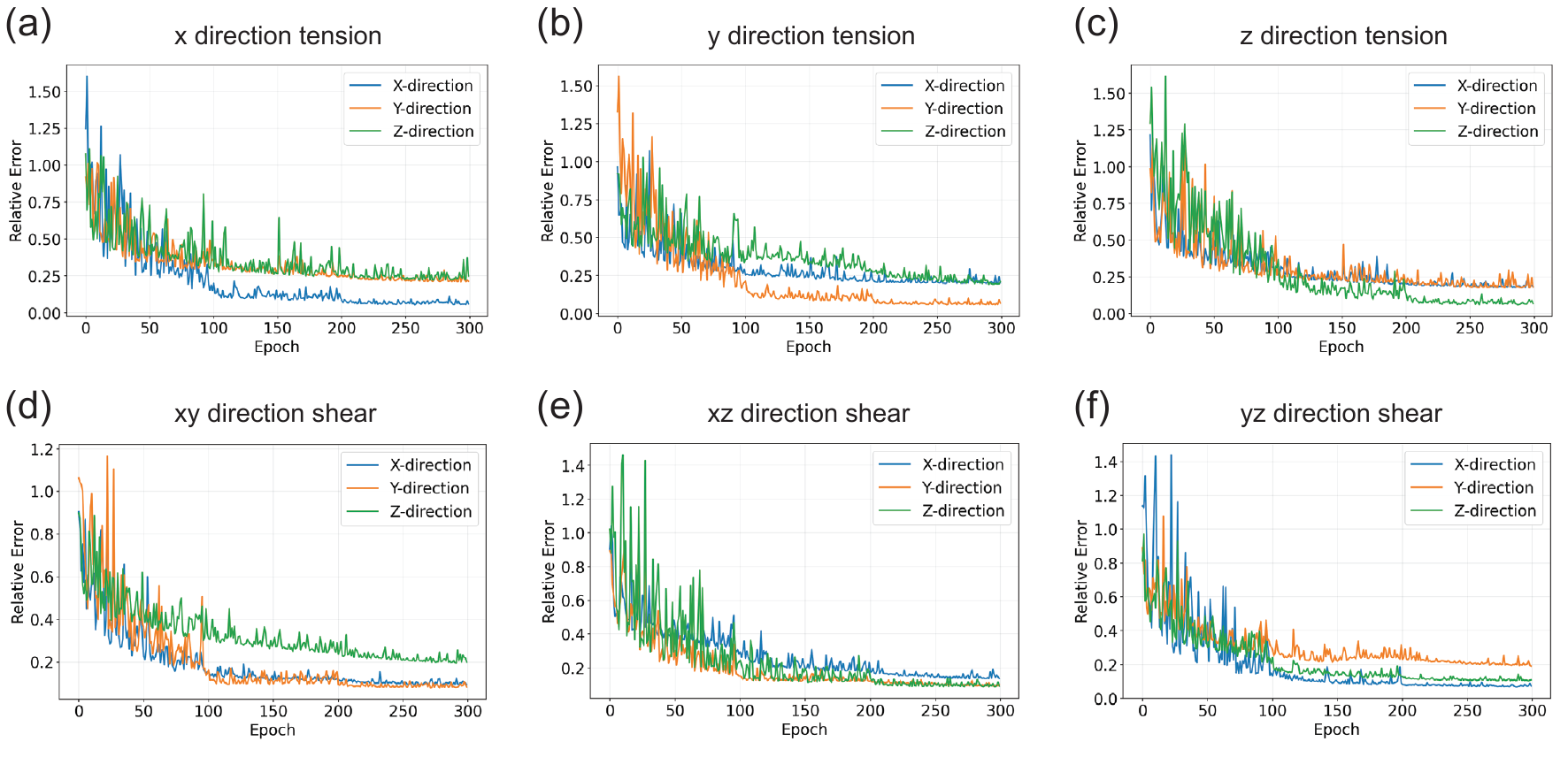}
		\par\end{centering}
	\caption{Testing errors of PFEM in the pretraining stage.
		(a--f) correspond to tensile loading in the $x$, $y$, and $z$ directions, and shear loading in the $xy$, $xz$, and $yz$ directions. 
		For each loading case, the relative errors of the three displacement components are reported.
		\label{fig:PFEM_error}}
\end{figure}

The testing errors of PFEM in the pretraining stage are shown in \Cref{fig:PFEM_error}, where the training process spans 300 epochs. 
The smallest errors are observed in the displacement component aligned with the principal loading direction. 
This behavior arises because the energy-based loss function assigns different sensitivities to displacement components under different loading cases. 
For example, under tensile loading in the $x$ direction, as shown in \Cref{fig:PFEM_error}a, the displacement error in the $x$ direction is smaller than those in the $y$ and $z$ directions. This is because the $x$-direction displacement is dominant under $x$-direction tensile loading, resulting in larger displacement magnitudes compared with the $y$ and $z$ directions. Consequently, reducing the error in the $x$-direction displacement contributes more significantly to the reduction of the total energy. As a result, during the optimization process PFEM tends to focus more on improving the accuracy of the $x$-direction displacement field.

\begin{figure}
	\begin{centering}
		\includegraphics[scale=0.55]{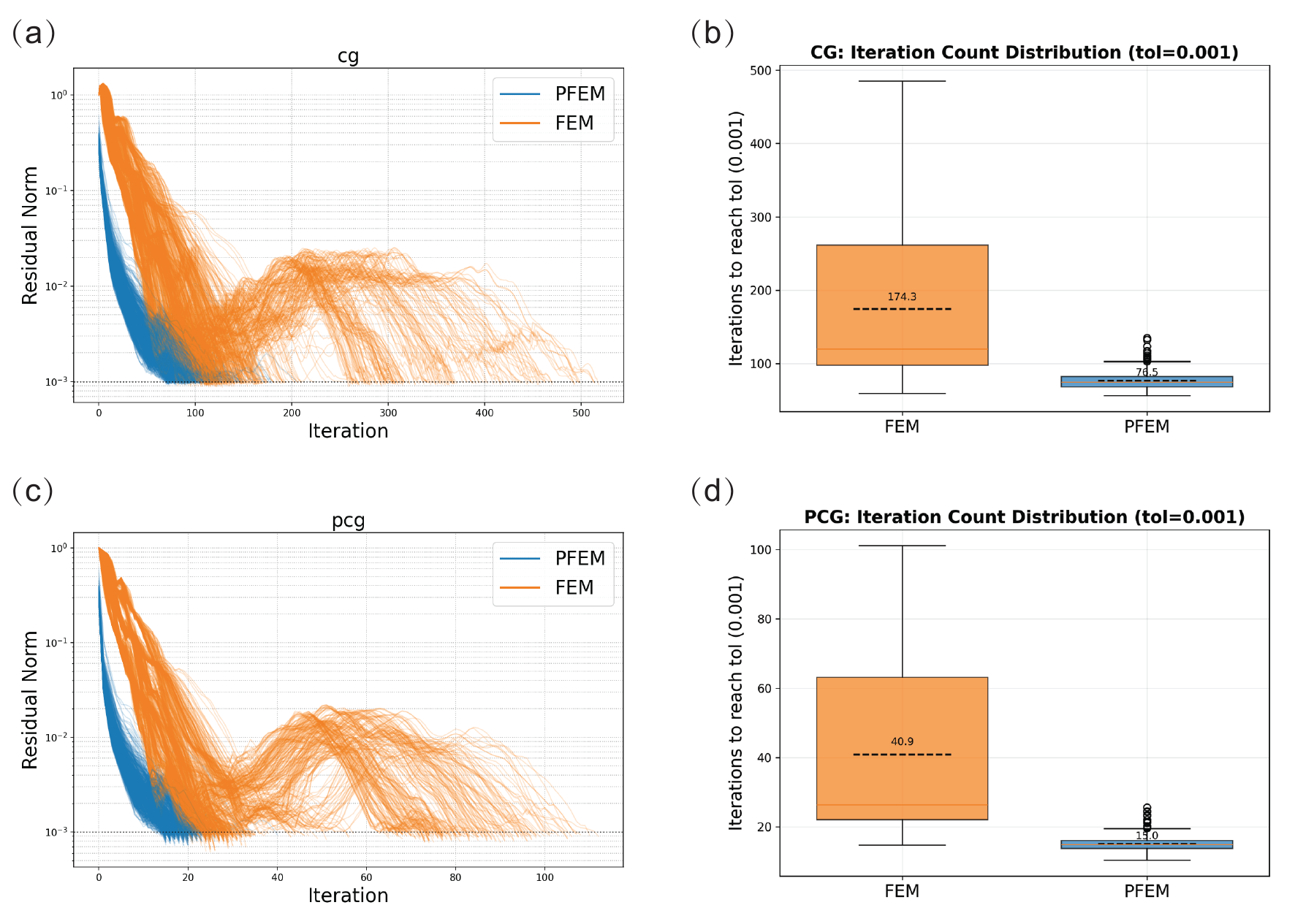}
		\par\end{centering}
	\caption{Warm-start results of PFEM for the 3D homogenization problem. 
		The first and second rows correspond to warm-start iterations using the conjugate gradient (CG) and preconditioned conjugate gradient (PCG) methods, respectively. 
		The first column shows the convergence histories of PFEM and FEM, while the second column compares the number of iterations required by PFEM and FEM.
		\label{fig:homo_warm_start}}
\end{figure}

Next, we investigate the performance gain of PFEM in the warm-start stage. 
\Cref{fig:homo_warm_start} presents the warm-start results for the three-dimensional homogenization problem. 
In the conventional FEM baseline, the zero vector is used as the initial guess. 
In contrast, PFEM uses the displacement field predicted by the Transolver as the initial guess for the iterative solver.

It can be clearly observed that PFEM significantly reduces the number of iterations required for convergence, regardless of whether the iterative solver in the warm-start stage is the conjugate gradient (CG) method or the preconditioned conjugate gradient (PCG) method. 
Compared with standard FEM, PFEM exhibits a more robust and smoother iterative convergence process. 
As shown in \Cref{fig:homo_warm_start}(a) and (c), the convergence rate is determined solely by the chosen iterative algorithm. 
Therefore, when the same iterative algorithm is used for both PFEM and FEM, the convergence rates remain identical. 
The reduction in iteration numbers achieved by PFEM originates entirely from the improved initial guess.

Overall, this subsection demonstrates that PFEM generalizes effectively to complex three-dimensional homogenization problems involving diverse geometries and material distributions. 
Notably, no special preprocessing is required, and the model only takes point-cloud inputs corresponding to material points. 
Furthermore, when combined with the warm-start stage, PFEM substantially reduces the number of iterations required to solve large-scale 3D homogenization problems.

\section{Discussion}

\subsection{Out-of-Distribution Generalization of the PFEM Pretraining Stage}

In the results presented above, the training and test datasets are sampled from the same underlying probability distributions. A natural and more challenging question is how the pretraining stage of PFEM performs when evaluated under distribution shifts.

To this end, we consider a test case with a geometric configuration that is not included in the training distribution, where the circular hole is located at $(2,2)$ with a radius of $0.5$. Moreover, the elastic modulus field, Poisson's ratio field, and boundary conditions are all generated from probability distributions that differ from those used during training.

\begin{figure} \begin{centering} \includegraphics[scale=0.70]{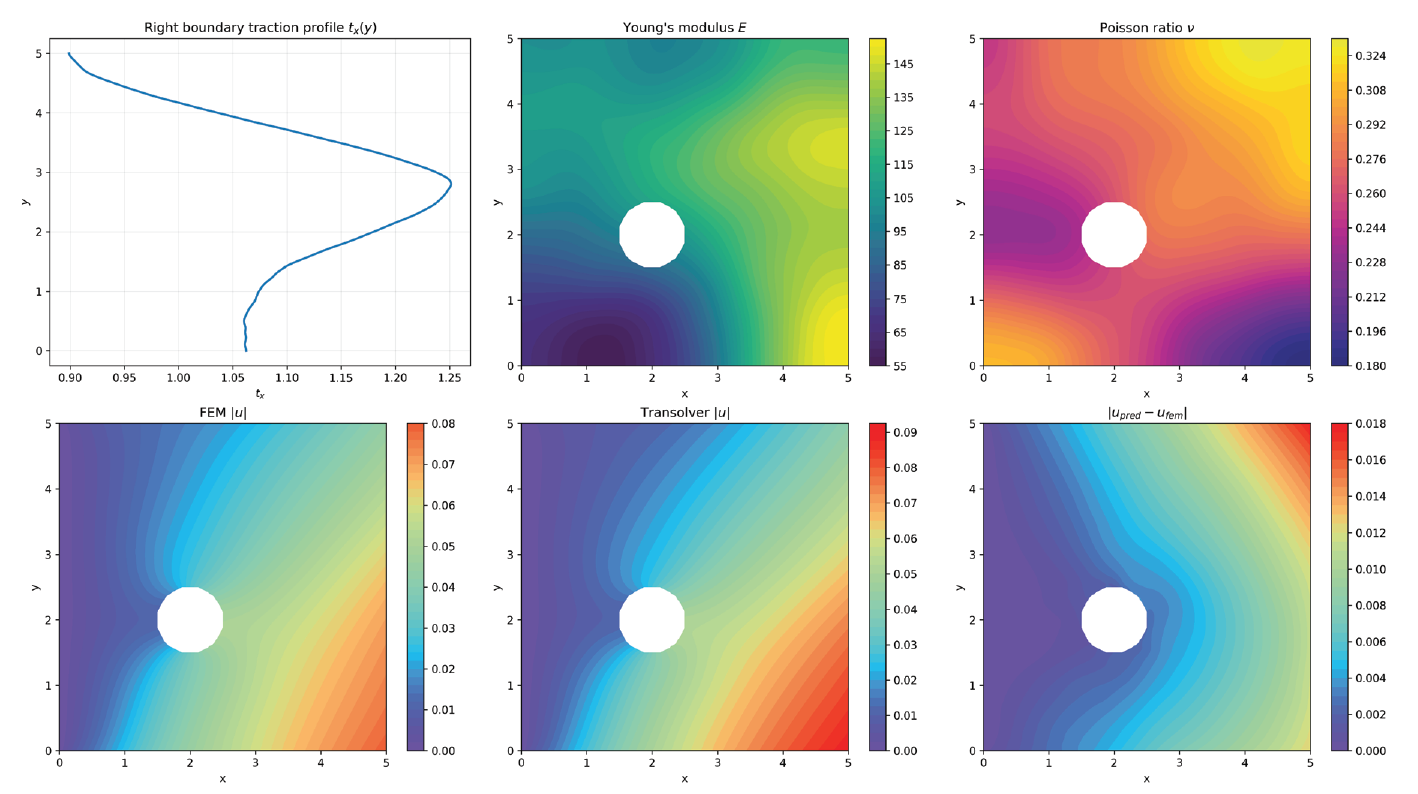} \par\end{centering} \caption{Out-of-Distribution Extrapolation Performance of PFEM in the linear elastic plate problem.\label{fig:PFEM_outdistribution}} \end{figure}

\Cref{fig:PFEM_outdistribution} illustrates the out-of-distribution generalization performance of the PFEM pretraining stage in the linear elastic plate problem. It can be observed that PFEM is able to accurately predict the displacement field, showing good agreement with the FEM reference solution. Compared with the in-distribution cases, a moderate increase in the prediction error is observed, which is expected under distribution shifts. Nevertheless, the overall solution quality remains high, demonstrating the robustness and strong extrapolation capability of the PFEM pretraining stage.

\subsection{Patch tests for PFEM}

In PFEM, the pretraining backbone is currently chosen as Transolver. However, it is expected that more advanced neural operator architectures will emerge in the future, and the pretraining model of PFEM can be naturally replaced. This raises an important practical question: how can one efficiently assess whether a candidate pretraining model possesses sufficient generalization potential within the PFEM framework? A straightforward approach would be extensive validation using large datasets, but such a strategy is computationally expensive and time-consuming.

To address this issue, we propose a patch test for PFEM, which serves as a fast and effective criterion to evaluate the potential of a pretraining model using only a very small amount of data, or even a single sample.

The classical patch test is a fundamental numerical experiment in the finite element method, designed to verify the consistency and convergence properties of finite elements. Its core idea is to prescribe an analytically known displacement field, from which the corresponding strain and stress fields are derived via constitutive relations and equilibrium equations. These fields are then imposed as boundary conditions on a small computational patch consisting of only a few elements. If the element formulation and numerical implementation are correct, the finite element solution should exactly or consistently reproduce the prescribed analytical solution. Due to its simplicity, low computational cost, and strong diagnostic capability, the patch test has become a standard tool in finite element development, particularly as a necessary condition for theoretical convergence.

Analogously, PFEM admits a corresponding patch-test-style numerical experiment during its pretraining stage. The basic procedure of the PFEM patch test is as follows: a single data sample associated with a known analytical solution or a high-fidelity numerical solution is selected and used as the sole training input. During training, PFEM is optimized exclusively through PDE constraints, without relying on any data-driven loss. After training, the PFEM prediction is compared against the reference solution. If the error reaches a sufficiently low level, this indicates that the PFEM model is capable of converging under PDE constraints for the given problem.

Furthermore, the convergence behavior of the error with respect to the number of training iterations is recorded. Suppose convergence is achieved after $N$ iterations for a single sample. When the dataset size is increased to $M$ samples with a batch size of one, the maximum number of training iterations required for convergence is typically bounded by $N \times M$, and in practice often significantly fewer.

It is important to emphasize that the PFEM patch test provides a necessary condition for model selection. If a pretraining model fails to converge even under single-sample training with PDE constraints, it is unlikely to generalize effectively to larger datasets and can therefore be discarded. Since this test involves only minimal data and limited computational cost, it offers an efficient and practical tool for rapidly screening candidate neural operator architectures within the PFEM framework.

\subsection{Probability density function of input point distributions}

The pretraining framework of PFEM is built upon Transolver, whose physics-attention mechanism fundamentally relies on weighted summation, which can be interpreted as a generalized integral operator. As a result, changes in the probability density function (PDF) of the input point distribution may introduce approximation errors.

In practice, we recommend using uniformly distributed sampling points, corresponding to a uniform probability density function, as the default strategy. 
At present, sensitivity to the input-point PDF constitutes a limitation of PFEM and represents a direction for future research and improvement. Nevertheless, our numerical experiments indicate that moderate deviations from uniform sampling in the training point distributions do not lead to significant accuracy degradation, suggesting that PFEM exhibits a reasonable degree of robustness with respect to sampling variations.

\subsection{PDEs as the complete data source}

The training of PFEM relies entirely on governing physical equations and does not require any labeled data, which essentially corresponds to a zero-shot learning paradigm. In contrast, conventional operator learning approaches are typically purely data-driven, requiring large datasets generated in advance by high-fidelity numerical solvers.

For problems in which the governing PDEs are explicitly known, such data-driven preprocessing is often unnecessary and introduces additional computational overhead. In fact, the PDEs themselves constitute the complete source of information, as traditional numerical solvers also obtain solutions solely by solving PDEs. PFEM leverages this observation by embedding the PDEs directly into the operator-learning process.

Beyond eliminating the need for precomputed training data, this strategy offers an additional advantage: it avoids discretization errors inherent in numerically generated datasets. Solutions obtained from conventional solvers are computed on fixed meshes and inevitably contain discretization-induced inaccuracies. Training neural operators directly through PDE constraints circumvents this issue.

A critical component of PDE-based training is the construction of derivatives. In PFEM, we strongly advocate the use of explicit function differentiation, such as derivatives derived from finite element shape functions, rather than automatic differentiation. Explicit function differentiation incurs significantly lower computational cost and achieves efficiency comparable to data-driven training. Moreover, it provides superior numerical robustness by exploiting the inherent structure and stability of finite element shape functions.

\subsection{Unlimited generalization and self-learning capability of PFEM}

PFEM possesses an inherent capability for unlimited generalization. This property arises from the fact that, in the warm-start stage, PFEM employs conventional FEM-based iterative solvers, which naturally exhibit extrapolation capability beyond the training distribution. Consequently, the performance of PFEM does not fundamentally deteriorate when encountering unseen geometries, material properties, or boundary conditions. In practice, the number of iterations required in the warm-start stage primarily depends on the quality of the initial solution provided by the pretraining stage.

To reduce the computational cost of the warm-start stage, it is therefore crucial that the pretraining phase delivers a high-quality initial guess. Notably, the pretraining phase of PFEM does not require labeled data; it only relies on point-cloud inputs encoding geometry, material properties, and boundary conditions. This observation highlights an important future research direction: how to optimally sample the input space during pretraining?

In principle, the distributions of geometry, material properties, and boundary conditions are infinite-dimensional. A practical objective is thus to train PFEM on a finite yet representative subset of this space, such that the learned operator covers the largest possible range of test scenarios. One promising strategy is to generate training inputs from a spectral perspective, for example by progressively enriching the training set from low-frequency to high-frequency components in the Fourier domain. This approach may provide a principled balance between finite training resources in practice and infinite generalization capability in theory.

In summary, the higher the accuracy achieved in the pretraining phase of PFEM, the fewer iterations are required in the warm-start phase. As PFEM is exposed to increasingly diverse geometries, material distributions, and boundary conditions over time, its pretraining model can continuously improve. In this sense, PFEM exhibits a self-learning capability: it becomes progressively more efficient as more scenarios are encountered. This characteristic fundamentally distinguishes PFEM from conventional FEM solvers, which do not possess the ability to improve.

\subsection{Comparison between PFEM warm-start and preconditioned iterative methods}

In the warm-start stage of PFEM, the prediction generated by the Transolver model is used as the initial guess for the iterative solver. While PFEM shares certain similarities with classical preconditioned iterative methods, the underlying mechanisms are fundamentally different.

The common objective of both approaches is to reduce the number of iterations required for convergence. However, traditional preconditioning techniques aim to approximate the inverse of the stiffness matrix in order to improve the conditioning of the linear system, thereby increasing the convergence rate of the iterative solver. In contrast, PFEM as illustrated in \Cref{fig:homo_warm_start} accelerates convergence by providing a more accurate initial guess, which effectively changes the starting point of the iteration rather than the spectral properties of the system matrix.

The relationship between PFEM warm-start and preconditioned iteration can be intuitively understood through a mountain-climbing analogy. Both methods seek to reach the summit, corresponding to the converged solution. A preconditioner reshapes the mountain, reducing steep slopes and making the ascent easier. PFEM, on the other hand, starts the climb from a higher elevation by supplying a better initial guess, thereby shortening the remaining distance to the summit.

\section{Conclusion}

In this work, we propose the Pretrain Finite Element Method (PFEM), a novel computational framework that integrates the efficiency of neural operator learning with the accuracy and robustness of classical finite element methods, while eliminating the need for labeled training data. PFEM consists of two stages: a pretraining stage and a warm-start stage. In the pretraining stage, a Transolver-based neural operator is trained purely through PDEs, using only point-cloud inputs that encode all essential information for numerical simulation, including geometry (spatial coordinates), material properties (e.g., Young’s modulus and Poisson’s ratio), and boundary conditions. In the warm-start stage, the prediction produced by the pretrained operator is used as the initial guess for a conventional iterative solver, thereby accelerating convergence without sacrificing accuracy.
We systematically validate PFEM on both linear and nonlinear problems, including elasticity and hyperelasticity benchmarks, under variations in geometry, material properties, and boundary conditions. The numerical results demonstrate that the pretraining stage of PFEM exhibits strong generalization capability across these variations, while the warm-start stage significantly reduces the number of iterations required by classical solvers. These results confirm that PFEM effectively combines fast operator-based prediction with high-fidelity numerical refinement.

Despite its promising performance, PFEM currently has several limitations and avenues for future extension. First, the present study focuses on quasi-static problems; extending PFEM to transient and time-dependent problems is an important direction for future research. While spatial derivatives are constructed using finite-element-based shape functions, temporal derivatives in transient problems could be incorporated using finite difference schemes. Second, similar to classical FEM, PFEM faces challenges related to the curse of dimensionality in high-dimensional problems. Recent work by Mandl et al.~\citep{mandl2025separable} addresses dimensionality issues in PINO through outer-product-based representations, and such ideas could be leveraged to further extend PFEM to higher-dimensional settings. Third, although the final accuracy of Transolver is comparable to, or even better than, that of FNO, we observe that Transolver may converge more slowly during training. Developing more efficient training strategies for Transolver is therefore an important topic for future work. Finally, since PFEM training relies solely on physical equations rather than data, the design of effective sampling strategies for input geometries, material distributions, and boundary conditions remains a key challenge. In principle, these inputs form an infinite-dimensional space; training on a finite yet representative subset potentially guided by spectral considerations such as progressive enrichment from low-frequency to high-frequency components—may offer a practical balance between finite computational resources and broad generalization.
Currently, we combine Transolver with FEM to utilize mesh-based initial solutions for the loss function’s energy and the initial values at FEM grid nodes. However, Transolver is inherently capable of directly handling point cloud data. In the future, Transolver would probably go well with some mesh-free methods, further expanding its versatility and enabling more flexible PDE solving without being constrained by grid requirements. This shift could open new possibilities for solving complex problems in computational mechanics.

In summary, we believe that PFEM represents the future of the finite element method, in much the same way that the advent of computers once gave rise to FEM itself. Today, artificial intelligence plays a role analogous to that of computers in the past, and we anticipate that it will give rise to a new paradigm in computational mechanics, namely PFEM.

\section*{Declaration of competing interest}
The authors declare that they have no known competing financial interests or personal relationships that could
have appeared to influence the work reported in this paper.

\section*{Acknowledgement}
The study was supported by the Key Project of the National Natural Science Foundation of China (12332005) and scholarship from Bauhaus University in Weimar.  We would like to thank Haixu Wu for the insightful discussions.

\section*{CRediT authorship contribution statement}

\textbf{Yizheng Wang}: Conceptualization, Methodology, Formal analysis, Investigation, Data curation, Validation, Visualization, Writing – original draft, Writing – review \& editing. 
\textbf{Zhongkai Hao}: Investigation. 
\textbf{Mohammad Sadegh Eshaghi}: Investigation.  
\textbf{Cosmin Anitescu}: Supervision, Writing – review \& editing.  
\textbf{Xiaoying Zhuang}: Supervision, Writing – review \& editing.  
\textbf{Timon Rabczuk}: Supervision, Writing – review \& editing, Funding acquisition.  
\textbf{Yinghua Liu}: Supervision, Funding acquisition.

\appendix

\section{Comparison of computational graphs for automatic and explicit differentiation}
\label{sec:computational-graph}

Automatic differentiation (AD) suffers from rapidly increasing computational cost as the order of derivatives increases, whereas the computational cost of explicit differentiation remains essentially independent of the derivative order.

To illustrate this difference, we consider the simple PDE shown in \Cref{eq:poisson_equation}:
\begin{equation}
	\nabla \cdot \big[ k(x,y)\nabla T(x,y) \big] = f(x,y),
	\label{eq:computational_graph}
\end{equation}
where the neural operator aims to learn the operator mapping $k(\boldsymbol{x}) \mapsto T(\boldsymbol{x})$. Spatial derivatives of $T(\boldsymbol{x})$ are then required to evaluate the differential operator.

\begin{figure}[t]
	\centering
	\includegraphics[scale=0.45]{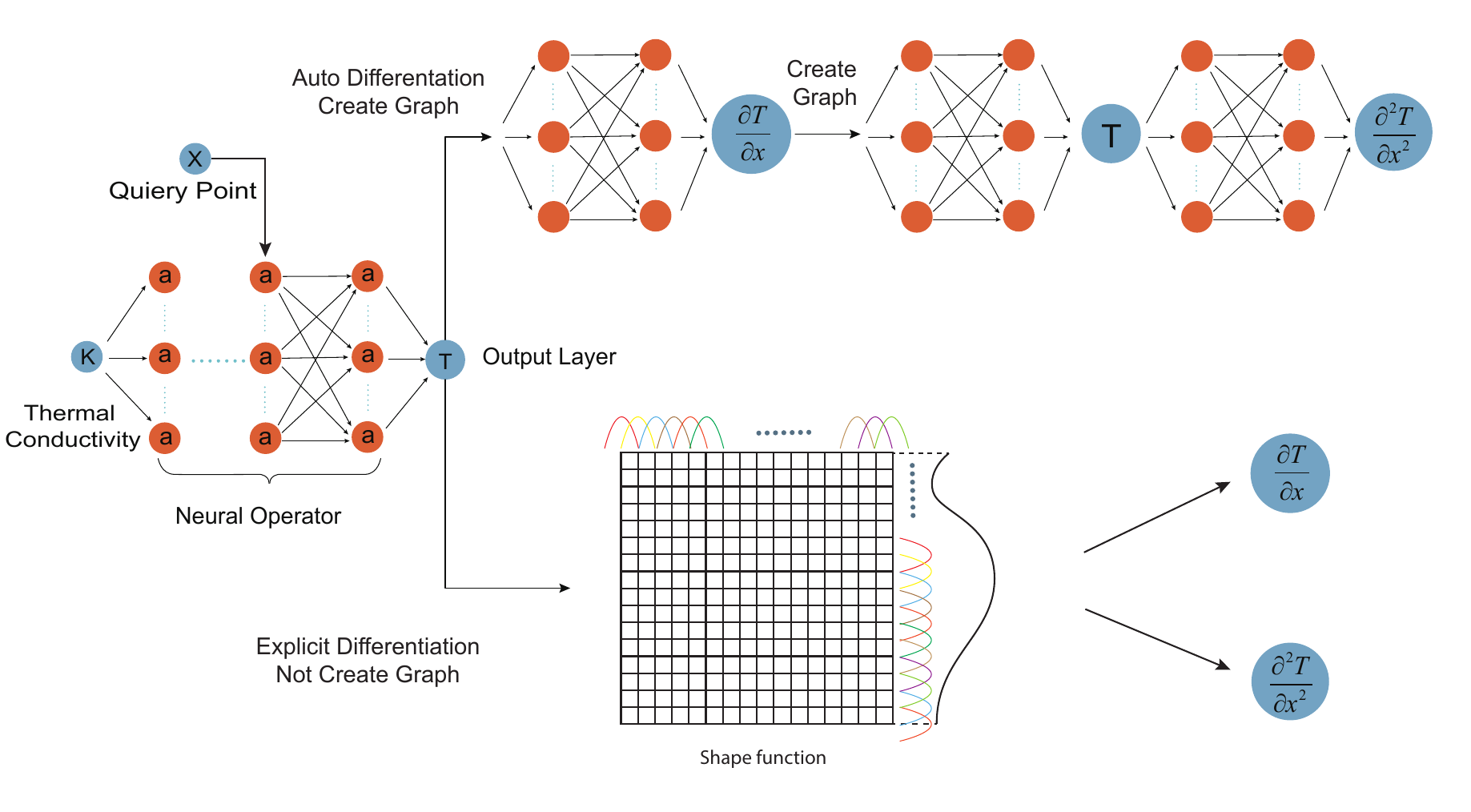}
	\caption{Comparison of computational graphs for automatic differentiation and explicit differentiation.}
	\label{fig:Computation-graph}
\end{figure}

Below, we analyze the computational complexity of automatic differentiation and explicit differentiation, respectively.

The upper part of \Cref{fig:Computation-graph} illustrates the computational graph associated with automatic differentiation of a neural operator. When constructing differential operators via AD, new computational graphs are generated recursively through the chain rule. Let $C_T$ denote the number of basic computational steps required to evaluate the mapping $\boldsymbol{x} \mapsto T(\boldsymbol{x})$ through the neural operator. Then:
\begin{itemize}
	\item The computational graph for $\boldsymbol{x} \mapsto \partial T / \partial \boldsymbol{x}$ requires approximately $2C_T$ steps;
	\item The computational graph for $\boldsymbol{x} \mapsto \partial^2 T / \partial \boldsymbol{x}^2$ requires approximately $4C_T$ steps.
\end{itemize}
In general, if the highest derivative order appearing in the PDE is $M$, the total computational complexity scales as
$
\mathcal{O}(2^M C_T),
$
where $C_T$ represents the base computational cost of the neural network mapping from coordinates to the target field. This exponential growth in complexity with respect to the derivative order makes automatic differentiation increasingly inefficient for PDE-based training involving higher-order derivatives.

In contrast, explicit differentiation reconstructs $T(\boldsymbol{x})$ using finite element shape functions. As a result, spatial derivatives can be derived analytically in closed form, regardless of the derivative order. Once the shape-function representation is established, derivatives of any order can be evaluated explicitly without introducing additional computational graphs.

Consequently, the computational cost of explicit differentiation remains constant with respect to the derivative order and does not grow with $M$. This property fundamentally distinguishes explicit differentiation from automatic differentiation in PDE-based operator learning.

In the pretraining phase of PFEM, differential operators must be evaluated repeatedly. If automatic differentiation is used (e.g., via PyTorch AD \citep{automatic_differential}), additional computational graphs are inevitably constructed, leading to a substantial increase in computational cost and memory usage. By contrast, explicit differentiation achieves both analytical accuracy and high computational efficiency.

Specifically, for an $M$-th order derivative evaluated at $N$ query points, the computational complexity of automatic differentiation and explicit differentiation scales as
$
\mathcal{O}(2^M C_T N)
\quad \text{and} \quad
\mathcal{O}(N),
$
respectively. Here, $N$ denotes the total number of query points. This analysis demonstrates that explicit differentiation is far more suitable for the pretraining phase of PFEM, offering efficiency comparable to purely data-driven approaches while retaining strong physical consistency.

\section{Iterative algorithms of the finite element method}
\label{sec:Iterative-FEM}
We first consider the linear finite element equation
\begin{equation}
	\boldsymbol{K}\boldsymbol{U}=\boldsymbol{F},
	\label{eq:linear_KXB}
\end{equation}
where $\boldsymbol{U}$ denotes the nodal values of the unknown field (e.g., displacement), $\boldsymbol{F}$ is the nodal external force vector, and $\boldsymbol{K}$ is the global stiffness matrix.

There are two major numerical approaches to solve \Cref{eq:linear_KXB}: direct methods and iterative methods. For large-scale problems, iterative methods are typically preferred due to their favorable computational complexity. The iterative scheme associated with \Cref{eq:linear_KXB} can be written as
\begin{equation}
	\boldsymbol{U}^{(k+1)}=\boldsymbol{\phi}(\boldsymbol{U}^{(k)};\boldsymbol{K},\boldsymbol{F}),
\end{equation}
where $\boldsymbol{U}^{(k)}$ denotes the solution at the $k$-th iteration. The iteration is terminated when the relative residual satisfies
\begin{equation}
	\mathcal{L}_{\text{rel}}=\frac{\|\boldsymbol{K}\boldsymbol{U}-\boldsymbol{F}\|}{\|\boldsymbol{F}\|}<\text{tol},
\end{equation}
where $\|\cdot\|$ denotes the Euclidean norm and $\text{tol}$ is a prescribed convergence threshold.

The essential difference between linear and nonlinear finite element methods lies in the fact that nonlinear FEM requires solving a nonlinear algebraic system
\begin{equation}
	\boldsymbol{f}^{\text{int}}(\boldsymbol{u})
	=
	\boldsymbol{f}^{\text{ext}}(\boldsymbol{u})
	+
	\boldsymbol{f}^{\text{kin}}(\boldsymbol{u}),
	\label{eq:nonlinear_FEM}
\end{equation}
where $\boldsymbol{f}^{\text{int}}(\boldsymbol{u})$ is the internal force vector as a nonlinear function of the displacement $\boldsymbol{u}$, $\boldsymbol{f}^{\text{ext}}(\boldsymbol{u})$ is the external force vector, and $\boldsymbol{f}^{\text{kin}}(\boldsymbol{u})$ denotes the inertial force. In linear FEM, $\boldsymbol{f}^{\text{int}}(\boldsymbol{u})=\boldsymbol{K}\boldsymbol{u}$, leading to a linear system. For nonlinear problems, the Newton--Raphson method is commonly employed. In this work, we restrict our discussion to quasi-static processes, for which $\boldsymbol{f}^{\text{kin}}=0$.

Nonlinear FEM involves both the material stiffness matrix $\boldsymbol{K}^{\text{mat}}$ and the geometric stiffness matrix $\boldsymbol{K}^{\text{geo}}$. To derive these terms, we consider the principle of virtual power:
\begin{equation}
	\int_{\Omega}\delta v_{i,j}\sigma_{ij}\,\mathrm{d}V
	=
	\int_{\Omega}\delta v_i f_i\,\mathrm{d}V
	+
	\int_{\Gamma^t}\delta v_i \bar{t}_i\,\mathrm{d}\Gamma,
	\label{eq:principle_of_virtual_power}
\end{equation}
where $v_i$ denotes the velocity field.

The current configuration $\boldsymbol{x}$ and the reference configuration $\boldsymbol{X}$ are discretized using shape functions:
\begin{align}
	x_i(\boldsymbol{X},t) &= N_I(\boldsymbol{X})x_{iI}(t),\\
	X_i &= N_I(\boldsymbol{X})X_{iI},
\end{align}
where Einstein summation is assumed. The displacement and velocity fields are expressed as
\begin{align}
	u_i(\boldsymbol{X},t) &= N_I(\boldsymbol{X})u_{iI}(t),\\
	v_i(\boldsymbol{X},t) &= N_I(\boldsymbol{X})\dot{u}_{iI}(t).
\end{align}

Substituting the above expressions into \Cref{eq:principle_of_virtual_power} yields
\begin{align}
	\delta\dot{u}_{iI}
	\Big(
	\int_{\Omega}\frac{\partial N_I}{\partial x_j}\sigma_{ij}\,\mathrm{d}V
	-
	\int_{\Omega}N_I f_i\,\mathrm{d}V
	-
	\int_{\Gamma^t}N_I\bar{t}_i\,\mathrm{d}\Gamma
	\Big)
	=0.
\end{align}
Since $\delta\dot{u}_{iI}$ is arbitrary, the discrete equilibrium equations follow as
\begin{align}
	f_{iI}^{\text{int}} - f_{iI}^{\text{ext}} &= 0,\\
	f_{iI}^{\text{int}} &= \int_{\Omega}\frac{\partial N_I}{\partial x_j}\sigma_{ij}\,\mathrm{d}V,\\
	f_{iI}^{\text{ext}} &= \int_{\Omega}N_I f_i\,\mathrm{d}V + \int_{\Gamma^t}N_I\bar{t}_i\,\mathrm{d}\Gamma.
\end{align}

Nonlinear FEM is commonly formulated using either the total Lagrangian or the updated Lagrangian approach \citep{belytschko2013nonlinear}. In this work, we adopt the total Lagrangian formulation, which is particularly suitable for hyperelastic problems. The internal and external forces are written as
\begin{align}
	f_{iI}^{\text{int}} &=
	\int_{\Omega_0}
	\frac{\partial N_I}{\partial X_l}
	\frac{\partial x_i}{\partial X_k}
	S_{kl}\,\mathrm{d}V_0,\\
	f_{iI}^{\text{ext}} &=
	\int_{\Omega_0}N_I f_i^0\,\mathrm{d}V_0
	+
	\int_{\Gamma_0^t}N_I\bar{t}_i^0\,\mathrm{d}\Gamma_0,
\end{align}
where $\boldsymbol{S}$ is the second Piola--Kirchhoff stress tensor, related to the Cauchy stress by $\boldsymbol{\sigma}=J^{-1}\boldsymbol{F}\boldsymbol{S}\boldsymbol{F}^T$.

Taking the incremental form leads to
\begin{align}
	\dot{\boldsymbol{f}}^{\text{int}} &=
	(\boldsymbol{K}^{\text{geo}}+\boldsymbol{K}^{\text{mat}})\dot{\boldsymbol{U}},\\
	\dot{\boldsymbol{f}}^{\text{ext}} &=
	\boldsymbol{K}^{\text{ext}}\dot{\boldsymbol{U}},
\end{align}
with
\begin{align}
	\boldsymbol{K}^{\text{geo}} &= \int_{\Omega_0}\boldsymbol{G}_0^T\boldsymbol{S}\boldsymbol{G}_0\,\mathrm{d}V_0,\\
	\boldsymbol{K}^{\text{mat}} &= \int_{\Omega_0}\boldsymbol{B}_0^T\boldsymbol{C}^{SE}\boldsymbol{B}_0\,\mathrm{d}V_0.
\end{align}

For nonlinear problems, the residual at the $(n+1)$-th load step is defined as
\begin{equation}
	\boldsymbol{r}(\boldsymbol{U}_{k+1}^{(n+1)})
	=
	\boldsymbol{f}^{\text{int}}(\boldsymbol{U}_{k+1}^{(n+1)})
	-
	\boldsymbol{f}^{\text{ext}}(\boldsymbol{U}_{k+1}^{(n+1)}).
\end{equation}
Applying a Taylor expansion at $\boldsymbol{U}_k^{(n+1)}$ yields
\begin{align}
	\boldsymbol{r}(\boldsymbol{U}_{k+1}^{(n+1)})
	=
	\boldsymbol{r}(\boldsymbol{U}_k^{(n+1)})
	+
	\frac{\partial\boldsymbol{r}}{\partial\boldsymbol{U}}
	\Big|_{\boldsymbol{U}_k^{(n+1)}}
	\Delta\boldsymbol{U}.
\end{align}
Setting the residual to zero leads to
\begin{align}
	\Delta\boldsymbol{U}
	&=
	-
	\left(\boldsymbol{K}^{\text{geo}}+\boldsymbol{K}^{\text{mat}}-\boldsymbol{K}^{\text{ext}}\right)^{-1}
	\boldsymbol{r}(\boldsymbol{U}_k^{(n+1)}),\\
	\boldsymbol{U}_{k+1}^{(n+1)}
	&=
	\boldsymbol{U}_k^{(n+1)}+\Delta\boldsymbol{U}.
\end{align}
The iteration terminates when $\|\boldsymbol{r}\|\leq \text{tol}$.

Nonlinear FEM involves two nested iterations: the outer iteration over load steps and the inner Newton--Raphson iteration at each load step. In PFEM, the outer iteration over the first $N_f-1$ load steps is bypassed. Instead, the pretrained neural operator provides an initial guess $\boldsymbol{U}^{NO}$ directly at the final load step, i.e.,
\begin{equation}
	\boldsymbol{U}_0^{(N_f)}=\boldsymbol{U}^{NO}.
\end{equation}
The warm-start phase of PFEM then follows the iterative scheme
\begin{equation}
	\boldsymbol{U}_{k+1}^{(N_f)}
	=
	\boldsymbol{\phi}
	\big(
	\boldsymbol{U}_k^{(N_f)};
	\boldsymbol{K}^{\text{geo}},
	\boldsymbol{K}^{\text{mat}},
	\boldsymbol{K}^{\text{ext}},
	\boldsymbol{U}_0^{(N_f)}=\boldsymbol{U}^{NO}
	\big).
\end{equation}
Since the pretrained PFEM solution provides an initial guess close to the true solution at the final load step, the number of Newton iterations required in nonlinear FEM is significantly reduced.

\section{Benefit of the warm-start stage in PFEM \label{sec:PFEM_warmstart_benefit}}

Compared with conventional numerical solvers, the warm-start stage of PFEM can significantly reduce the number of iterations required for convergence. Below, we provide a convergence analysis for iterative methods applied to linear systems to explain this phenomenon.

We consider the fixed-point iteration
\begin{equation}
	\boldsymbol{X}^{(k)}=\boldsymbol{B}\boldsymbol{X}^{(k-1)}+\boldsymbol{f},
	\label{eq:iterative_way}
\end{equation}
where the exact solution $\boldsymbol{X}^{*}$ satisfies
\begin{equation}
	\boldsymbol{X}^{*}=\boldsymbol{B}\boldsymbol{X}^{*}+\boldsymbol{f}.
	\label{eq:exact_solution}
\end{equation}
We define the iteration error as
\begin{equation}
	\boldsymbol{e}^{(k)}=\boldsymbol{X}^{(k)}-\boldsymbol{X}^{*}.
	\label{eq:iter_error}
\end{equation}
Substituting \Cref{eq:iterative_way,eq:exact_solution} into \Cref{eq:iter_error} yields
\begin{equation}
	\begin{aligned}
		\boldsymbol{e}^{(k)}
		&=(\boldsymbol{B}\boldsymbol{X}^{(k-1)}+\boldsymbol{f})-(\boldsymbol{B}\boldsymbol{X}^{*}+\boldsymbol{f})\\
		&=\boldsymbol{B}(\boldsymbol{X}^{(k-1)}-\boldsymbol{X}^{*})
		=\boldsymbol{B}\boldsymbol{e}^{(k-1)}.
	\end{aligned}
\end{equation}
By recursion, we obtain
\begin{equation}
	\boldsymbol{e}^{(k)}=\boldsymbol{B}^{k}\boldsymbol{e}^{(0)}.
\end{equation}
Note that the quality of the initial guess is fully characterized by $\boldsymbol{e}^{(0)}=\boldsymbol{X}^{(0)}-\boldsymbol{X}^{*}$. A better initial guess implies a smaller error norm $\|\boldsymbol{e}^{(0)}\|$. Using the definition of the induced (subordinate) matrix norm, we have
\begin{equation}
	\|\boldsymbol{B}^{k}\|
	=\max_{\boldsymbol{e}^{(0)}}\frac{\|\boldsymbol{B}^{k}\boldsymbol{e}^{(0)}\|}{\|\boldsymbol{e}^{(0)}\|}
	\geq
	\frac{\|\boldsymbol{B}^{k}\boldsymbol{e}^{(0)}\|}{\|\boldsymbol{e}^{(0)}\|},
\end{equation}
and thus
\begin{equation}
	\|\boldsymbol{e}^{(k)}\|
	=\|\boldsymbol{B}^{k}\boldsymbol{e}^{(0)}\|
	\leq
	\|\boldsymbol{B}^{k}\|\;\|\boldsymbol{e}^{(0)}\|.
	\label{eq:iter_error_evolution}
\end{equation}
This inequality indicates that, to reach the same target accuracy $\|\boldsymbol{e}^{(k)}\|$, a better initial guess reduces the upper bound of the required iteration number $k$. This provides a theoretical explanation for why a more accurate initial solution from the pretraining stage of PFEM can reduce the iteration count in the warm-start stage, especially in scenarios where conventional iterative solvers cannot easily obtain a good initial guess.

The benefit brought by a good initial guess becomes less pronounced when the convergence tolerance (\texttt{tol}) is set to be very small. Rearranging \Cref{eq:iter_error_evolution} gives
\begin{equation}
	\frac{\|\boldsymbol{e}^{(k)}\|}{\|\boldsymbol{e}^{(0)}\|}
	\leq
	\|\boldsymbol{B}^{k}\|.
	\label{eq:ratio_error}
\end{equation}
Let $\varepsilon=\|\boldsymbol{B}^{k}\|$, then \Cref{eq:ratio_error} can be rewritten as
\begin{equation}
	\frac{\|\boldsymbol{e}^{(k)}\|}{\|\boldsymbol{e}^{(0)}\|}
	\leq
	\varepsilon.
	\label{eq:elison}
\end{equation}
From $\|\boldsymbol{B}^{k}\|=\varepsilon$, we have
\begin{align}
	\|\boldsymbol{B}^{k}\|^{\frac{1}{k}} &= \varepsilon^{\frac{1}{k}},\nonumber\\
	\ln\|\boldsymbol{B}^{k}\|^{\frac{1}{k}} &= \ln\varepsilon^{\frac{1}{k}},\nonumber\\
	k &= \frac{-\ln\varepsilon}{-\ln\|\boldsymbol{B}^{k}\|^{\frac{1}{k}}}.
	\label{eq:iter_num_discussion}
\end{align}
We define the asymptotic convergence rate as
\begin{equation}
	R(\boldsymbol{B})=\lim_{k\rightarrow+\infty}-\ln\|\boldsymbol{B}^{k}\|^{\frac{1}{k}}.
	\label{eq:liter_ratio}
\end{equation}
The rate $R(\boldsymbol{B})$ depends only on the iteration matrix $\boldsymbol{B}$ and is independent of the iteration index $k$ and the initial guess. For sufficiently large $k$, substituting \Cref{eq:liter_ratio} into \Cref{eq:iter_num_discussion} yields
\begin{equation}
	k=\frac{-\ln\varepsilon}{R(\boldsymbol{B})}.
	\label{eq:K_elison}
\end{equation}

Different choices of the convergence tolerance essentially correspond to different target magnitudes of $\|\boldsymbol{e}^{(k)}\|$. Let $\texttt{tol}=10^{-m}$. In practice, one may regard $\|\boldsymbol{e}^{(k)}\|$ as being on the order of \texttt{tol} at convergence. Then \Cref{eq:elison} becomes
\begin{align}
	\frac{10^{-m}}{\|\boldsymbol{e}^{(0)}\|} &\leq \varepsilon,\nonumber\\
	\ln\frac{10^{-m}}{\|\boldsymbol{e}^{(0)}\|} &\leq \ln\varepsilon,\nonumber\\
	-m\ln 10-\ln\|\boldsymbol{e}^{(0)}\| &\leq \ln\varepsilon.
	\label{eq:elison_tol}
\end{align}
Substituting \Cref{eq:elison_tol} into \Cref{eq:K_elison} yields
\begin{align}
	k &= \frac{-\ln\varepsilon}{R(\boldsymbol{B})},\nonumber\\
	&\leq
	\frac{m\ln 10+\ln\|\boldsymbol{e}^{(0)}\|}{R(\boldsymbol{B})}
	=
	\frac{m\ln 10}{R(\boldsymbol{B})}
	+
	\frac{\ln\|\boldsymbol{e}^{(0)}\|}{R(\boldsymbol{B})}.
	\label{eq:k_R}
\end{align}

From \Cref{eq:k_R}, the term $m\ln 10/R(\boldsymbol{B})$ characterizes the influence of the convergence tolerance on the iteration number, whereas $\ln\|\boldsymbol{e}^{(0)}\|/R(\boldsymbol{B})$ quantifies the influence of the initial guess. Therefore, as the tolerance becomes tighter (i.e., $m$ increases), the relative contribution of the initial guess to the overall iteration count becomes smaller, implying a diminishing speedup benefit from improved initialization when extremely small tolerances are required.

\section{Numerical homogenization theory \label{sec:homo}}

In this section, we briefly introduce the numerical homogenization framework proposed by Andreassen et al.~\citep{andreassen2014determine}, which aims to compute the effective macroscopic elastic tensor of periodic composite materials. For simplicity, we consider a two-phase material composed of solid and void regions. The elastic tensor of the solid phase is denoted by $E_{ijkl}^{S}$.

To compute the homogenized elastic tensor $E_{IJKL}^{H}$, six independent unit strain fields are applied to the unit cell, including three normal strains and three shear strains. The six macroscopic unit strain tensors $\varepsilon_{pq}^{(0)(IJ)}$ are defined as

\begin{equation}
	\begin{aligned}
		\boldsymbol{\varepsilon}^{(0)(11)} &=
		\begin{bmatrix}
			1 & 0 & 0\\
			0 & 0 & 0\\
			0 & 0 & 0
		\end{bmatrix},
		\quad
		\boldsymbol{\varepsilon}^{(0)(22)} =
		\begin{bmatrix}
			0 & 0 & 0\\
			0 & 1 & 0\\
			0 & 0 & 0
		\end{bmatrix},
		\quad
		\boldsymbol{\varepsilon}^{(0)(33)} =
		\begin{bmatrix}
			0 & 0 & 0\\
			0 & 0 & 0\\
			0 & 0 & 1
		\end{bmatrix},
		\\
		\boldsymbol{\varepsilon}^{(0)(12)} &=
		\begin{bmatrix}
			0 & 1 & 0\\
			1 & 0 & 0\\
			0 & 0 & 0
		\end{bmatrix},
		\quad
		\boldsymbol{\varepsilon}^{(0)(13)} =
		\begin{bmatrix}
			0 & 0 & 1\\
			0 & 0 & 0\\
			1 & 0 & 0
		\end{bmatrix},
		\quad
		\boldsymbol{\varepsilon}^{(0)(23)} =
		\begin{bmatrix}
			0 & 0 & 0\\
			0 & 0 & 1\\
			0 & 1 & 0
		\end{bmatrix}.
	\end{aligned}
	\label{eq:unit strain}
\end{equation}

\begin{figure}
	\begin{centering}
		\includegraphics[scale=0.45]{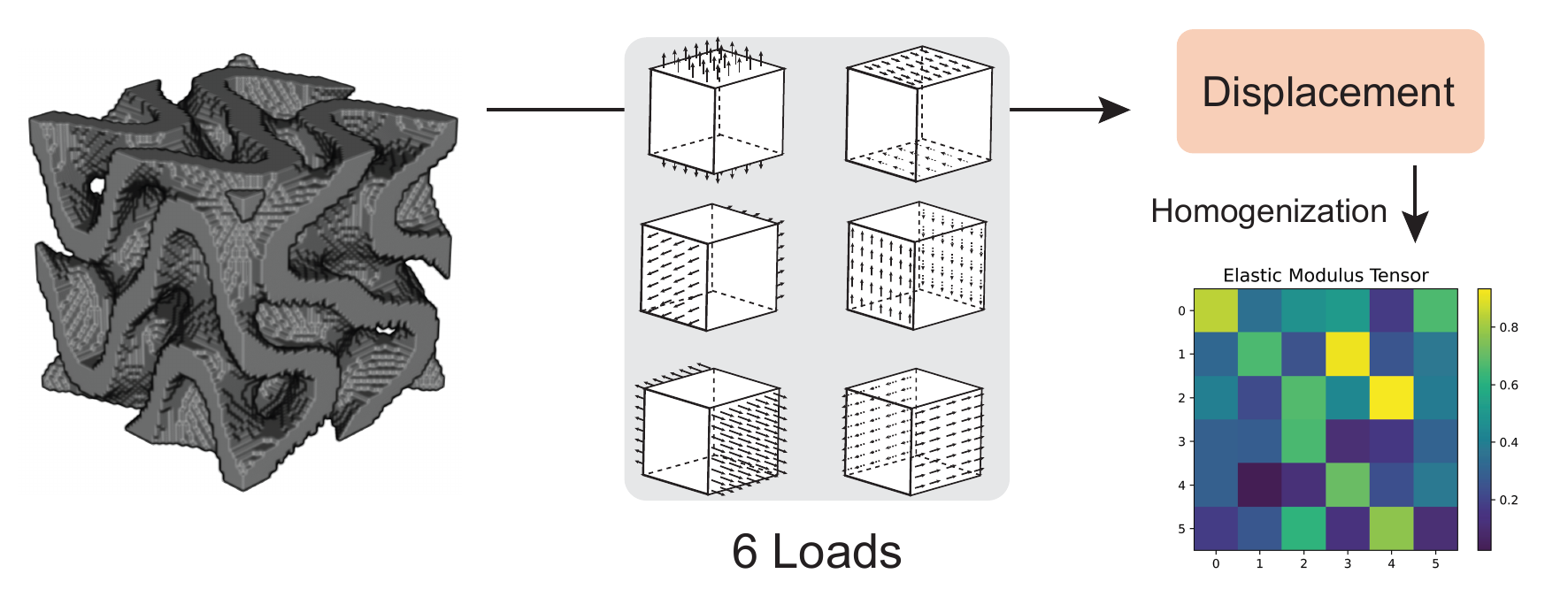}
		\par\end{centering}
	\caption{Schematic illustration of numerical homogenization. 
		Six macroscopic unit strain fields are applied to the unit cell, producing the corresponding displacement fields. The homogenized elastic tensor is then computed through the homogenization formula based on the resulting displacement fields.}
	\label{fig:Homogenization_schematic}
\end{figure}

Here $\varepsilon_{pq}^{(0)(IJ)}$ denotes a macroscopic unit strain field imposed on the unit cell domain $\Omega$. The capital indices $(IJ)$ correspond to the six loading cases and do not follow the Einstein summation convention. The six loading configurations are illustrated in \Cref{fig:Homogenization_schematic}.

The homogenized elastic tensor $E_{IJKL}^{H}$ is computed as

\begin{equation}
	E_{IJKL}^{H}
	=
	\frac{1}{|\Omega|}
	\int_{\Omega}
	E_{pqrs}
	\left(
	\varepsilon_{pq}^{(0)(IJ)} - \varepsilon_{pq}^{(v)(IJ)}
	\right)
	\left(
	\varepsilon_{rs}^{(0)(KL)} - \varepsilon_{rs}^{(v)(KL)}
	\right)
	\, d\Omega ,
	\label{eq:homogenized elastic matrix}
\end{equation}
where $E_{pqrs}$ denotes the local elastic tensor of the TPMS unit cell,

\begin{equation}
	E_{ijkl} =
	\begin{cases}
		E_{ijkl}^{S}, & \text{in the solid phase}, \\
		0, & \text{in the void phase}.
	\end{cases}
\end{equation}
where $|\Omega|$ denotes the volume of the unit cell. The fluctuation strain field $\varepsilon_{pq}^{(v)(IJ)}$ is defined as

\begin{equation}
	\varepsilon_{pq}^{(v)(IJ)}
	=
	\varepsilon_{pq}\big(\boldsymbol{X}^{(v)(IJ)}\big)
	=
	\frac{1}{2}
	\left(
	X_{p,q}^{(v)(IJ)} + X_{q,p}^{(v)(IJ)}
	\right),
\end{equation}
where $\boldsymbol{X}^{(v)(IJ)}$ denotes the fluctuation displacement field induced by the imposed macroscopic unit strain field. The total displacement field can therefore be written as

\begin{equation}
\boldsymbol{X}^{(IJ)} =
\boldsymbol{X}^{(0)(IJ)} - \boldsymbol{X}^{(v)(IJ)},
\end{equation}
where $\boldsymbol{X}^{(0)(IJ)}$ represents the macroscopic affine displacement field corresponding to the uniform strain field, while $\boldsymbol{X}^{(v)(IJ)}$ represents the fluctuation displacement caused by material heterogeneity.

We consider the total potential energy of the unit cell

\begin{equation}
	\begin{aligned}
		\Pi^{(IJ)}
		&=
		\int_{\Omega}
		\frac{1}{2}
		E_{klpq}
		\varepsilon_{kl}^{(IJ)}
		\varepsilon_{pq}^{(IJ)}
		\, d\Omega ,
		\\
		\varepsilon_{kl}^{(IJ)}
		&=
		\frac{1}{2}
		\left(
		X_{k,l}^{(IJ)} + X_{l,k}^{(IJ)}
		\right)
		=
		\varepsilon_{kl}^{(0)(IJ)} -
		\varepsilon_{kl}^{(v)(IJ)} .
	\end{aligned}
	\label{eq:homo_energy}
\end{equation}

Taking the first variation of $\Pi^{(IJ)}$ and using the symmetry $E_{klpq}=E_{pqkl}$ yields

\begin{equation}
	\delta \Pi^{(IJ)}
	=
	\int_{\Omega}
	E_{klpq}
	\varepsilon_{kl}^{(IJ)}
	\delta \varepsilon_{pq}^{(IJ)}
	\, d\Omega
	=
	-
	\int_{\Omega}
	E_{klpq}
	\big(
	\varepsilon_{kl}^{(0)(IJ)} -
	\varepsilon_{kl}^{(v)(IJ)}
	\big)
	\delta \varepsilon_{pq}^{(v)(IJ)}
	\, d\Omega .
\end{equation}

By enforcing $\delta \Pi^{(IJ)} = 0$, the fluctuation displacement field satisfies the following weak form of the equilibrium equation:

\begin{equation}
	\int_{\Omega}
	E_{klpq}
	\varepsilon_{kl}(\boldsymbol{X}^{(v)(IJ)*})
	\varepsilon_{pq}(\boldsymbol{X}^{(v)(IJ)})
	\, d\Omega
	=
	\int_{\Omega}
	E_{klpq}
	\varepsilon_{kl}(\boldsymbol{X}^{(v)(IJ)*})
	\varepsilon_{pq}^{(0)(IJ)}
	\, d\Omega ,
	\label{eq:weak_form}
\end{equation}
where $\boldsymbol{X}^{(v)(IJ)*}$ denotes a virtual fluctuation displacement field. In practice, the fluctuation displacement field $\boldsymbol{X}^{(v)(IJ)}$ is typically solved using the finite element method. The fluctuation strain $\boldsymbol{\varepsilon}^{(v)(IJ)}$ can then be computed and substituted into \Cref{eq:homogenized elastic matrix} to obtain the homogenized elastic tensor.

To ensure a unique solution and correctly represent the periodic microstructure, the weak form \eqref{eq:weak_form} is solved subject to periodic boundary conditions (PBCs). These conditions enforce identical fluctuation displacements on opposite faces of the unit cell. Assuming $\Omega = [0,L]^3$, the PBCs can be expressed as

\begin{equation}
	\begin{aligned}
		\boldsymbol{X}^{(v)(IJ)}(0,x_2,x_3)
		&=
		\boldsymbol{X}^{(v)(IJ)}(L,x_2,x_3),
		\\
		\boldsymbol{X}^{(v)(IJ)}(x_1,0,x_3)
		&=
		\boldsymbol{X}^{(v)(IJ)}(x_1,L,x_3),
		\\
		\boldsymbol{X}^{(v)(IJ)}(x_1,x_2,0)
		&=
		\boldsymbol{X}^{(v)(IJ)}(x_1,x_2,L).
	\end{aligned}
	\label{eq:PBCs}
\end{equation}

In conventional finite element implementations, the PBCs are enforced by constraining the degrees of freedom of corresponding nodes on opposite faces of the unit cell.

In PFEM, periodic constraints are imposed through a ring-type coupling strategy. Specifically, pairs of corresponding nodes on opposite faces are connected by virtual elements, and the associated virtual strain energy is incorporated into the overall loss function. It should be emphasized that this virtual strain energy formulation is essentially equivalent to a penalty method, where the stiffness matrix of the virtual elements acts as the penalty parameter controlling the strength of the periodic constraint.

The fluctuation displacement field $\boldsymbol{X}^{(v)(IJ)}$ must be solved for all six loading cases during the homogenization process. In traditional approaches, this requires repeated finite element simulations, which can be computationally expensive, especially when different geometries and material distributions are considered. Therefore, accelerating the computation of displacement fields directly improves the overall efficiency of the homogenization process. This observation motivates the use of PFEM to accelerate numerical homogenization.

\section{Generation of TPMS geometries\label{sec:TPMS}}

In this section, we briefly introduce the mathematical description of triply periodic minimal surfaces (TPMS). A minimal surface is characterized by zero mean curvature at every point. TPMS are periodic minimal surfaces that extend infinitely in three-dimensional space. A common way to represent TPMS is through a level-set formulation. In this framework, the geometry is described by a scalar field combined with trigonometric functions, which enables flexible control of both topology and mechanical properties. The level-set equation is generally written as
\begin{equation}
	\phi(x,y,z)=c,
\end{equation}
where $x$, $y$, and $z$ denote the spatial coordinates. The scalar field $\phi(x,y,z)$ defines a continuous field in the three-dimensional domain, and the level-set value $c$ partitions a TPMS unit cell into two connected subdomains, corresponding to the solid phase and the void phase, respectively. By varying $c$, the porosity, or equivalently the volume fraction, of the TPMS unit cell can be adjusted. The TPMS surface is thus represented by the iso-surface associated with the level-set value $c$.

One widely used TPMS level-set model is the Schwarz Primitive surface, which is defined as
\begin{equation}
	\cos(2\pi x)+\cos(2\pi y)+\cos(2\pi z)=c.
\end{equation}
This level-set equation represents a zero-thickness shell-like minimal surface. To construct load-bearing TPMS lattice materials from such minimal surfaces, two common types of TPMS architectures are considered.

The first type is referred to as \emph{solid networks}, defined by
\begin{equation}
	\phi(x,y,z)>c.
\end{equation}
For solid-network TPMS cells, the region satisfying $\phi(x,y,z)>c$ is identified as the solid phase, while the remaining region is treated as void.

The second type is referred to as \emph{sheet networks}, defined by
\begin{equation}
	-c\leq \phi(x,y,z)\leq c.
\end{equation}
As the name suggests, sheet-network TPMS selects a thin volumetric layer around the zero-thickness minimal surface, thereby producing a sheet-like topology within the unit cell. Solid-network and sheet-network TPMS exhibit distinct topological characteristics and effective elastic properties, providing greater flexibility for topology-driven material design.

Another important topological parameter of TPMS is the volume fraction, also known as the relative density. The volume fraction is defined as the ratio of the solid volume to the total volume of the TPMS unit cell. This parameter directly influences the effective elastic tensor and stiffness of the unit cell, and therefore plays a key role in the homogenized elastic response.

In this study, we generated a series of TPMS unit cells to train and validate the proposed neural-operator-based homogenization framework. Three representative TPMS families were considered, namely the \emph{Schoen Gyroid}, \emph{Schwarz Diamond}, and \emph{Fischer Koch S}. Their corresponding level-set equations are given by
\begin{equation}
	\begin{aligned}
		&\text{Schoen Gyroid:}\quad
		\sin(2\pi x)\cos(2\pi y)+\sin(2\pi y)\cos(2\pi z)
		+\sin(2\pi z)\cos(2\pi x)=c,\\
		&\text{Schwarz Diamond:}\quad
		\cos(2\pi x)\cos(2\pi y)\cos(2\pi z)
		-\sin(2\pi x)\sin(2\pi y)\sin(2\pi z)=c,\\
		&\text{Fischer Koch S:}\quad
		\cos(4\pi x)\sin(2\pi y)\cos(2\pi z)
		+\cos(2\pi x)\cos(4\pi y)\sin(2\pi z)\\
		&\qquad\qquad\qquad\quad
		+\sin(2\pi x)\cos(2\pi y)\cos(4\pi z)=c.
	\end{aligned}
	\label{eq:TPMStype}
\end{equation}

\section{Supplementary code}
The code of this work will be available at \url{https://github.com/yizheng-wang/Research-on-Solving-Partial-Differential-Equations-of-Solid-Mechanics-Based-on-PINN} after accepted.

\bibliographystyle{elsarticle-num}
\addcontentsline{toc}{section}{\refname}\bibliography{reference.bib}

\end{document}